\documentclass{article}
\usepackage{enumitem}

\usepackage{amsthm}
\usepackage{bm}
\usepackage[a4paper,margin=2.5cm,nohead]{geometry}

\usepackage{natbib}
\usepackage{amsmath,amssymb}
\usepackage{graphicx}
\usepackage{algorithm,algorithmic}
\usepackage{textcomp}
\usepackage{colortbl}
\usepackage[breakable]{tcolorbox}
\tcbset{boxrule=0.2mm}
\usepackage{tabularray}

\usepackage{amsmath}
\usepackage{amssymb}

\usepackage{mathtools}

\usepackage{pifont}

\newtheorem{theorem}{Theorem}[section]

\newtheorem{lemma}[theorem]{Lemma}
\newtheorem{corollary}[theorem]{Corollary}
\newtheorem{assumption}[theorem]{Assumption}

\usepackage[utf8]{inputenc} 
\usepackage[T1]{fontenc}    
\usepackage{hyperref}       
\hypersetup{
    colorlinks=true,
    linkcolor=blue,
    filecolor=magenta,      
    urlcolor=cyan,
    citecolor=blue,
    }

\usepackage{url}            
\usepackage{booktabs}       
\usepackage{amsfonts}       
\usepackage{nicefrac}       
\usepackage{microtype}      
\usepackage{xcolor}         

\title{A Stochastic Approximation Approach for Efficient Decentralized Optimization on Random Networks}

\author{%
  Chung-Yiu Yau \\
  The Chinese University of Hong Kong\\
  \texttt{cyyau@se.cuhk.edu.hk}
  \and
  Haoming Liu \\
  Peking University\\
  \texttt{haomingliu@pku.edu.cn}
  \and
  Hoi-To Wai \\
  The Chinese University of Hong Kong\\
  \texttt{htwai@se.cuhk.edu.hk}
  \thanks{This work was supported by the Shun Hing Institute of Advanced Engineering, The Chinese University of Hong Kong, under Project \#MMT-p5-23.}
}

\def\prm{\mathbf{x}}
\def\avgprm{\bar{\mathbf{x}}}
\def\dprm{\bm{\lambda}}

\def\bC{\mathbf{C}}
\def\wta{{\bf A}}

\def\wtq{{\bf Q}}
\def\wtk{{\bf K}}
\def\sampL{L_{\sf s}}

\def\sigavg{{\bar{\sigma}^2}}
\def\data{{\mathbb{P}}}
\def\bbC{{\mathbb{C}}}
\def\cmthread{\texttt{communication\_thread()}}
\def\cpthread{\texttt{computation\_thread()}}
\def\oneotimes{{\bf 1}_{\otimes}^\top}
\def\oneotimesT{{\bf 1}_{\otimes}}
\def\ARA{{\bf L}}

\def\serr{{\bf e}_s}

\def\gerr{{\bf e}_g}

\def\rxi{{\bf R}}
\def\rxii{{\bf R}}
\def\rlamb{ }

\def\a{{\tt a}}
\def\b{{\tt b}}
\def\c{{\tt c}}
\def\d{{\tt d}}
\def\e{{\tt e}}
\def\f{{\tt f}}
\def\g{{\tt g}}
\def\mom{{\bf m}_x}
\def\momi{{\bf m}_{x,i}}
\def\momdual{{\bf m}_\lambda}
\def\momidual{{\bf m}_{\lambda,i}}
\def\avgmom{\overline{\bf m}_x}

\def\algname{{\tt FSPDA}}
\def\algnamesa{{\tt FSPDA-SA}}
\def\algnamevr{{\tt FSPDA-STORM}}

\newcommand{\dotp}[2]{\left\langle{#1}\ \middle|\ {#2}\right\rangle}
\newcommand{\expec}[1]{\mathbb{E}\left[ {#1} \right]}

\newcommand*{\oscarhide}[1]{}

\newlength\myindent
\begin{document}

\maketitle

\begin{abstract}
A challenging problem in decentralized optimization is to develop algorithms with fast convergence on random and time varying topologies under unreliable and bandwidth-constrained communication network. This paper studies a stochastic approximation approach with a \underline{F}ully \underline{S}tochastic \underline{P}rimal \underline{D}ual \underline{A}lgorithm ({\tt FSPDA}) framework. Our framework relies on a novel observation that randomness in time varying topology can be incorporated in a stochastic augmented Lagrangian formulation, whose expected value admits saddle points that coincide with stationary solutions of the decentralized optimization problem. With the {\tt FSPDA} framework, we develop two new algorithms supporting efficient sparsified communication on random time varying topologies --- {\tt FSPDA-SA} allows agents to execute multiple local gradient steps depending on the time varying topology to accelerate convergence, and {\tt FSPDA-STORM} further incorporates a variance reduction step to improve sample complexity. For problems with smooth (possibly non-convex) objective function, within $T$ iterations, we show that {\tt FSPDA-SA} (resp.~{\tt FSPDA-STORM}) finds an $\mathcal{O}( 1/\sqrt{T} )$-stationary (resp.~$\mathcal{O}( 1/T^{2/3} )$) solution. Numerical experiments show the benefits of the {\tt FSPDA} algorithms.
\end{abstract}

\section{Introduction}
Consider $n$ agents that communicate on an undirected and connected graph/network $\mathcal{G} = (\mathcal{V}, \mathcal{E})$ with ${\cal V} = [n] := \{1, \ldots, n\}, {\cal E} \subseteq {\cal V} \times {\cal V}$. Each agent $i \in [n]$ has access to a continuously differentiable (possibly non-convex) local objective function $f_i: \mathbb{R}^d \rightarrow \mathbb{R}$ and maintains a local decision variable $\prm_i \in \mathbb{R}^d$. Denote $\prm = [\prm_1^\top, ..., \prm_n^\top ]^\top \in \mathbb{R}^{nd}$. Our aim is to tackle:
\begin{equation}
\textstyle \min_{\prm \in \mathbb{R}^{nd}} ~\frac{1}{n} \sum_{i=1}^n f_i(\prm_i) \quad {\rm s.t.}  \quad \prm_i = \prm_j,~\forall~(i,j) \in {\cal E}. 
\label{eq:main_problem}
\end{equation}
In other words, \eqref{eq:main_problem} seeks a $\prm^\star \in \mathbb{R}^d$ that minimizes $F( \prm ) := (1/n) \sum_{i=1}^n f_i( \prm )$. We are interested in the stochastic optimization setting where each $f_i( \prm_i )$ is given by (with slight abuse of notation)
\begin{equation} \label{eq:stoc-fi}
f_i( \prm_i ) := \mathbb{E}_{ \xi_i \sim \data_i } [ f_i ( \prm_i ; \xi_i ) ]
\end{equation}
where $\data_i$ represents the $i$-th data distribution. Problem~\eqref{eq:main_problem} is relevant to the distributed learning problem especially in the decentralized case where a central server is absent. Prior works \citep{nedic2009distributed, lian2017can,nedic2017achieving,qu2017harnessing} demonstrated that {\it decentralized} algorithms can tackle \eqref{eq:main_problem} efficiently through repeated message exchanges among the neighbors and local stochastic gradient updates.
 
Towards an efficient decentralized algorithm for \eqref{eq:main_problem}, an important direction is to consider a {\it time varying graph topology} setting where the {\it active edge set} in ${\cal G}$ changes over time. This is a generic setting covering cases when the communication links are unreliable, or the agents choose not to communicate in a certain round (a.k.a.~local updates) \citep{koloskova2019decentralized2,nadiradze2021asynchronous}. 
By assuming that a random topology is drawn at each iteration, the convergence of decentralized stochastic gradient (DSGD) has been studied in \citep{lobel2010distributed, nadiradze2021asynchronous} and is later on unified by \citep{koloskova2020unified} with tighter bounds for local updates, periodic sampling, etc. An alternative \citep{sundhar2010distributed} is to analyze DSGD for the $B$-connectivity setting which requires the union of every $B$ consecutive time varying topologies to yield a connected graph. Nevertheless, these works focused on vanilla DSGD that may have slow convergence (in transient stage) and is limited to bounded data heterogeneity. The prior restrictions can be relaxed using advanced algorithms such as gradient tracking \citep{qu2017harnessing}, EXTRA \citep{shi2015extra} and primal-dual framework \citep{hong2017prox,hajinezhad2019perturbed,yi2021linear}.

\begin{table*}[t]
    \def\arraystretch{1.1}
    \centering
    \begin{tabular}{l|c c c c c}
    \toprule 
    Prior Works & {\tt SG} & {\tt TV} & {\tt w/o\,BH} & {\tt Rate} \\
    \midrule
    {Prox-GPDA} {\citep{hong2017prox}} & \ding{55} & \ding{55} & $\checkmark$ & {\tt Asympt.} \\
    {NEXT} {\citep{di2016next}} & \ding{55} & $\checkmark$ & $\checkmark$ & {\tt Asympt.} \\
    {DSGD} \citep{koloskova2020unified} & $\checkmark$ & $\checkmark$ & \ding{55} & $\mathcal{O} ({\sigma} / {\sqrt{nT}} )$  \\
    {Swarm-SGD} {\citep{nadiradze2021asynchronous}} & $\checkmark$ & $\checkmark$ & \ding{55} & $\mathcal{O} ( {\sigma^2} / {\sqrt{T}} )$  \\
    {CHOCO-SGD} {\citep{koloskova2019decentralized2}} & $\checkmark$ & \ding{55}$^\ddagger$ & \ding{55} & $\mathcal{O} ({ {\sigma} } / {\sqrt{nT}} )$ \\
    {Decen-Scaffnew} {\citep{mishchenko2022proxskip}} & $\checkmark$ & \ding{55}$^\dagger$ & $\checkmark$ & $\mathcal{O} ( {\sigma} / {\sqrt{nT}} )$ \\
    Local-GT \citep{liu2024decentralized} & $\checkmark$ & \ding{55}$^\dagger$ & $\checkmark$ & $\mathcal{O} ( { {\sigma}} / {\sqrt{nT}} )$ \\ 
    LED \citep{alghunaim2024local} & $\checkmark$ & \ding{55}$^\dagger$ & $\checkmark$ & $\mathcal{O} ( { {\sigma}} / {\sqrt{nT}} )$ \\ 
    \algnamesa ~({\bf This Work}) & $\checkmark$ & $\checkmark$ & $\checkmark$ & $\mathcal{O} ( { {\sigma}} / {\sqrt{nT}} )$ \\
    \algnamevr ~({\bf This Work}) & $\checkmark$ & $\checkmark$ & $\checkmark$ & $\mathcal{O} ( { \sigma^{2/3}} / {T^{2/3}} )$ \\
    \bottomrule
    \end{tabular}
    \caption{Comparison of decentralized algorithms for {\bf non-convex} optimization. 
    In the table, `{\tt SG}' is `Stochastic Gradient', `{\tt TV}' is `Time Varying Graph', `{\tt w/o\,BH}' is `Without Bounded Heterogeneity', and `{\tt Rate}' is the expected squared gradient norm $\mathbb{E}[ \| \nabla F( \avgprm ) \|^2 ]$ after $T$ iterations. Note that $\sigma^2$ is the variance of stochastic gradient. 
    $^\ddagger$CHOCO-SGD incorporates broadcast gossip as a special case of compression.  
    $^\dagger$ProxSkip, Local-GT, LED consider local updates with periodic communication. 
    }
    \label{tab:rate_compare_ncvx}
\end{table*}

As noted by \citep{koloskova2021improved}, analyzing the convergence of sophisticated algorithms with time varying topology, such as gradient tracking \citep{qu2017harnessing} is challenging due to the non-symmetric product of two (or more) mixing matrices. Existing works considered various restrictions on the time varying topology ${\cal G}^{(t)} = ( {\cal V}, {\cal E}^{(t)} )$ {and/or} the problem \eqref{eq:main_problem}: \citep{koloskova2021improved,liu2024decentralized} studied gradient tracking with local updates that essentially takes ${\cal E}^{(t)} = {\cal E}$ periodically and ${\cal E}^{(t)} = \emptyset$ otherwise, also see \citep{mishchenko2022proxskip, guo2023revisiting, alghunaim2024local} for a similar result and note that such algorithms require extra synchronization overhead; \citep{kovalev2021lower, kovalev2024lower} considered a setting where ${\cal G}^{(t)}$ is connected for any $t$; \citep{nedic2017achieving, li2024accelerated} focused on (accelerated) gradient tracking with deterministic gradient when $F(\prm)$ is (strongly) convex; \citep{di2016next} also considered deterministic gradient with possibly non-convex $F(\prm)$ but only provides asymptotic convergence guarantees; \citep{lei2018asymptotic, yau2023fully} considered asymptotic convergence guarantees in the case of strictly (or strongly) convex $F(\prm)$.
We provide a non-exhaustive list summarizing the convergence of existing works in Table~\ref{tab:rate_compare_ncvx}.

The above discussion highlights a gap in the existing literature ---
\begin{center}\vspace{-.1cm}
\emph{Is there any algorithm that achieves fast convergence on time varying (random) topology?}
\vspace{-.1cm}
\end{center}
This paper gives an affirmative answer through developing the \underline{F}ully \underline{S}tochastic \underline{P}rimal \underline{D}ual \underline{A}lgorithm ({\algname}) framework that leads to efficient decentralized algorithms tackling \eqref{eq:main_problem} in its general form. The framework features the design of a new stochastic augmented Lagrangian function. 

As pointed out by \citep{chang2020distributed}, many decentralized algorithms (including gradient tracking) can be interpreted as primal-dual algorithms finding a saddle point of the augmented Lagrangian function. However, its extension to time varying topology is not straightforward due to the inconsistency in dual variables updates. 
To overcome this challenge, we propose a stochastic equality constrained reformulation of \eqref{eq:main_problem} to model randomness in topology. Then, the latter yields a stochastic augmented Lagrangian function. Applying stochastic approximation (SA) to solve the latter leads to the {\algname} framework. Our contributions are 
\begin{itemize}[leftmargin=*, topsep=0pt, itemsep=0pt]
    \item We propose two new algorithms: (i) {\algnamesa} is derived by vanilla SA that applies primal-dual stochastic gradient descent-ascent on the stochastic augmented Lagrangian, (ii) {\algnamevr} uses an additional control variate / momentum term to reduce the drift term's variance in a recursive manner. Both algorithms are fully stochastic as the random time varying topology is treated as a part of randomness. Additionally, our framework supports sparsified communication, i.e., the agents can choose to communicate a subset of primal coordinates at each iteration.
    \item We show that after $T$ iterations, {\algnamesa} (resp.~{\algnamevr}) finds in expectation a solution whose squared gradient norm is ${\cal O}(1/\sqrt{T})$ (resp.~${\cal O}(1/T^{2/3})$). The convergence analysis is derived from a new Lyapunov function design that involves an unsigned inner product term and incorporates a variance condition on the random time varying topologies. 
    Interestingly, we show empirically that 
    using momentum in dual updates benefits the consensus error convergence.
    \item We also demonstrate that both {\algnamesa} and {\algnamevr} can be implemented in a fully asynchronous manner, i.e., the agents can communicate and compute at different time slots, and supports local update as the algorithms allow for arbitrary time varying topology. That said, we remark that the convergence rates with local updates of {\algnamesa} and {\algnamevr} are only suboptimal. 
\end{itemize}
We provide numerical experiments to show that {\algnamesa} and {\algnamevr} outperform existing algorithms in terms of iteration and communication complexity.

\noindent \textbf{Notations.} Let ${\bf W} \in \mathbb{R}^{d \times d}$ be a symmetric (not necessarily positive semidefinite) matrix, the ${\bf W}$-weighted (semi) inner product of vectors ${\bf a,b} \in \mathbb{R}^d$ is denoted as $\dotp{{\bf a}}{{\bf b}}_{{\bf W}} := {\bf a}^\top {\bf W} {\bf b}$. Similarly, the ${\bf W}$-weighted (semi)  norm is denoted by $\| {\bf a} \|_{{\bf W}}^2 := \dotp{{\bf a}}{{\bf a}}_{{\bf W}}$. The subscript notation is omitted for ${\bf I}$-weighted inner products. For any square matrix ${\bf X}$, $( {\bf X} )^\dagger$ denotes its pseudo inverse.

\section{The Fully Stochastic Primal Dual Algorithm ({\algname}) Framework} \label{sec:fspda}\vspace{-.1cm}
This section develops the {\algname} framework for tackling \eqref{eq:main_problem} and describes two variants of the framework leading to decentralized stochastic optimization of \eqref{eq:main_problem}. 
Let $\widetilde{\wta} \in \{-1, 0, 1 \}^{| {\cal E} | \times n}$ be an incidence matrix of ${\mathcal{G}}$. By defining $\wta = \widetilde{\wta} \otimes {\bf I}_d \in \{-1,0,1\}^{|{\cal E}|d \times nd}$, we observe that the consensus constraint in \eqref{eq:main_problem} is equivalent to $\wta \prm = {\bf 0}$.

Our first step is to model the randomness in the time varying topology using the random variable (r.v.) $\xi_a \sim \mathbb{P}_a$. For each realization $\xi_a$, we define the random incidence matrix $\wta(\xi_a) := {\bf I}( \xi_a ) \wta \in \{-1,0,1\}^{|{\cal E}|d \times nd}$ where ${\bf I}(\xi_a) \in \{0,1\}^{|{\cal E}|d \times |{\cal E}|d}$ is a binary diagonal matrix. In addition to selecting each edge of ${\cal G}$ randomly, ${\bf I}(\xi_a)$ selects a random subset of $d$ coordinates. As we will see later, this allows our approach to simultaneously achieve random sparsification for communication compression. 

Assume that $\mathbb{E}_{\xi_a \sim \mathbb{P}_a} [ {\bf I}(\xi_a) ]$ is a positive diagonal matrix, \eqref{eq:main_problem} is equivalent to:
\begin{equation} \label{eq:main_problem2}
\textstyle \min_{\prm \in \mathbb{R}^{nd}} ~\frac{1}{n} \sum_{i=1}^n \mathbb{E}_{ \xi_i \sim \mathbb{P}_i } [ f_i(\prm_i; \xi_i ) ] \quad {\rm s.t.} \quad \mathbb{E}_{\xi_a \sim \mathbb{P}_a }[ \wta(\xi_a) ] \, \prm = {\bf 0}.
\end{equation} 
Denote $\xi = ( \xi_1, \ldots, \xi_n, \xi_a )$, {\algname} hinges on the following \emph{augmented Lagrangian} function of \eqref{eq:main_problem2}:
\begin{equation} \label{eq:stoc-lag} 
\begin{split}
& {\cal L}( \prm, \dprm) := \mathbb{E}_{\xi}[ {\cal L}( \prm, \dprm ; \xi) ] \\
& \text{with} ~~ \textstyle {\cal L}( \prm, \dprm; \xi ) := \sum_{i=1}^n f_i(\prm_i; \xi_i) + \tilde\eta \dotp{\dprm}{\wta(\xi_a) \prm} + \frac{\tilde\gamma}{2}  \| \wta(\xi_a) \prm \|^2 ,
\end{split}
\end{equation}
where $\tilde{\eta} > 0, \tilde{\gamma} > 0$ are penalty parameters.
It can be verified that the saddle points of ${\cal L}( \prm, \dprm)$ correspond to the KKT points of \eqref{eq:main_problem} \citep{bertsekas2016nonlinear}. For brevity, in the rest of this paper, we may drop the subscript in $\xi$ whenever the notation is clear from the context.

{\algname} is developed from applying stochastic approximation (SA) to seek a saddle point of \eqref{eq:stoc-lag}. By recognizing $\wta( \xi )^\top \wta( \xi ) = \wta^\top\wta( \xi )$, we consider the stochastic gradients: 
\begin{equation} \label{eq:stoc-grad}
    \begin{split}
        & \nabla_{\prm} {\cal L}( \prm, \dprm; \xi ) := \nabla {\bf f}( \prm ; \xi ) + \tilde\eta \wta^\top \rlamb \dprm + \tilde\gamma \wta^\top \wta( \xi ) \prm , ~~ \nabla_{\dprm} {\cal L}( \prm, \dprm; \xi ) := \tilde\eta \wta( \xi ) \prm ,
    \end{split}
\end{equation}
where $\nabla {\bf f}( \prm ; \xi ) = [ \nabla f_1( \prm_1; \xi_1 );  \ldots; \nabla f_n( \prm_n; \xi_n ) ] \in \mathbb{R}^{nd}$. Notice that to facilitate algorithm development, we have taken a deterministic $\wta$ for the term in $\nabla_{\prm} {\cal L}$ related to $\dprm$. Now observe the $i$th $d$-dimensional block of $\wta^\top \wta( \xi ) \prm$ which can be aggregated within ${\cal N}_i(\xi)$ the neighborhood of the $i$th agent as:
\begin{equation} \label{eq:local-compute}
\textstyle \big[ \wta^\top \wta( \xi ) \prm \big]_i = \sum_{j \in {\cal N}_i(\xi)} {\bf C}_{ij}(\xi) ( \prm_j - \prm_i ) ,
\end{equation}
where ${\bf C}_{ij}(\xi) \in \{0,1\}^{d \times d}$ is diagonal and depends on the selected coordinates for the edge $(i,j)$ under randomness $\xi$. 
Eq. \eqref{eq:local-compute} {\it only} relies on $\prm_j$ from neighbor $j$ that is connected on the time varying topology ${\cal G}(\xi)$. For illustration, an example of the above random graph model is given by Figure \ref{fig:extended_graph} in Appendix \ref{app:async}. Importantly, \eqref{eq:stoc-grad} shows that with the stochastic augmented Lagrangian function, the time varying topology can be treated implicitly as a part of the randomness in the stochastic primal-dual gradients. The framework is thus described as being \emph{fully stochastic} as in \citep{bianchi2021fully}, and departs from \citep{liu2024decentralized,alghunaim2024local} that treat the topology as fixed during the derivation of primal-dual algorithm(s).
From \eqref{eq:stoc-grad}, \eqref{eq:local-compute}, we derive {\it two} variants of {\algname}. 

\vspace{.1cm} {\bf {\algnamesa} Algorithm.} The first variant of {\algname} is derived from a direct application of stochastic gradient descent-ascent (SGDA) updates. Take $\alpha>0, \beta>0$ as the step sizes, we have
\begin{equation}
    \prm^{t+1} = \prm^t - \alpha \nabla_{\prm} {\cal L}( \prm^t, \dprm^t; \xi^t ) , ~~ \dprm^{t+1} = \dprm^t + \beta \nabla_{\dprm} {\cal L}( \prm^t, \dprm^t; \xi^t ) .
\end{equation}
Taking the variable substitution $\widehat{\dprm} := \wta^\top \dprm$ yields the following recursion:
\begin{tcolorbox}[boxsep=2pt,left=4pt,right=4pt,top=3pt,bottom=3pt]
\textbf{\sf \underline{{\algnamesa}:}} for any $t \geq 0$ and any $i \in [n]$,
\begin{subequations} \label{eq:fspda-sa}
\begin{align}
    \prm_i^{t+1} & = \prm_i^t - \alpha \nabla f_i( \prm_i^t; \xi_i^t ) - \eta \widehat{\dprm}_i^t  \textstyle + \gamma \sum_{j \in {\cal N}_i(\xi_a^t)} {\bf C}_{ij}(\xi_a^t) ( \prm_j^t - \prm_i^t ), \label{eq:fspdasa_primal} \\
    \label{eq:fspdasa_dual}
    \widehat{\dprm}_i^{t+1} & \textstyle = \widehat{\dprm}_i^t + \beta \sum_{j \in {\cal N}_i(\xi_a^t)} {\bf C}_{ij}(\xi_a^t) ( \prm_j^t - \prm_i^t ).
\end{align}
\end{subequations}
\end{tcolorbox}
Note that $\prm^0, \widehat{\dprm}^0$ can be {initialized arbitrarily}. 

\vspace{.1cm} {\bf {\algnamevr} Algorithm.} The second variant of {\algname} reduces the variance of the stochastic gradient term in \eqref{eq:stoc-grad} using the recursive momentum variance reduction technique \citep{cutkosky2019momentum}. Herein, the key idea is to utilize a control variate in estimating the (primal-dual) gradients of ${\cal L}(\prm, \dprm)$. Take $\alpha,\beta>0$ and $a_x, a_{\lambda} \in [0,1]$ as the momentum parameters, we have $\prm^{t+1} = \prm^t - \alpha \mom^t , \dprm^{t+1} = \dprm^t + \beta \momdual^t$ as the primal-dual updates, and 
\begin{equation} \label{eq:fspdavr_mom_main}
    \begin{split}
    & \mom^{t+1} = \nabla_{\prm} {\cal L}(\prm^{t+1}, \dprm^{t+1}; \xi^{t+1}) + (1 - a_x) (\mom^t - \nabla_{\prm}  {\cal L}(\prm^t, \dprm^{t}; \xi^{t+1}) ), \\
    & \momdual^{t+1} = \nabla_{\dprm} {\cal L}(\prm^{t+1}, \dprm^{t+1}; \xi^{t+1}) + (1 - a_\lambda ) (\momdual^t - \nabla_{\dprm} {\cal L}(\prm^t, \dprm^{t}; \xi^{t+1}) ).
    \end{split}
\end{equation}
The aim of $\mom^{t+1}$ is to estimate $\nabla_{\prm} {\cal L}( \prm^{t+1}, \dprm^{t+1} )$. Now, instead of the straightforward estimator $\nabla_{\prm} {\cal L}(\prm^{t+1}, \dprm^{t+1}; \xi^{t+1})$, we include an extra zero-mean term $\mom^t - \nabla_{\prm}  {\cal L}(\prm^t, \dprm^{t}; \xi^{t+1})$ to reduce the variance of the stochastic gradient estimation. The latter is a control variate that is computed recursively. Particularly, it has been shown in \citep{cutkosky2019momentum} that it can effectively reduce variance with a carefully designed parameter $a_x$, provided that the stochastic gradient map satisfies a mean-square Lipschitz condition. 
We summarize the algorithm as follows. 
\begin{tcolorbox}[boxsep=2pt,left=4pt,right=4pt,top=3pt,bottom=3pt]
\textbf{\sf \underline{{\algnamevr}:}} for any $t \geq 0$ and any $i \in [n]$,
\begin{subequations} \label{eq:fspda-vr}
\begin{align}
    \prm_i^{t+1} & = \prm_i^t - \alpha {\momi^t}, \label{eq:fspdavr_primal} \\
    \label{eq:fspdavr_dual}
    \widehat{\dprm}_i^{t+1} & \textstyle = \widehat{\dprm}_i^t + \beta {\momidual^t}, \\
    {\momi^{t+1}} & = (1 \! - \! a_x) \big[ {\momi^t} + \nabla f_i( \prm_i^t; \xi_i^{t+1} ) - \eta \widehat{\dprm}_i^t  \textstyle + \gamma \sum_{j \in {\cal N}_i(\xi_a^{t+1})} {\bf C}_{ij}(\xi_a^{t+1}) ( \prm_j^t - \prm_i^t ) \big] \label{eq:fspdavr_mom} \\
    & \quad + \nabla f_i( \prm_i^{t+1}; \xi_i^{t+1} ) - \eta \widehat{\dprm}_i^{t+1}  \textstyle + \gamma \sum_{j \in {\cal N}_i(\xi_a^{t+1})} {\bf C}_{ij}(\xi_a^{t+1}) ( \prm_j^{t+1} - \prm_i^{t+1} ) \notag \\
    \momidual^{t+1} & = (1-a_\lambda) \big[ \momidual^t + \textstyle \sum_{j \in {\cal N}_i(\xi_a^{t+1})} {\bf C}_{ij}(\xi_a^{t+1}) ( \prm_j^t - \prm_i^t ) \big] \label{eq:fspdavr_momdual} \\
    & \quad + \textstyle \sum_{j \in {\cal N}_i(\xi_a^{t+1})} {\bf C}_{ij}(\xi_a^{t+1}) ( \prm_j^{t+1} - \prm_i^{t+1} ) \notag
\end{align}
\end{subequations}
\end{tcolorbox}
Note that to achieve the theoretical performance (see later in Sec.~\ref{sec:convergence}), $\prm^0, \widehat{\dprm}^0, \mom^0, \momdual^0$ shall be initialized as $\prm_i^0 = \avgprm^0$, $\widehat{\dprm}_i^{0} = (\alpha/\eta) n^{-1} (\nabla F(\avgprm^0) - \nabla f_i(\avgprm^0))$, $\momi^0 = \nabla F(\avgprm^0) $, $\momidual^0 = {\bf 0}$ according to \eqref{eq:storm_grad_bound}. We remark that a simple initialization choice $\widehat{\dprm}^{0} = \momi^{0} = \momidual^{0} = {\bf 0}$ works well in practice.

Both {\algnamesa} and {\algnamevr} are decentralized algorithms that can be implemented on random time varying topology, and support randomized sparisification for further communication compression. The key is to observe that in \eqref{eq:fspda-sa}, \eqref{eq:fspda-vr}, the only information required for agent $i$ is to obtain $\sum_{j \in {\cal N}_i(\xi_a^t)} {\bf C}_{ij}(\xi_a^t) ( \prm_j^t - \prm_i^t )$, and in addition $\sum_{j \in {\cal N}_i(\xi_a^t)} {\bf C}_{ij}(\xi_a^t) ( \prm_j^{t-1} - \prm_i^{t-1} )$ for {\algnamevr}, at iteration $t$. 

\subsection{Implementation Details and Connection to Existing Works}

We discuss several features of the {\algname} algorithms and their connections to existing works.

\vspace{.1cm} {\bf Local \& Asynchronous Updates.} 
The \emph{local update} scheme where each agent $i$ is allowed to update its own local variables $\prm_i, \dprm_i$ for multiple iterations without a communication step is a common practice in decentralized optimization \citep{liu2024decentralized,li2024accelerated,alghunaim2024local,mishchenko2022proxskip}. 
As discussed before, such scheme can be seen as a special case of the {\algname} framework where the time varying topology ${\cal E}^{(t)}$ is chosen such that the latter alternates between ${\cal E}^{(t)} = {\cal E}$ and ${\cal E}^{(t)} = \emptyset$. 

Furthermore, {\algnamesa} allows for the general case of {\it asynchronous} updates. This is done so by taking the stochastic gradient as $\nabla f_i( \prm_i^t; \xi^t ) = b_i( \xi^t ) \, \overline{b}_i \, \nabla f_i( \prm_i^t; {\xi}^t )$ such that $b_i (\xi^t) \in \{0,1\}$ with $\mathbb{E}[ b_i (\xi^t) ] = 1/ \overline{b}_i$ for some constant $\overline{b}_i > 0$. Detailed discussions for a fully asynchronous implementation of {\algnamesa} can be found in Appendix \ref{app:async}.

\vspace{.1cm} {\bf Connection to Existing Works.}
Evaluating $\prm^{t+2} - \prm^{t+1}$ from the {\algnamesa} sequence and observe that the combination of \eqref{eq:fspdasa_primal} and \eqref{eq:fspdasa_dual} is equivalent to the second order recursion:
\begin{equation} \label{eq:fspda_single}
\begin{aligned}
    \prm^{t+2} &= 2\left({\bf I} - \frac{\gamma}{2} \wta^\top \wta(\xi^{t+1}) \right) \prm^{t+1}  - \left( {\bf I} - (\gamma - \eta \beta) \wta^\top \rlamb \wta(\xi^t)  \right) \prm^t \\
    & \quad - \alpha \left( \nabla {\bf f}(\prm^{t+1}; \xi^{t+1}) - \nabla {\bf f}(\prm^t; \xi^t) \right).
\end{aligned}
\end{equation}
This reduces the {\algnamesa} recursion into a primal-only sequence by eliminating the dual sequence $\dprm^t$. 
In the deterministic optimization setting when $\wta(\xi) \equiv \wta$ and $\nabla {\bf f}(\prm; \xi) \equiv \nabla {\bf f}( \prm )$,  \eqref{eq:fspda_single} is equivalent to the EXTRA algorithm \citep{shi2015extra} using the mixing matrix ${\bf W} = {\bf I} - \gamma {\rm Diag}( \tilde{\bf W} {\bf 1} ) + \gamma \tilde{\bf W}$ where $\tilde{\bf W}$ is the 0-1 adjacency matrix of ${\cal G}$. Here, with an appropriate choice of $\gamma$, ${\bf W}$ will be doubly stochastic and satisfies the convergence requirement in \citep{shi2015extra}.
Similar observations have been made in \citep{nedic2017achieving} for the gradient tracking and DIGing algorithms.

On the other hand, for stochastic optimization on random networks, \eqref{eq:fspda_single} suggests each agent to keep the current and previous iterates received from neighbors in the corresponding time varying topology. In this case, \eqref{eq:fspda_single} yields an extension of the EXTRA/GT algorithms to time varying topology. 

\section{Convergence Analysis of {\algname}} \label{sec:convergence}
This section presents the convergence rate analysis of {\algname} for \eqref{eq:main_problem}. Unless otherwise specified, we focus on the case with smooth but possibly non-convex objective function. Specifically, we consider:
\begin{assumption} \label{assm:lip}
    Each $f_i$ is $L$-smooth, i.e., for $i = 1, \ldots, n$,
    \begin{equation} \label{eq:f_lip}
        \|\nabla f_i({\bf x}) - \nabla f_i({\bf y}) \| \leq L \|{\bf x} - {\bf y} \| ~\forall~ {\bf x},{\bf y} \in \mathbb{R}^{d}.
    \end{equation}
    There exists $f_\star > -\infty$ such that $f_i( {\bf x} ) \geq f_\star$ for any ${\bf x} \in \mathbb{R}^d$.
\end{assumption}
\noindent Note this implies that the global objective function $F(\cdot)$ is $L$-smooth but possibly non-convex. 

We further assume that the random network ${\cal G} (\xi_a) $ is connected in expectation, yet each realization ${\cal G}(\xi_a)$ may not be connected. Let $\rxii = \expec{ {\bf I}(\xi_a) }$, this leads to the following property concerning the expected graph Laplacian matrix $\wta^\top \rxii \wta = \expec{ \wta(\xi_a)^\top \wta }$. Defining the matrix $\wtk := ( {\bf I}_n - {\bf 1}{\bf 1}^\top/n ) \otimes {\bf I}_d$, we have
\begin{assumption} \label{assm:rand-graph}
    There exists $\rho_{\max} \geq \rho_{\min} > 0$ and $\bar{\rho}_{\max} \geq \bar{\rho}_{\min} > 0$ such that 
    \begin{equation}
    \begin{aligned} 
    & \rho_{\min} \wtk \preceq \wta^\top \rxii \wta \preceq \rho_{\max} \wtk \quad  \text{and} \quad \bar{\rho}_{\min} \wtk \preceq \wta^\top \wta \preceq \bar{\rho}_{\max} \wtk.
    \end{aligned}
    \end{equation}
\end{assumption}
\noindent It holds that $\wta^\top \rxii \wta \wtk = \wta^\top \rxii \wta = \wtk \wta^\top \rxii \wta$. The above assumption can be satisfied if ${\cal G}$ is connected \citep{yi2021linear}, \cite[Lemma 2]{yi2018distributed} and ${\rm diag}(\rxii) > {\bf 0}$ such that each edge is selected with a positive probability. As an important consequence, if $\gamma \leq \rho_{\min} / \rho_{\max}^2$, we have 
\[
\| ( {\bf I} - \gamma \wta^\top \rxii \wta ) \prm \|_{\wtk}^2 \leq (1 - \gamma \rho_{\min} ) \| \prm \|_{\wtk}^2,~\forall~\prm \in \mathbb{R}^{nd}.
\]
We thus observe that the operator $({\bf I} - \gamma \wta^\top \rxii \wta)$ serves a similar purpose as the mixing matrix in a average consensus algorithms and $\rho_{\min}$ can be interpreted as the spectral radius of ${\cal G}$ similar to \citep[Eq.~(12)]{koloskova2020unified}. 
Moreover, if we define  
$\wtq := ( \wta^\top \rxii \wta )^\dagger$
such that it holds $\wtq \wta^\top \rxii \wta = \wta^\top \rxii \wta \wtq = \wtk$, Assumption~\ref{assm:rand-graph} implies that $\rho_{\max}^{-1} \wtk \preceq \wtq \preceq \rho_{\min}^{-1} \wtk$. 

Next we consider several assumptions on the noise variance of the random quantities in {\algname}:
\begin{assumption} \label{assm:f_var} For any fixed $\prm_i \in \mathbb{R}^d$, $i \in [n]$, there exists $\sigma_i \geq 0$ such that
    \begin{equation} \label{eq:f_var}
        \mathbb{E}_{ \xi_i \sim \data_i} [\| \nabla f_i(  \prm_i ; \xi_i ) - \nabla f_i( \prm_i ) \|^2] \leq \sigma_i^2.
    \end{equation}
    To simplify notations, we define $\sigavg := (1/n) \sum_{i=1}^n \sigma_i^2$.
\end{assumption}
\begin{assumption} \label{assm:graph_var}
    For any fixed $\prm \in \mathbb{R}^{nd}$, there exists $\sigma_A \geq 0$ such that
    \begin{equation}
        \mathbb{E}_{\xi_a \sim \data_a} [ \|\wta(\xi_a)^\top \wta \prm - \wta^\top \rlamb \rxii \wta \prm \|^2 ] \leq \sigma_A^2 \| \prm \|^2_{\wtk}. \label{eq:graph_var}
    \end{equation}
\end{assumption}
\noindent Assumption \ref{assm:f_var} is standard. Meanwhile for Assumption \ref{assm:graph_var}, the variance term $\sigma_A^2$ measures the quality of the random topology ${\cal G}(\xi_a)$ in approximating the expected graph Laplacian $\wta^\top \rxii \wta$. The latter is important as it contributes to the variance in the drift term of {\algname}. Observe that $\sigma_A^2$ decreases with the proportion of edges selected in each random subgraph ${\cal G}( \xi_a )$.  

To facilitate our discussions, we define the following quanitites:
\begin{equation} \textstyle
\textstyle \avgprm^t := \frac{1}{n} \sum_{i=1}^n \prm_i^t, \quad \sum_{i=1}^n \| \prm_i^t - \avgprm^t \|^2  = \| \prm^t \|_{\wtk}^2.
\end{equation}

\textbf{Convergence of {\algnamesa}.}
We summarize the convergence rate for {\algnamesa} as follows. The proof can be found in Appendix~\ref{sec:proof}:
\begin{tcolorbox}[boxsep=2pt,left=4pt,right=4pt,top=3pt,bottom=3pt]
\begin{theorem}
\label{thm:main}
Under Assumptions~\ref{assm:lip}, \ref{assm:rand-graph}, \ref{assm:f_var}, \ref{assm:graph_var}. Suppose that the step sizes satisfy
the conditions defined in \eqref{eq:ss_cond}. 
Then, for any $T \geq 1$ with the random stopping iteration ${\sf T} \sim {\rm Unif} \{0,...,T-1\}$, the iterates generated by {\algnamesa} satisfy 
\begin{align}
\expec{ \| \nabla F( \avgprm^{\sf T} ) \|^2 } & \leq \frac{ F_0 - f_\star }{ \alpha T / 8 } + 8 \alpha \bbC_{\sigma } \frac{\sigavg}{n}, \\
\textstyle \expec{ \sum_{i=1}^n \| \prm_i^{\sf T} - \avgprm^{\sf T} \|^2 } & \leq \frac{ F_0 - f_\star }{ \a \gamma \rho_{\min} T / 8 } + \frac{8 \alpha^2 \bbC_{\sigma } \sigavg}{\a \gamma \rho_{\min} n}, 
\end{align}
for any $\a>0$,
where 
$F_0$, $\bbC_{\sigma}$ are defined in \eqref{eq:ft_def_restated}, \eqref{eq:bbC_bound}.
\end{theorem}
\end{tcolorbox}
Setting $\a = {\cal O}(n/\sqrt{T \sigavg})$, $\alpha = \sqrt{n / (T \sigavg)}$ (and assuming $\bar{\sigma}>0$), we have
\begin{equation} \label{eq:crude_bound}
\expec{ \| \nabla F( \avgprm^{\sf T} ) \|^2 } = {\cal O} \left( {\bar{\sigma}} / {\sqrt{nT}} \right),
\end{equation}
which is the same \emph{asymptotic convergence rate} as a centralized SGD algorithm that takes $n$ stochastic gradient samples uniformly from each agent, i.e., linear speedup \citep{lian2017can}. 
Also, using $\a = 1$, the consensus error converges as a rate of $\expec{ \sum_{i=1}^n \| \prm_i^{\sf T} - \avgprm^{\sf T} \|^2 }= \mathcal{O}(n^2 \sigma_A^2 \rho_{\max} / (T \rho_{\min}^2 ))$ under the same step size choice used in \eqref{eq:crude_bound}. Notice that for $T \gg 1$, the effect of random topology only degrades the convergence of consensus error, keeping the transient rate in \eqref{eq:crude_bound} unaffected.
If the gradients are deterministic ($\bar{\sigma}=0$), setting $\a = (L^2 \eta_{\infty} \rho_{\min})^{1/3}$, $\alpha = \alpha_{\infty}$ will yield a better convergence rate as $\expec{ \| \nabla F( \avgprm^{\sf T} ) \|^2 } = {\cal O}(\sigma_A^4 \sqrt{n}/ T)$. Without a transient phase, the error due to random graph and coordinate sparsification is persistent through $\sigma_A^4$ in the above convergence rate.

We further show that the convergence of {\algnamesa} can be accelerated if the objective function of \eqref{eq:main_problem} satisfies the Polyak-Lojasiewicz (PL) condition:
\begin{assumption}
    \label{assm:pl_main} 
    There exists a constant $\mu > 0$ such that $2 \mu (F(\prm) - f_\star) \leq \| \nabla F(\prm) \|^2,~\forall \prm \in \mathbb{R}^d$.
\end{assumption}
\noindent Assumption~\ref{assm:pl_main} includes strongly convex functions as a special case, but also includes other non-convex functions; see \citep{karimi2016linear}. We observe:
\begin{tcolorbox}[boxsep=2pt,left=4pt,right=4pt,top=3pt,bottom=3pt]
\begin{corollary}
\label{cor:pl_main}
Suppose the assumptions and step size conditions in Theorem \ref{thm:main} hold. Furthermore, with Assumption \ref{assm:pl_main}, there exists $\delta \in (0,1)$ such that for any $t \ge 0$, 
\begin{equation}
\mathbb{E}_t [ F_{t+1} - f_\star ] \leq (1-\delta) (F_t - f_\star) + \bbC_{\sigma} \alpha^2 \sigavg / n \label{eq:pl_ineq_main}
\end{equation}
for $F_t, \bbC_\sigma$ defined in \eqref{eq:ft_def_restated}, \eqref{eq:lyapunov_main}, and $\delta = \min \{ {\alpha \mu} / {4}, ~{\gamma \rho_{\min}} / {16},~ {\eta \beta} / ({3 \rho_{\min}}), {\eta} / {12} \}$.
\end{corollary}
\end{tcolorbox}
\noindent The proof can be found in Appendix \ref{app:pl_proof}.
By setting $\alpha = c \ln(T)/(n^2 T)$
in \eqref{eq:pl_ineq_main}, with a carefully chosen $c$ and a sufficiently large $T$ such that $\alpha \leq \alpha_{\infty}$, we can ensure that 
\begin{equation}
    \expec{F(\avgprm^T) - f_\star + \| \prm^T \|^2_{\wtk}} = \mathcal{O}\left( {\sigavg \ln(T)} / ({\mu nT}) \right)
\end{equation}
In the case of deterministic gradient, i.e., $\sigavg = 0$, by setting $\alpha = \alpha_{\infty}$, \eqref{eq:pl_ineq_main} ensures a linear convergence rate of 
$\expec{F(\avgprm^T) - f_\star + \| \prm^T \|^2_{\wtk}} = \mathcal{O}((1-\delta)^T)$,
which shows that the performance of {\algnamesa} is on par with \citep{nedic2017achieving, xu2017convergence}, despite it only requires one round of (sparsified) transmission per iteration.


\vspace{.1cm} \textbf{Convergence  of {\algnamevr}.}
To exploit the benefits of control variates, we need an additional assumption on the stochastic gradient map:
\begin{assumption} \label{assm:sample_lip}
    Each stochastic function $f_i(\cdot; \xi)$ is $\sampL$-smooth in expectation, i.e., for $i = 1, \ldots, n$,
    \begin{equation} \label{eq:f_sample_lip}
        \mathbb{E}_\xi\left[ \|\nabla f_i({\bf x}; \xi) - \nabla f_i({\bf y}; \xi) \|^2 \right] \leq \sampL^2 \|{\bf x} - {\bf y} \|^2 ~\forall~ {\bf x},{\bf y} \in \mathbb{R}^{d}.
    \end{equation}
\end{assumption}
\noindent The above assumption is also known as the mean-square smoothness condition, see \citep{cutkosky2019momentum}, which is strictly stronger than Assumption~\ref{assm:lip}. 
We observe the following convergence guarantee for {\algnamevr}, whose proof can be found in Appendix \ref{app:thm_main_vr_proof}.
\begin{tcolorbox}[boxsep=2pt,left=4pt,right=4pt,top=3pt,bottom=3pt]
    \begin{theorem} \label{thm:main_vr}
    Under Assumptions \ref{assm:lip}, \ref{assm:rand-graph}, \ref{assm:f_var}, \ref{assm:graph_var}, \ref{assm:sample_lip}. Suppose that the step sizes satisfy the conditions in \eqref{eq:storm_ss_cond_start} - \eqref{eq:storm_ss_cond_end}. Then, for any $T \geq 1$ with the random stopping iteration ${\sf T} \sim {\rm Unif} \{0,...,T-1\}$, the iterates generated by {\algnamevr} satisfy
    \begin{align}
        \expec{\| \nabla F(\avgprm^{\sf T}) \|^2 } &\leq \frac{F_0 - f_\star}{T \alpha / 4} + \frac{(\e \cdot 2a_x^2 + \f \cdot 4a_x^2 n) \sigavg}{\alpha / 4},  \label{eq:storm_grad_bound} \\
        \textstyle \expec{ \sum_{i=1}^n \| \prm_i^{\sf T} - \avgprm^{\sf T} \|^2 } &\leq \frac{F_0 - f_\star}{T \a \gamma \rho_{\min} / 8} + \frac{(\e \cdot 2a_x^2 + \f \cdot 4a_x^2 n) \sigavg}{\a \gamma \rho_{\min} / 8} ,
    \end{align}
    where the constants $F_0$, $\a, \e, \f$ are defined in \eqref{eq:potential_storm_def}.
\end{theorem}
\end{tcolorbox}

Setting $\alpha = \mathcal{O}(\bar{\sigma}^{-2/3} T^{-1/3})$, $\eta = \mathcal{O}(n)$, $\gamma = \mathcal{O}(T^{-1/3})$, $\beta = \mathcal{O}(n^{-1} T^{-2/3})$, $a_x = \mathcal{O}(\bar{\sigma}^{-4/3} T^{-2/3})$, $a_\lambda = \mathcal{O}( T^{-1/3} )$, $\f = {\cal O}( n^{-1} T^{1/3} )$ (see \eqref{eq:ss_choice_storm_start} - \eqref{eq:ss_choice_storm_end}), and initializing the algorithm such that $\| {\bf v}^0 \|_{\wtk}^2 = \mathcal{O}(T^{-2/3})$, $ \| \avgmom^0 - (1/n) \oneotimes \nabla {\bf f}(\prm^0) \|^2 = \mathcal{O}(T^{-1/3})$ and $\| \mom^0 - \nabla_\prm {\cal L}(\prm^0, \dprm^0) \|^2 = \mathcal{O}(T^{-1/3})$, we have
\begin{align} \label{eq:storm_rate}
    \expec{\| \nabla F(\avgprm^{\sf T}) \|^2 } &= \mathcal{O}\big( {\bar{\sigma}^{2/3} } / {T^{2/3}} \big).
\end{align}
In regard to the order of $\bar{\sigma}$ and $T$, provided that $n$ is small, the convergence rate of {\algnamevr} matches the lower bound \citep{arjevani2023lower} for non-convex functions under the same smoothness assumption. 
Moreover, by the same choice of step sizes, the consensus error converges at the rate of $\expec{ \sum_{i=1}^n \| \prm_i^{\sf T} - \avgprm^{\sf T} \|^2 } = \mathcal{O}( \bar{\sigma}^{2/3} n \rho_{\min}^{-1} T^{-2/3} )$.
We remark that in \eqref{eq:storm_rate}, the rate remains constant as $n$ increases such that {\algnamevr} does not offer the same {\it linear speedup} observed in Theorem~\ref{thm:main} for {\algnamesa}.
Nevertheless, as $T \gg 1$, the rate of {\algnamevr} will surpass that of {\algnamesa} and other decentralized algorithms on time varying topologies. 

Lastly, we provide detailed discussions on the convergence rates above, e.g., transient time, effects of random topology, etc., {in Appendix~\ref{app:detail_analysis}}.

\subsection{Insight from Analysis: Fixed Point Iteration of {\algnamesa}}

From \eqref{eq:fspdasa_primal}, the following recursive relationship holds for $\avgprm^t$: using the relation ${\bf 1}^\top \wta^\top  = {\bf 0}$, we have
\begin{equation} \label{eq:one-step-descent}
    \textstyle  \avgprm^{t+1} = \avgprm^t - \frac{\alpha}{n} \sum_{i=1}^n \nabla {f}_i( \prm_i^t; \xi_i^{t} ).
\end{equation}
This shows that the evolution of $\{ \avgprm^t \}_{t \geq 0}$ is similar to that of `centralized' SGD applied on \eqref{eq:main_problem} except that the local gradients  are evaluated on the local iterates. However, it is still not straightforward to analyze the convergence of {\algnamesa} as the update of $\prm^t$ involves the dual variable $\dprm^t$ which lacks an intuitive interpretation for constructing the right Lyapunov function. 

To this end, we study the fixed point(s) of \eqref{eq:fspda-sa} to gain insights.
Suppose that for some $t_\star$, the fixed point conditions $\mathbb{E}[ \dprm^{t_\star+1} ~|~ \xi^{:t_\star}] = \dprm^{t_\star}, \mathbb{E}[ \prm^{t_\star+1} ~|~ \xi^{:t_\star}] = \prm^{t_\star}$ hold. Since $\rxii$ is a diagonal matrix with positive diagonal elements, we observe
\begin{equation} \label{eq:dual_fixed}
    \mathbb{E}[ \dprm^{t_\star+1} ~|~ \xi^{:t_\star}] = \dprm^{t_\star} \Longleftrightarrow \rxii \wta \prm^{t_\star} = {\bf 0} \Longleftrightarrow \wta \prm^{t_\star} = {\bf 0},
\end{equation}
On the other hand, the primal update yields
\begin{equation} \label{eq:primal_update_at_fixed}
    \mathbb{E}[ \prm^{t_\star+1} ~|~ \xi^{:t_\star}] = \prm^{t_\star} - \alpha \nabla {\bf f}( \prm^{t_\star} ) - \eta \wta^\top \dprm^{t_\star}.
\end{equation}
Since $\prm_1^{t_\star} = \prm_2^{t_\star} = \cdots = \prm_n^{t_\star}$ at the fixed point (due to \eqref{eq:dual_fixed}), by the consensus condition across two time steps, it implies
\begin{equation} \label{eq:dual_meaning}
    \begin{aligned}
    & \mathbb{E}[ \prm^{ {t_\star} +1} ~|~ \xi^{:t_\star}] - \prm^{t_\star} = ({\bf 1} \otimes {\bf I}_d )(\avgprm^{{t_\star}+1} - \avgprm^{{t_\star}}) \\
    & \textstyle \quad \Longleftrightarrow  \alpha \nabla {\bf f}( \prm^{t_\star} ) + \eta \wta^\top \dprm^{t_\star} = \frac{\alpha}{n} ({\bf 1} {\bf 1}^\top \otimes {\bf I}_d ) \nabla {\bf f}( \prm^{t_\star} ) \\
    & \textstyle \quad \Longleftrightarrow \eta \wta^\top \dprm^{t_\star} = {\alpha} \left( \frac{1}{n} {\bf 1}{\bf 1}^\top - {\bf I}_n \right) \otimes {\bf I}_d ~\nabla {\bf f}( ({\bf 1} \otimes {\bf I}) \avgprm^{{t_\star}} ). 
    \end{aligned}
\end{equation}
From \eqref{eq:dual_meaning}, we see that $\widehat{\dprm}^t$ shall converge to the difference between global and local gradient. Inspired by the above, to facilitate the analysis later, we define 
\begin{equation} \label{eq:v_def_main}
    \textstyle \mathbf{v}^t := \wta^\top \dprm^{t} + \frac{\alpha}{\eta} \nabla {\bf f}(({\bf 1} \otimes {\bf I}) \avgprm^{t}),
\end{equation}
for any $t \geq 0$.
In particular, we see that $\| {\bf v}^t \|^2_{\wtk}$ measures the violation of \eqref{eq:dual_meaning} in tracking the average deterministic gradient using the dual variables. 
The latter will be instrumental in analyzing the consensus error bound, as revealed in Lemma \ref{lemma:consensus}.

\section{Numerical Experiments} \label{sec:num}
This section reports the numerical experiments on practical performance of {\algname}. 
For the time varying topology, we take an extreme setting where for each realization ${\cal G}(\xi_a)$, {only one edge will be selected uniformly at random from ${\cal G}$}. We evaluate the performance with the worst-agent metric, i.e., we present the training loss as $\max_{i \in [n]} F(\prm_i^t)$,  and the stationarity/gradient-norm measure as $\max_{i \in [n]} \| \nabla F(\prm_i^t) \|^2$. This captures the worst-case of the solutions produced by the algorithms.
Unless otherwise specified, all algorithms are initialized with $\prm_i^0 = \avgprm^0$, and for {\algname} we initialize $\widehat{\dprm}^{0} = \momi^{0} = \momidual^{0} = {\bf 0}$, and the stochastic gradients are estimated with a batch size of 256. In the interest of space, omitted details and hyperparameters of the experiments can be found in Appendix~\ref{app:exp}.
The source codes for running all of our experiments are available in \url{https://github.com/OscarYau525/FSPDA}.

\vspace{.1cm} \textbf{MNIST Experiments.} The first set of experiments considers a moderate-scale setting of training a one hidden layer feed-forward neural network with 100 hidden neurons (total number of parameters $d=\text{79,510}$) on the MNIST dataset with $m=60,000$ samples of $784$-dimensional features. 

In the first experiment, we consider the static topology ${\cal G}$ as an Erdos-Renyi graph with connectivity of $p=0.5$ and $n=10$ agents. We compare the proposed {\algnamesa}, {\algnamevr} with six benchmark algorithms utilizing different types of time-varying topology. 
Among them, {\tt DSGD} \citep{koloskova2020unified} and {\tt Swarm-SGD} \citep{nadiradze2021asynchronous} use the general time varying topology setting as {\algname} where each edge of ${\cal G}(\xi_a)$ is active uniformly at random, in addition to random sparsification used {\algnamesa} and adaptive quantized used in {\tt Swarm-SGD}; {\tt CHOCO-SGD} \citep{koloskova2019decentralized} takes ${\cal G}(\xi_a)$ as an broadcasting subgraph where one agent selects all his/her neighbors; {\tt Decen-Scaffnew} \citep{mishchenko2022proxskip}, {\tt LED} \citep{alghunaim2024local}, and {\tt K-GT} \citep{liu2024decentralized} utilize local updates where ${\cal G}(\xi_a)$ is either taken as an empty topology, or as the static topology ${\cal G}$. 
We configure these algorithms such that they have the same communication cost (in terms of bits transmitted over network) \emph{on average}. For instance, the local update algorithms ({\tt Decen-Scaffnew}, {\tt LED}, {\tt K-GT}) only communicate once using ${\cal G}$ every $\mathcal{O}\left(\frac{|{\cal E}|d}{k}\right)$ iterations to match the communication cost of $k$-coordinate sparse one-edge random graph used in {\algname}.

The local objective function held by each agent is the cross-entropy classification loss on a local dataset with $m_i = 6000$ samples, plus a regularization loss $\frac{\lambda}{2} \| \prm_i\|^2$ with $\lambda = 10^{-4}$, where $\prm_i$ are the weight parameters of the feed-forward neural network classifier. We split the training set into $n=10$ disjoint sets such that each set contains only one class label and assign each set to one agent as its local dataset. Note that as we do not shuffle the data samples across local datasets, the local objective function held by different agents will become highly heterogeneous. 

\begin{figure}[t]
    \centering
    \includegraphics[width=0.97\textwidth]{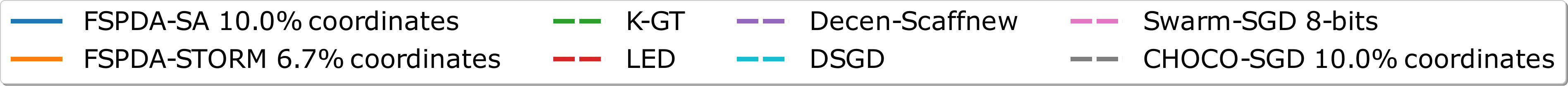} \\
    \vspace{0.25cm}
    \includegraphics[width=0.25\textwidth]{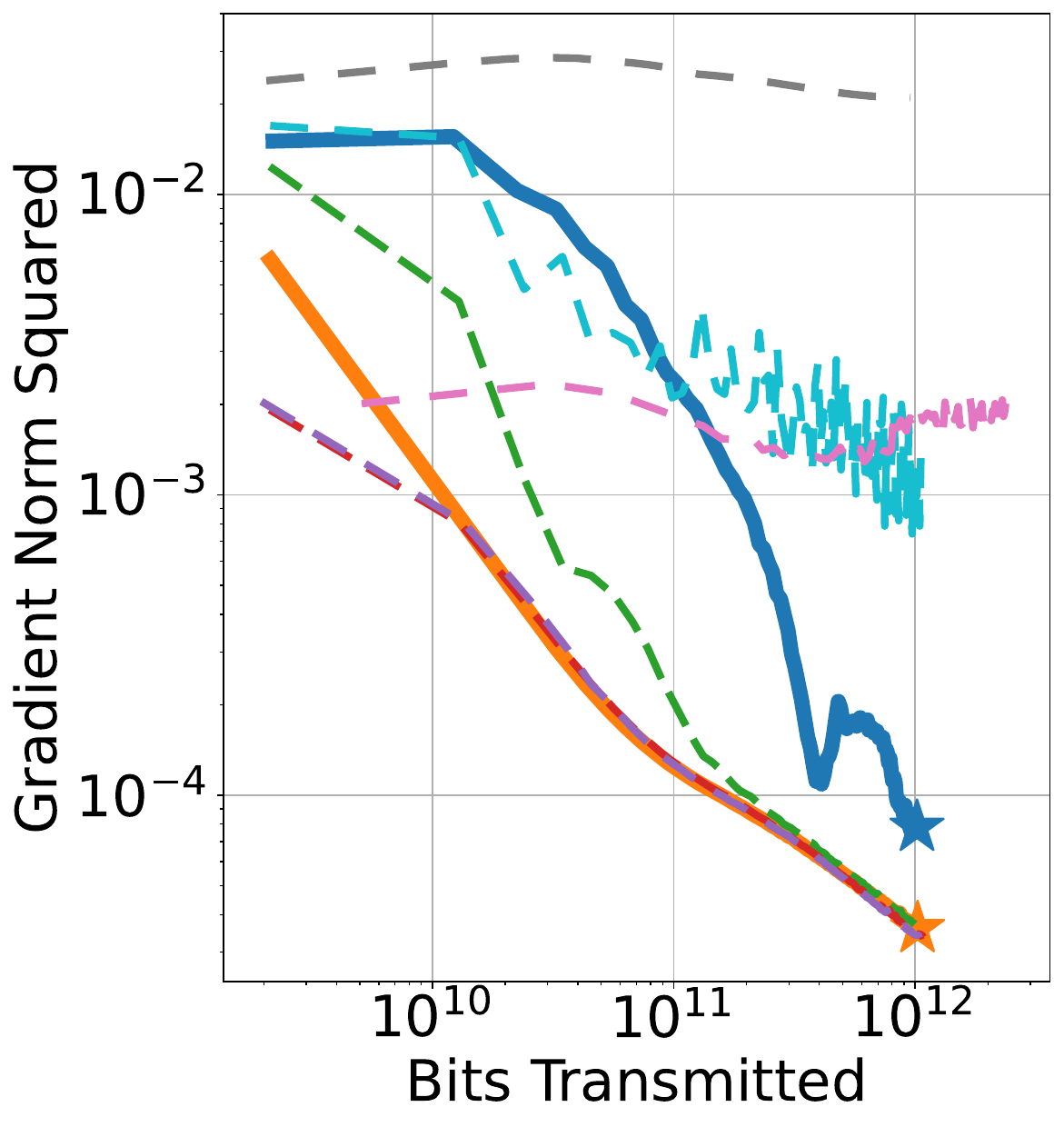}
    \includegraphics[width=0.25\textwidth]{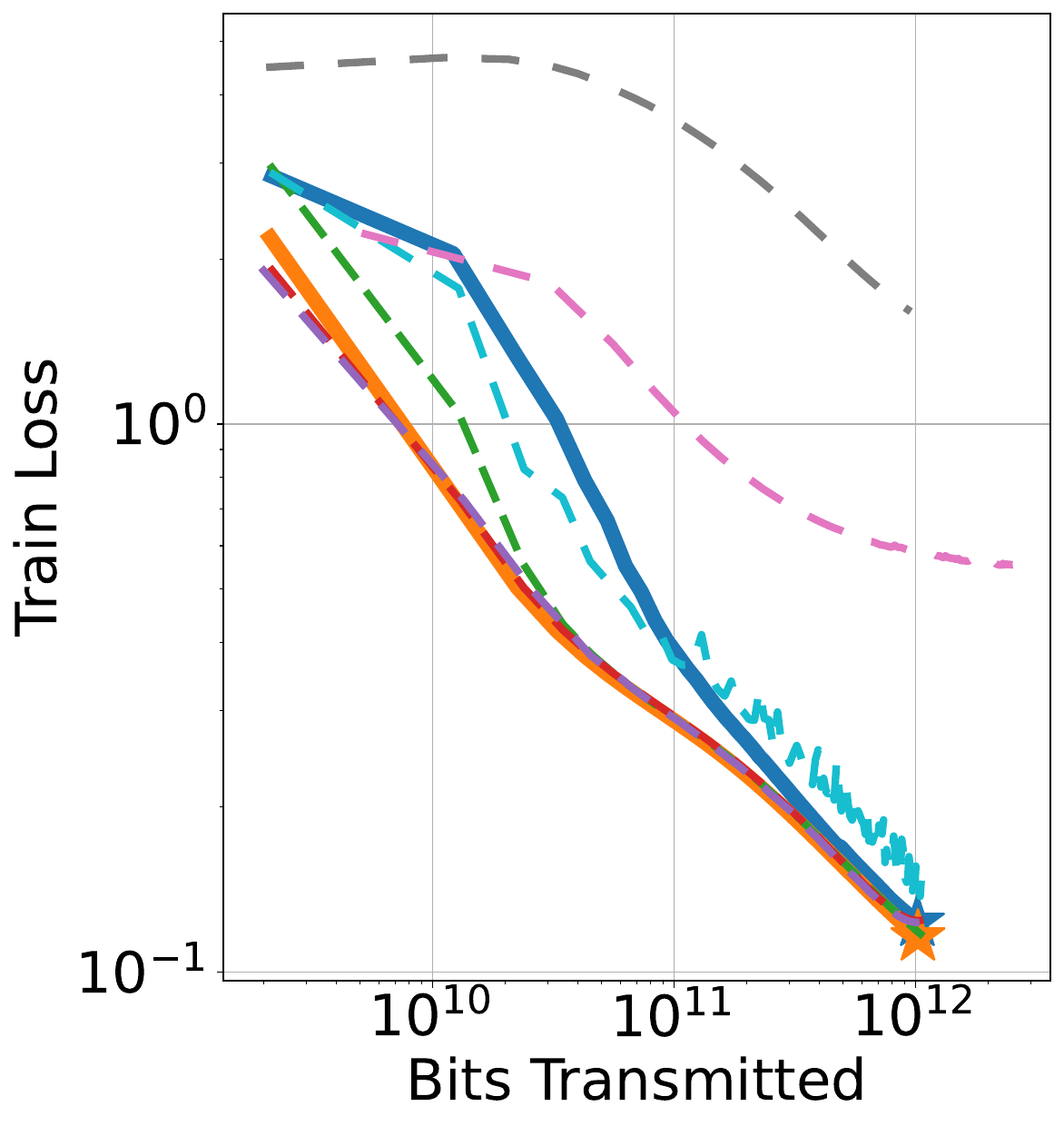}
    \includegraphics[width=0.25\textwidth]{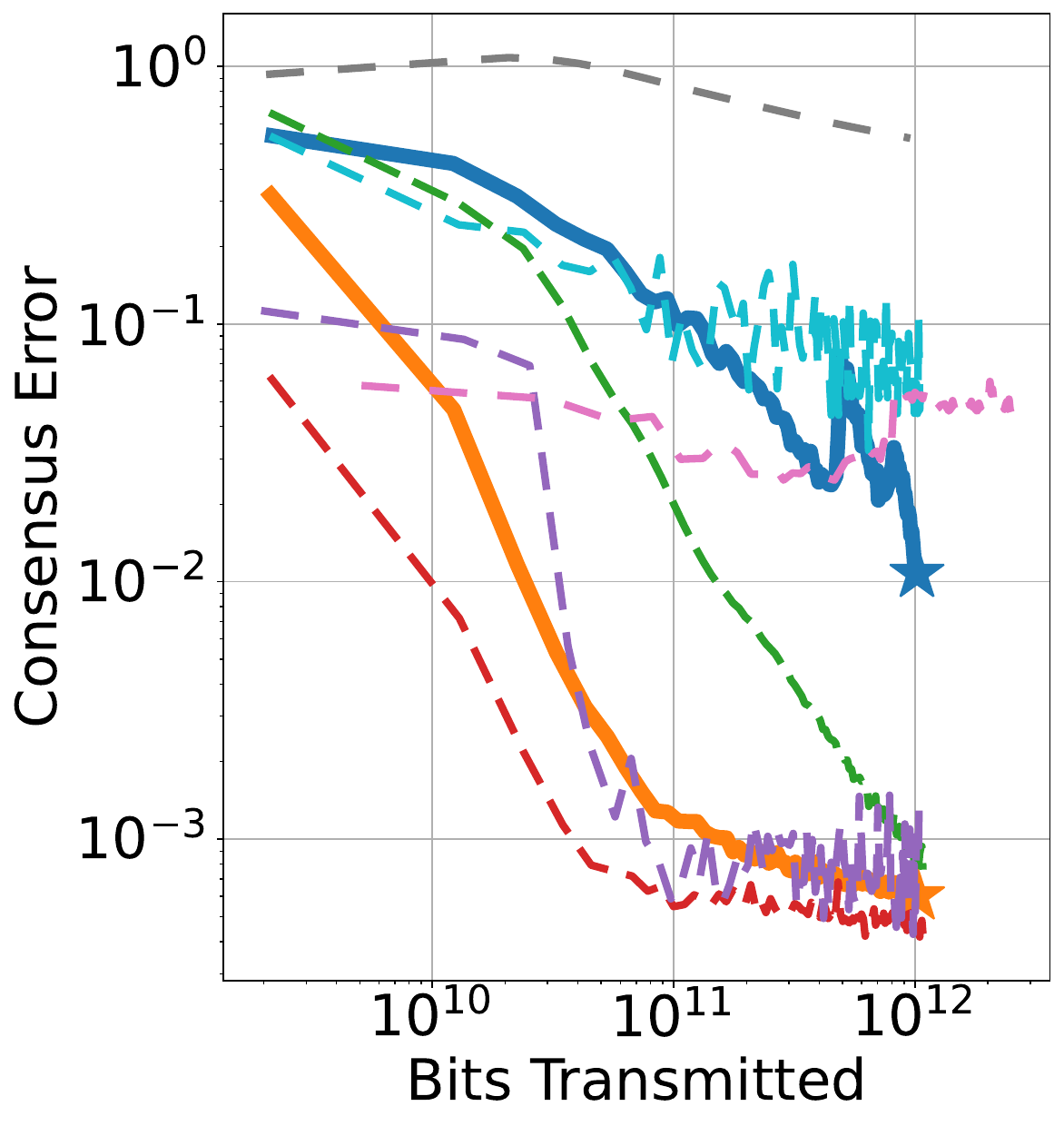}
    \caption{Feed-forward neural network classification training on MNIST using $10^6$ iterations.
    }\vspace{-.2cm}
    \label{fig:mnist_new}
\end{figure}

Fig.~\ref{fig:mnist_new} compares the squared gradient norm, training loss, consensus error of the benchmarked algorithms. 
We first note that both {\algname} algorithms have significantly outperformed {\tt DSGD}, {\tt Swarm-SGD} on the general time varying topology as well as {\tt CHOCO-SGD}. Meanwhile, the performance of {\algname} is comparable to the local update algorithms {\tt Decen-Scaffnew}, {\tt LED}, {\tt K-GT}. Notice that the latter require additional synchronization steps which may not be suitable for random networks. 
Lastly, we notice that as $T \gg 1$, {\algnamevr} can slightly outperform {\algnamesa} due to its ${\cal O}(1/T^{2/3})$ rate as shown in our analysis. 
We further expand the experiments by a series of ablation studies over data heterogeneity, sparsity levels, graph topologies, gradient noise and dual momentum in Appendix \ref{app:ablation}.

\vspace{.1cm} \textbf{Imagenet Experiments.}
The second set of experiments consider a large-scale setting for training a Resnet-50 network (total number of parameters $d=\text{25,557,032}$) on the Imagenet dataset (training dataset of 1,281,168 images from 100 classes, re-scaled and cropped to 256 $\times$ 256 image dimensions).
We consider cross-entropy classification loss plus the same L2 norm regularization loss as in the previous setup.
We split the dataset across a network of $n=8$ nodes where the static graph ${\cal G}$ is taken as the fully connected topology. The performance metrics are measured at the network average iterate $\avgprm^t$. Inspired by \cite[Eq.~(5)]{loshchilov2016sgdr} we adopt a cosine learning rate scheduling with 5 epochs of linear warm up for every algorithm. In particular, the step sizes $\alpha, \eta$ of {\algnamesa} are scheduled simultaneously such that $\alpha_t / \eta_t$ remains constant, as illustrated in Appendix \ref{app:exp}. We draw a batch of 128 samples to estimate the stochastic gradient.

\begin{figure}
    \centering
    \includegraphics[width=0.6\textwidth]{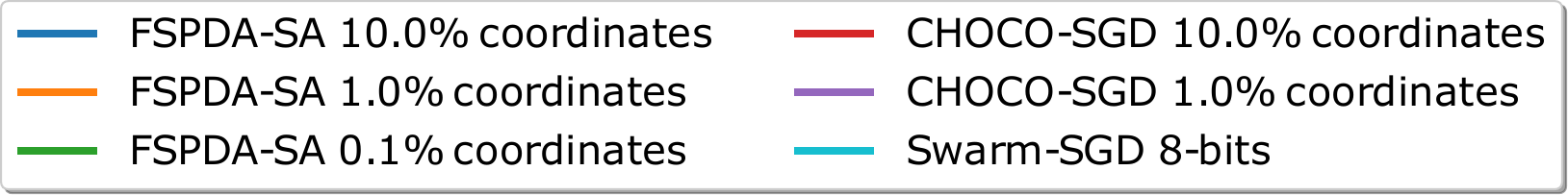}\\
    \includegraphics[width=0.22\textwidth]{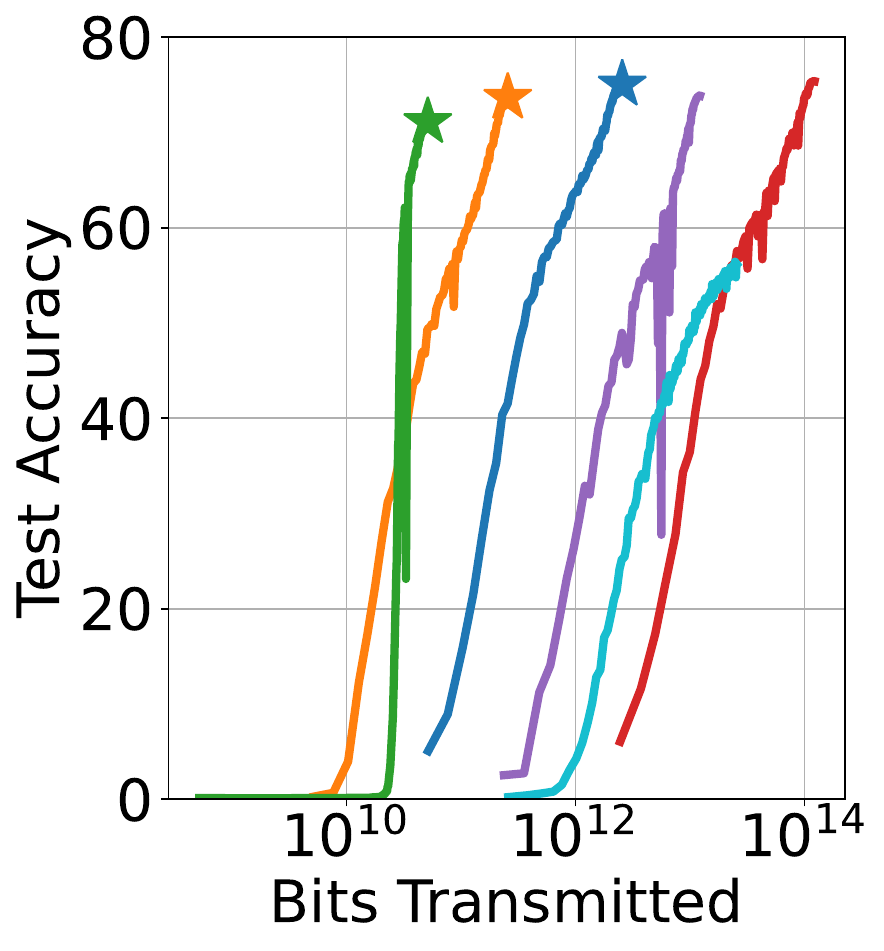}
    \includegraphics[width=0.22\textwidth]{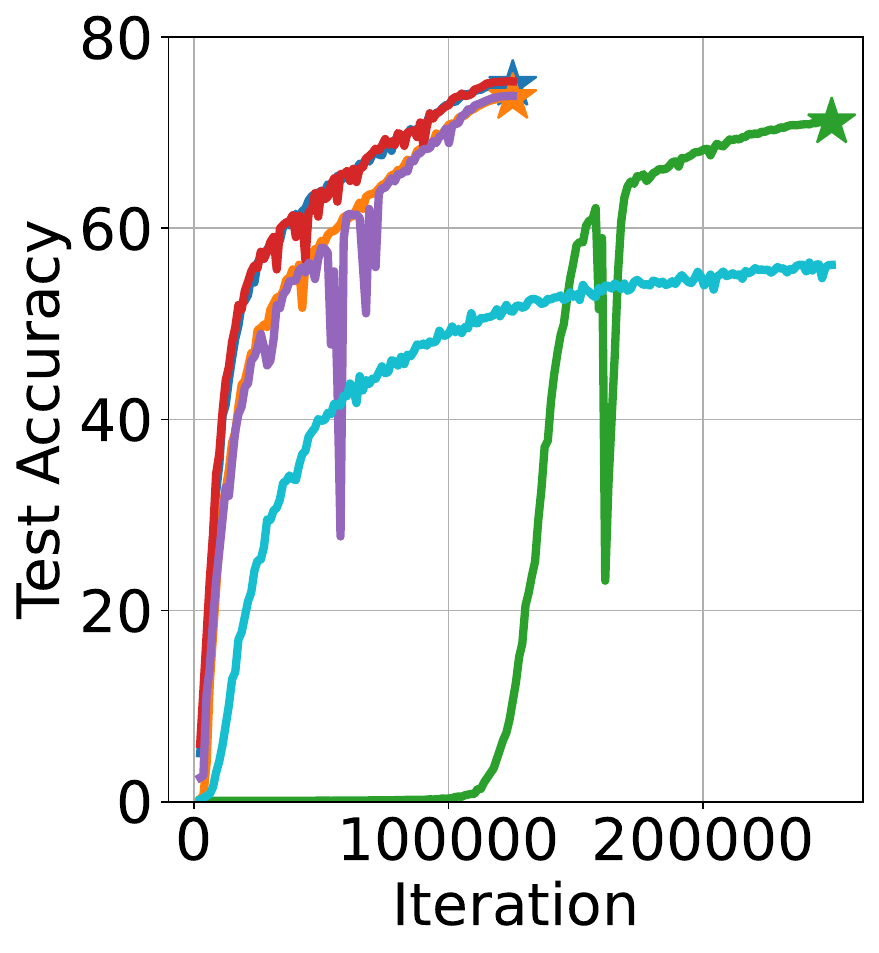}
    \includegraphics[width=0.22\textwidth]{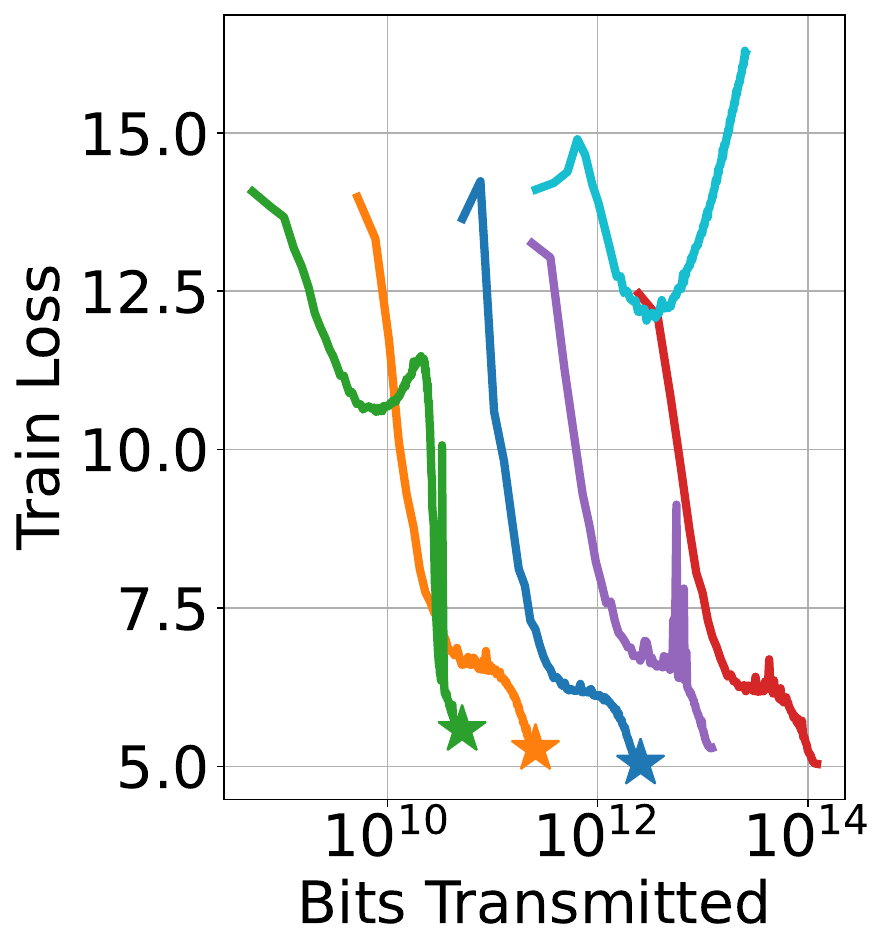}
    \includegraphics[width=0.22\textwidth]{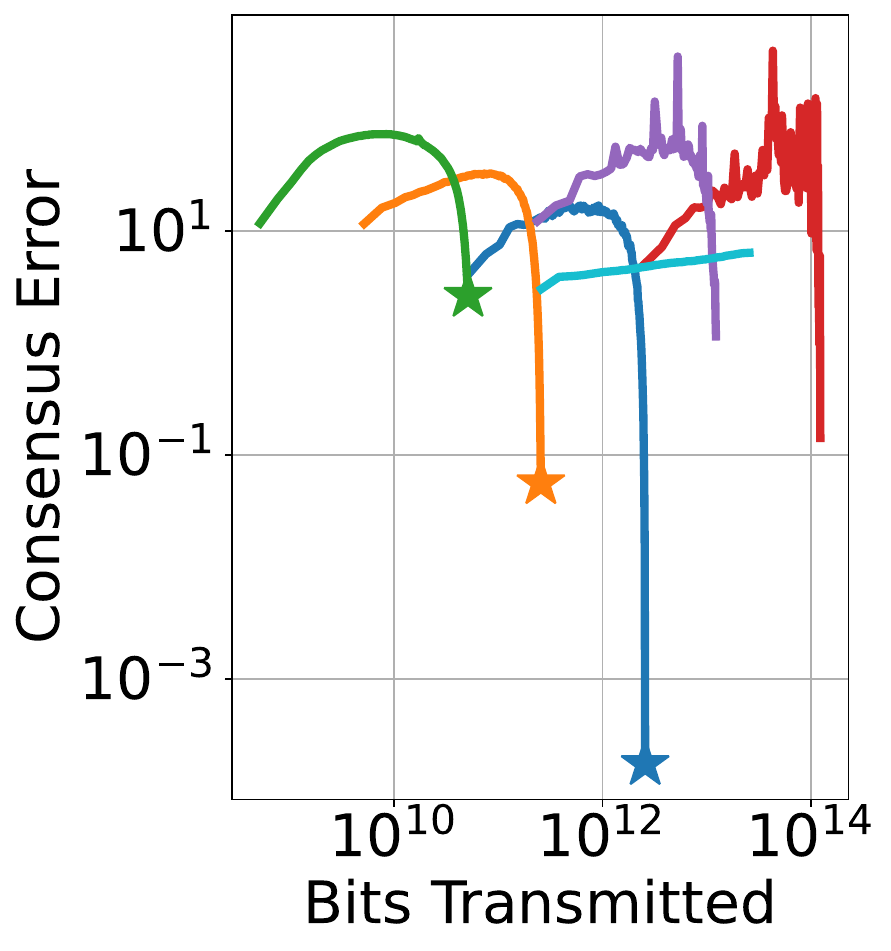}
    \caption{Resnet-50 classification training on Imagenet.
    }\vspace{-.4cm}
    \label{fig:imagenet}
\end{figure}

We focus on the communication efficiency and only compare {\algnamesa}, {\tt CHOCO-SGD}, {\tt Swarm-SGD} in this experiment due to limited resources. The results are reported in Figure \ref{fig:imagenet} that compare the test accuracy and training loss against iteration number and bits transmitted. When compared with {\tt CHOCO-SGD}, {\algnamesa} achieves almost the same accuracy using one-edge random graphs with at least 100x reduction in communication cost on 100 epoch training. Also notice that further compressing the communication to 0.1\% sparse coordinates in {\algnamesa} requires more training epochs to recover the same level of accuracy.

\section{Conclusion}
This paper proposed a fully stochastic primal dual gradient algorithm ({\algname}) framework for decentralized optimization over arbitrarily time varying random networks. We utilize a new stochastic augmented Lagrangian function and apply SA to search for its saddle point. We develop two algorithms, one is by plain SA ({\algnamesa}), and one uses control variates for variance reduction ({\algnamevr}). We prove that both algorithms achieve state-of-the-art convergence rates, while relaxing assumptions on both bounded heterogeneity and the type of time varying topologies.

\bibliographystyle{plainnat}
\bibliography{ref.bib}

\begin{thebibliography}{38}
\providecommand{\natexlab}[1]{#1}
\providecommand{\url}[1]{\texttt{#1}}
\expandafter\ifx\csname urlstyle\endcsname\relax
  \providecommand{\doi}[1]{doi: #1}\else
  \providecommand{\doi}{doi: \begingroup \urlstyle{rm}\Url}\fi

\bibitem[Alghunaim(2024)]{alghunaim2024local}
Sulaiman~A Alghunaim.
\newblock Local exact-diffusion for decentralized optimization and learning.
\newblock \emph{IEEE Transactions on Automatic Control}, 2024.

\bibitem[Arjevani et~al.(2023)Arjevani, Carmon, Duchi, Foster, Srebro, and Woodworth]{arjevani2023lower}
Yossi Arjevani, Yair Carmon, John~C Duchi, Dylan~J Foster, Nathan Srebro, and Blake Woodworth.
\newblock Lower bounds for non-convex stochastic optimization.
\newblock \emph{Mathematical Programming}, 199\penalty0 (1):\penalty0 165--214, 2023.

\bibitem[Bertsekas(2016)]{bertsekas2016nonlinear}
Dimitri Bertsekas.
\newblock \emph{Nonlinear Programming}, volume~4.
\newblock Athena Scientific, 2016.

\bibitem[Bianchi et~al.(2021)Bianchi, Hachem, and Salim]{bianchi2021fully}
Pascal Bianchi, Walid Hachem, and Adil Salim.
\newblock A fully stochastic primal-dual algorithm.
\newblock \emph{Optimization Letters}, 15\penalty0 (2):\penalty0 701--710, 2021.

\bibitem[Chang et~al.(2020)Chang, Hong, Wai, Zhang, and Lu]{chang2020distributed}
Tsung-Hui Chang, Mingyi Hong, Hoi-To Wai, Xinwei Zhang, and Songtao Lu.
\newblock Distributed learning in the nonconvex world: From batch data to streaming and beyond.
\newblock \emph{IEEE Signal Processing Magazine}, 37\penalty0 (3):\penalty0 26--38, 2020.

\bibitem[Cutkosky and Orabona(2019)]{cutkosky2019momentum}
Ashok Cutkosky and Francesco Orabona.
\newblock Momentum-based variance reduction in non-convex sgd.
\newblock \emph{Advances in neural information processing systems}, 32, 2019.

\bibitem[Guo et~al.(2023)Guo, Alghunaim, Yuan, Condat, and Cao]{guo2023revisiting}
Luyao Guo, Sulaiman~A Alghunaim, Kun Yuan, Laurent Condat, and Jinde Cao.
\newblock Revisiting decentralized proxskip: Achieving linear speedup.
\newblock \emph{arXiv preprint arXiv:2310.07983}, 2023.

\bibitem[Hajinezhad and Hong(2019)]{hajinezhad2019perturbed}
Davood Hajinezhad and Mingyi Hong.
\newblock Perturbed proximal primal--dual algorithm for nonconvex nonsmooth optimization.
\newblock \emph{Mathematical Programming}, 176\penalty0 (1):\penalty0 207--245, 2019.

\bibitem[Hong et~al.(2017)Hong, Hajinezhad, and Zhao]{hong2017prox}
Mingyi Hong, Davood Hajinezhad, and Ming-Min Zhao.
\newblock Prox-pda: The proximal primal-dual algorithm for fast distributed nonconvex optimization and learning over networks.
\newblock In \emph{International Conference on Machine Learning}, pages 1529--1538. PMLR, 2017.

\bibitem[Kairouz et~al.(2021)Kairouz, McMahan, Avent, Bellet, Bennis, Bhagoji, Bonawitz, Charles, Cormode, Cummings, et~al.]{kairouz2021advances}
Peter Kairouz, H~Brendan McMahan, Brendan Avent, Aur{\'e}lien Bellet, Mehdi Bennis, Arjun~Nitin Bhagoji, Kallista Bonawitz, Zachary Charles, Graham Cormode, Rachel Cummings, et~al.
\newblock Advances and open problems in federated learning.
\newblock \emph{Foundations and trends{\textregistered} in machine learning}, 14\penalty0 (1--2):\penalty0 1--210, 2021.

\bibitem[Karimi et~al.(2016)Karimi, Nutini, and Schmidt]{karimi2016linear}
Hamed Karimi, Julie Nutini, and Mark Schmidt.
\newblock Linear convergence of gradient and proximal-gradient methods under the polyak-{\l}ojasiewicz condition.
\newblock In \emph{Machine Learning and Knowledge Discovery in Databases: European Conference, ECML PKDD 2016, Riva del Garda, Italy, September 19-23, 2016, Proceedings, Part I 16}, pages 795--811. Springer, 2016.

\bibitem[Koloskova et~al.(2019{\natexlab{a}})Koloskova, Lin, Stich, and Jaggi]{koloskova2019decentralized2}
Anastasia Koloskova, Tao Lin, Sebastian~U Stich, and Martin Jaggi.
\newblock Decentralized deep learning with arbitrary communication compression.
\newblock In \emph{International Conference on Learning Representations}, 2019{\natexlab{a}}.

\bibitem[Koloskova et~al.(2019{\natexlab{b}})Koloskova, Stich, and Jaggi]{koloskova2019decentralized}
Anastasia Koloskova, Sebastian Stich, and Martin Jaggi.
\newblock Decentralized stochastic optimization and gossip algorithms with compressed communication.
\newblock In \emph{International Conference on Machine Learning}, pages 3478--3487. PMLR, 2019{\natexlab{b}}.

\bibitem[Koloskova et~al.(2020)Koloskova, Loizou, Boreiri, Jaggi, and Stich]{koloskova2020unified}
Anastasia Koloskova, Nicolas Loizou, Sadra Boreiri, Martin Jaggi, and Sebastian Stich.
\newblock A unified theory of decentralized sgd with changing topology and local updates.
\newblock In \emph{International Conference on Machine Learning}, pages 5381--5393. PMLR, 2020.

\bibitem[Koloskova et~al.(2021)Koloskova, Lin, and Stich]{koloskova2021improved}
Anastasiia Koloskova, Tao Lin, and Sebastian~U Stich.
\newblock An improved analysis of gradient tracking for decentralized machine learning.
\newblock \emph{Advances in Neural Information Processing Systems}, 34:\penalty0 11422--11435, 2021.

\bibitem[Kovalev et~al.(2021)Kovalev, Gasanov, Gasnikov, and Richtarik]{kovalev2021lower}
Dmitry Kovalev, Elnur Gasanov, Alexander Gasnikov, and Peter Richtarik.
\newblock Lower bounds and optimal algorithms for smooth and strongly convex decentralized optimization over time-varying networks.
\newblock \emph{Advances in Neural Information Processing Systems}, 34:\penalty0 22325--22335, 2021.

\bibitem[Kovalev et~al.(2024)Kovalev, Borodich, Gasnikov, and Feoktistov]{kovalev2024lower}
Dmitry Kovalev, Ekaterina Borodich, Alexander Gasnikov, and Dmitrii Feoktistov.
\newblock Lower bounds and optimal algorithms for non-smooth convex decentralized optimization over time-varying networks.
\newblock \emph{arXiv preprint arXiv:2405.18031}, 2024.

\bibitem[Lei et~al.(2018)Lei, Chen, and Fang]{lei2018asymptotic}
Jinlong Lei, Han-Fu Chen, and Hai-Tao Fang.
\newblock Asymptotic properties of primal-dual algorithm for distributed stochastic optimization over random networks with imperfect communications.
\newblock \emph{SIAM Journal on Control and Optimization}, 56\penalty0 (3):\penalty0 2159--2188, 2018.

\bibitem[Li and Lin(2024)]{li2024accelerated}
Huan Li and Zhouchen Lin.
\newblock Accelerated gradient tracking over time-varying graphs for decentralized optimization.
\newblock \emph{Journal of Machine Learning Research}, 25\penalty0 (274):\penalty0 1--52, 2024.

\bibitem[Lian et~al.(2017)Lian, Zhang, Zhang, Hsieh, Zhang, and Liu]{lian2017can}
Xiangru Lian, Ce~Zhang, Huan Zhang, Cho-Jui Hsieh, Wei Zhang, and Ji~Liu.
\newblock Can decentralized algorithms outperform centralized algorithms? a case study for decentralized parallel stochastic gradient descent.
\newblock \emph{Advances in neural information processing systems}, 30, 2017.

\bibitem[Liu et~al.(2024)Liu, Lin, Koloskova, and Stich]{liu2024decentralized}
Yue Liu, Tao Lin, Anastasia Koloskova, and Sebastian~U Stich.
\newblock Decentralized gradient tracking with local steps.
\newblock \emph{Optimization Methods and Software}, pages 1--28, 2024.

\bibitem[Lobel and Ozdaglar(2010)]{lobel2010distributed}
Ilan Lobel and Asuman Ozdaglar.
\newblock Distributed subgradient methods for convex optimization over random networks.
\newblock \emph{IEEE Transactions on Automatic Control}, 56\penalty0 (6):\penalty0 1291--1306, 2010.

\bibitem[Lorenzo and Scutari(2016)]{di2016next}
Paolo~Di Lorenzo and Gesualdo Scutari.
\newblock Next: In-network nonconvex optimization.
\newblock \emph{IEEE Transactions on Signal and Information Processing over Networks}, 2\penalty0 (2):\penalty0 120--136, 2016.

\bibitem[Loshchilov and Hutter(2016)]{loshchilov2016sgdr}
Ilya Loshchilov and Frank Hutter.
\newblock Sgdr: Stochastic gradient descent with warm restarts.
\newblock \emph{arXiv preprint arXiv:1608.03983}, 2016.

\bibitem[Lu et~al.(2019)Lu, Zhang, Sun, and Hong]{lu2019gnsd}
Songtao Lu, Xinwei Zhang, Haoran Sun, and Mingyi Hong.
\newblock Gnsd: A gradient-tracking based nonconvex stochastic algorithm for decentralized optimization.
\newblock In \emph{2019 IEEE Data Science Workshop (DSW)}, pages 315--321. IEEE, 2019.

\bibitem[Mishchenko et~al.(2022)Mishchenko, Malinovsky, Stich, and Richt{\'a}rik]{mishchenko2022proxskip}
Konstantin Mishchenko, Grigory Malinovsky, Sebastian Stich, and Peter Richt{\'a}rik.
\newblock Proxskip: Yes! local gradient steps provably lead to communication acceleration! finally!
\newblock In \emph{International Conference on Machine Learning}, pages 15750--15769. PMLR, 2022.

\bibitem[Nadiradze et~al.(2021)Nadiradze, Sabour, Davies, Li, and Alistarh]{nadiradze2021asynchronous}
Giorgi Nadiradze, Amirmojtaba Sabour, Peter Davies, Shigang Li, and Dan Alistarh.
\newblock Asynchronous decentralized sgd with quantized and local updates.
\newblock \emph{Advances in Neural Information Processing Systems}, 34:\penalty0 6829--6842, 2021.

\bibitem[Nedic and Ozdaglar(2009)]{nedic2009distributed}
Angelia Nedic and Asuman Ozdaglar.
\newblock Distributed subgradient methods for multi-agent optimization.
\newblock \emph{IEEE Transactions on Automatic Control}, 54\penalty0 (1):\penalty0 48--61, 2009.

\bibitem[Nedic et~al.(2017)Nedic, Olshevsky, and Shi]{nedic2017achieving}
Angelia Nedic, Alex Olshevsky, and Wei Shi.
\newblock Achieving geometric convergence for distributed optimization over time-varying graphs.
\newblock \emph{SIAM Journal on Optimization}, 27\penalty0 (4):\penalty0 2597--2633, 2017.

\bibitem[Pu et~al.(2021)Pu, Olshevsky, and Paschalidis]{pu2021sharp}
Shi Pu, Alex Olshevsky, and Ioannis~Ch Paschalidis.
\newblock A sharp estimate on the transient time of distributed stochastic gradient descent.
\newblock \emph{IEEE Transactions on Automatic Control}, 67\penalty0 (11):\penalty0 5900--5915, 2021.

\bibitem[Qin et~al.(2021)Qin, Etesami, and Uribe]{qin2021communication}
Tiancheng Qin, S~Rasoul Etesami, and C{\'e}sar~A Uribe.
\newblock Communication-efficient decentralized local sgd over undirected networks.
\newblock In \emph{2021 60th IEEE Conference on Decision and Control (CDC)}, pages 3361--3366. IEEE, 2021.

\bibitem[Qu and Li(2017)]{qu2017harnessing}
Guannan Qu and Na~Li.
\newblock Harnessing smoothness to accelerate distributed optimization.
\newblock \emph{IEEE Transactions on Control of Network Systems}, 5\penalty0 (3):\penalty0 1245--1260, 2017.

\bibitem[Ram et~al.(2010)Ram, Nedi{\'c}, and Veeravalli]{sundhar2010distributed}
S~Sundhar Ram, Angelia Nedi{\'c}, and Venugopal~V Veeravalli.
\newblock Distributed stochastic subgradient projection algorithms for convex optimization.
\newblock \emph{Journal of optimization theory and applications}, 147:\penalty0 516--545, 2010.

\bibitem[Shi et~al.(2015)Shi, Ling, Wu, and Yin]{shi2015extra}
Wei Shi, Qing Ling, Gang Wu, and Wotao Yin.
\newblock Extra: An exact first-order algorithm for decentralized consensus optimization.
\newblock \emph{SIAM Journal on Optimization}, 25\penalty0 (2):\penalty0 944--966, 2015.

\bibitem[Xu et~al.(2017)Xu, Zhu, Soh, and Xie]{xu2017convergence}
Jinming Xu, Shanying Zhu, Yeng~Chai Soh, and Lihua Xie.
\newblock Convergence of asynchronous distributed gradient methods over stochastic networks.
\newblock \emph{IEEE Transactions on Automatic Control}, 63\penalty0 (2):\penalty0 434--448, 2017.

\bibitem[Yau and Wai(2023)]{yau2023fully}
Chung-Yiu Yau and Hoi-To Wai.
\newblock Fully stochastic distributed convex optimization on time-varying graph with compression.
\newblock In \emph{2023 62nd IEEE Conference on Decision and Control (CDC)}, pages 145--150. IEEE, 2023.

\bibitem[Yi et~al.(2018)Yi, Yao, Yang, George, and Johansson]{yi2018distributed}
Xinlei Yi, Lisha Yao, Tao Yang, Jemin George, and Karl~H Johansson.
\newblock Distributed optimization for second-order multi-agent systems with dynamic event-triggered communication.
\newblock In \emph{2018 IEEE Conference on Decision and Control (CDC)}, pages 3397--3402. IEEE, 2018.

\bibitem[Yi et~al.(2021)Yi, Zhang, Yang, Chai, and Johansson]{yi2021linear}
Xinlei Yi, Shengjun Zhang, Tao Yang, Tianyou Chai, and Karl~H Johansson.
\newblock Linear convergence of first-and zeroth-order primal--dual algorithms for distributed nonconvex optimization.
\newblock \emph{IEEE Transactions on Automatic Control}, 67\penalty0 (8):\penalty0 4194--4201, 2021.

\end{thebibliography}

\newpage

\appendix
\newpage
\section{Asynchronous Implementation of {\algname}} \label{app:async}
We show that {\algname} algorithms can be implemented in an asynchronous manner.
For illustration purpose, we concentrate on {\algnamesa} and a scenario where any sparsified communication round is almost delay-free and only one edge is active at a time, also known as the pairwise gossip setting, while the stochastic gradient computation constitutes the major synchronization overhead in the algorithm. 
We describe the implementation of \eqref{eq:fspdasa_primal}, \eqref{eq:fspdasa_dual} by the pseudo-code in Algorithm~\ref{alg:async_fsppd} from the perspective of the $i$-th agent. Notice that the threads {\cmthread} and {\cpthread} are persistently running \emph{in parallel} at each agent. Moreover, the local variables ${\cal B}_i$, $g_i$, $t_i$ are held by agent $i$ and updated by both threads. 

The thread {\cmthread} aims at preparing the message $\sum_{j \in {\cal N}_i(\xi_a^t)} {\bf C}_{ij}(\xi_a^t) ( \prm_j^t - \prm_i^t )$ needed at \eqref{eq:fspdasa_primal}, \eqref{eq:fspdasa_dual} using a gossip-like step. 
At each run by agent $i$, the agent checks if s/he is connected to any active neighbor on the current graph ${\cal G}^t = {\cal G}( \xi_a^t ) := ( {\cal V}, {\cal E}( \xi_a^t ) )$. If the pair of agents $(i,j)$ are connected, they will exchange the sparsified decision variable $\{ (\prm^{t_i}_{i,k}, \prm^{t_j}_{j,k}) \}_{k \in {\cal I}_{{ij}}(\xi_a^t)}$ and make preparation for the computation thread. The neighbor $j$ whose communicated with agent $i$ will be added to the buffer ${\cal B}_i$. Notice that the protocol is designed such that the asynchronous implementation of {\algnamesa} aligns with \eqref{eq:fspdasa_primal}, \eqref{eq:fspdasa_dual}. For instance, our asynchronous gradient model with $b_i(\xi^t) = 0$ allows agent $i$ to skip gradient computation and the sparse extended graph model \eqref{eq:local-compute} allows agent $i$ to skip communication, as implemented in [L\ref{line:asyncgrad}, Alg. \ref{alg:comm_thread}] which serves as an update after the \emph{idle} state of agent $i$ to compensate for the missed iterations led by neighbor $j \in {\cal B}_i$ using the locally stored $\widehat{\dprm}_i^{t_i}$.

The thread {\cpthread} aims at executing the primal-dual steps \eqref{eq:fspdasa_primal}, \eqref{eq:fspdasa_dual} with local updates, i.e., updating using local stochastic gradient before the next round of communication. Note that in federated learning, local update has been used extensively which led to significant performance improvement; see \citep{kairouz2021advances, qin2021communication}. As mentioned, the stochastic gradient (SG) computation in [L\ref{line:SGcompute}, Alg. \ref{alg:comp_thread}] is the major bottleneck of the algorithm. 
Upon the completion of [L\ref{line:SGcompute}, Alg. \ref{alg:comp_thread}], the two cases in [L\ref{line:caseA}, Alg. \ref{alg:comp_thread}] \& [L\ref{line:caseB}, Alg. \ref{alg:comp_thread}] essentially implement \eqref{eq:fspdasa_primal}, \eqref{eq:fspdasa_dual}. For the case of [L\ref{line:caseB}, Alg. \ref{alg:comp_thread}] where the communication buffer ${\cal B}_i$ is non-empty, we perform a sparse gossip with SG update.
Upon completing this round of primal-dual update, the communication buffer ${\cal B}_i$ will be cleared. 

Overall, by running the two threads persistently at each agent, the decentralized system effectively implements {\algnamesa} in \eqref{eq:fspdasa_primal}, \eqref{eq:fspdasa_dual} as an asynchronous algorithm. The same asynchronous implementation can be easily extended to {\algnamevr}.

\begin{algorithm}
\caption{{\algname} from Agent $i$'s Perspective} 
\begin{algorithmic}[1]\label{alg:async_fsppd}
\STATE {\bfseries input:} Iteration number $T$. 
\STATE {\bf local variable (initialize):} Communication buffer ${\cal B}_i$; gradient counter $g_i = 0$; iteration counter $t_i = 0$.
\WHILE{$\max_{j \in [n]} t_j < T$ \textbf{in parallel}}
       \STATE \textit{communication\_thread()}; see Algorithm~\ref{alg:comm_thread}.
       \STATE \textit{computation\_thread()}; see Algorithm~\ref{alg:comp_thread}.
\ENDWHILE
\end{algorithmic}
\end{algorithm}

\begin{algorithm}
\caption{\textit{communication\_thread()} of Agent $i$} 
\begin{algorithmic}[1]\label{alg:comm_thread}
\STATE {\bf local variable:} buffer ${\cal B}_i$; counters $g_i$, $t_i$.
\FOR{$(i,j) \in {\cal E}( \xi_a^t )$}
    \IF{${\cal B}_i = \emptyset$ and ${\cal B}_j = \emptyset$}
    \STATE Agents $i,j$ exchanges $t_i, t_j$.
    \IF{$t_i < t_j$}
        \STATE Interrupt [L\ref{line:SGcompute}, Alg. \ref{alg:comp_thread}] of \textit{computation\_thread()} and run [L\ref{line:local_grad_step}, Alg. \ref{alg:comp_thread}] to consume (any) SG buffer. \\
        \STATE \ding{229} \emph{Local non-SG step:} evaluate
        \begin{align}
        \prm_i^{t_j} = \prm_i^{t_i} - (t_j - t_i) \eta \widehat{\dprm}^{t_i}_i, 
        \end{align}
        and set $\widehat{\dprm}^{{t_j}}_i = \widehat{\dprm}^{t_i}_i$, $t_i \leftarrow t_j$. \label{line:asyncgrad}
    \ENDIF
    \STATE {\tt //* {begin streaming index-value pairs} *//}
    \FOR{$k \in {\cal I}_{{ij}}(\xi_a^t) $}
        \STATE Agent $j$ receives $(k, \prm_{i,k}^{t_i})$ from agent $i$.
        \STATE Agent $i$ receives $\prm_{j,k}^{t_j}$ from agent $j$.
    \ENDFOR
    \STATE {\tt //* {end streaming upon time-out or interruption} *//}
    \IF{gossip is successful}
        \STATE Update ${\cal B}_i \leftarrow \{j\}$, ${\cal B}_j \leftarrow \{i\}$.
    \ENDIF
    \ENDIF
\ENDFOR
\end{algorithmic}
\end{algorithm}

\begin{algorithm}
\caption{\textit{computation\_thread()} of Agent $i$} 
\begin{algorithmic}[1]\label{alg:comp_thread}
\STATE {\bf local variable:} buffer ${\cal B}_i$; counters $g_i$, $t_i$.
\STATE {\tt //* {begin compute SG} *//}

Compute $\nabla f_i(\prm_i^{t_i}; \xi_i^{t_i })$ and increment $g_i \leftarrow g_i + 1$. \label{line:SGcompute}

{\tt //* {end compute SG upon completion or interruption} *//}
\IF{communication buffer ${\cal B}_i = \emptyset$} \label{line:caseA}
    \STATE \label{line:local_grad_step} \ding{229} \emph{Local SG step:} with $\hat{c}_i = g_i / (t_i + 1)$, evaluate
    \begin{align}
        \prm_i^{{t_i}+1} &= \prm_i^{t_i} - \eta \widehat{\dprm}^{t_i}_i - \alpha \hat{c}_i \nabla f_i(\prm_i^{t_i}; \xi_i^{{t_i}}),
    \end{align}
    and set $\widehat{\dprm}^{{t_i}+1}_i = \widehat{\dprm}^{t_i}_i$,
        $t_i \leftarrow t_i + 1$.
\ELSE \label{line:caseB}
    \STATE Identify the quantities:
    \begin{align}
        & \textstyle t_i' = \max\{ t_i, ~\max_{ j\in {\cal B}_i} t_j  \}, \quad \textstyle d_i = 1 + t_i' - t_i, \quad {\bf C}_{ij}(\xi^{t_i'}) = {\rm BinDiag}({\cal I}_{ij}(\xi^{t_i'})), \\
        & \hat{c}_i = \begin{cases}
            g_i / (t_i' + 1) & \text{if $\nabla f_i(\prm_i^{t_i}; \xi_i^{t_i'})$ is ready}, \\
            0& \text{otherwise}.
        \end{cases}
    \end{align}
    \STATE \ding{229} \emph{Gossip with SG step}: evaluate
    \begin{align}
       \prm_i^{t_i' + 1} & \textstyle = \prm_i^{t_i} - \gamma \sum_{j \in {\cal B}_i} {\bf C}_{ij}(\xi^{t_i'}) (\prm_i^{t_i} - \prm_j^{t_j}) - d_i \eta \widehat{\dprm}^{t_i}_i - \alpha \hat{c}_i \nabla f_i(\prm_i^{t_i}; \xi_i^{t'_i}), \\
       \widehat{\dprm}_i^{t_i'+1} &= \textstyle \widehat{\dprm}_i^{t_i} + \beta \sum_{j \in {\cal B}_i} {\bf C}_{ij}(\xi^{t'_i}) (\prm_i^t - \prm_j^t),
    \end{align}
    and set $t_i \leftarrow t_i' + 1$, ${\cal B}_i \leftarrow \emptyset$.
\ENDIF
\end{algorithmic}
\end{algorithm}

\begin{figure}[h]
    \centering
    \includegraphics[width=0.5\textwidth]{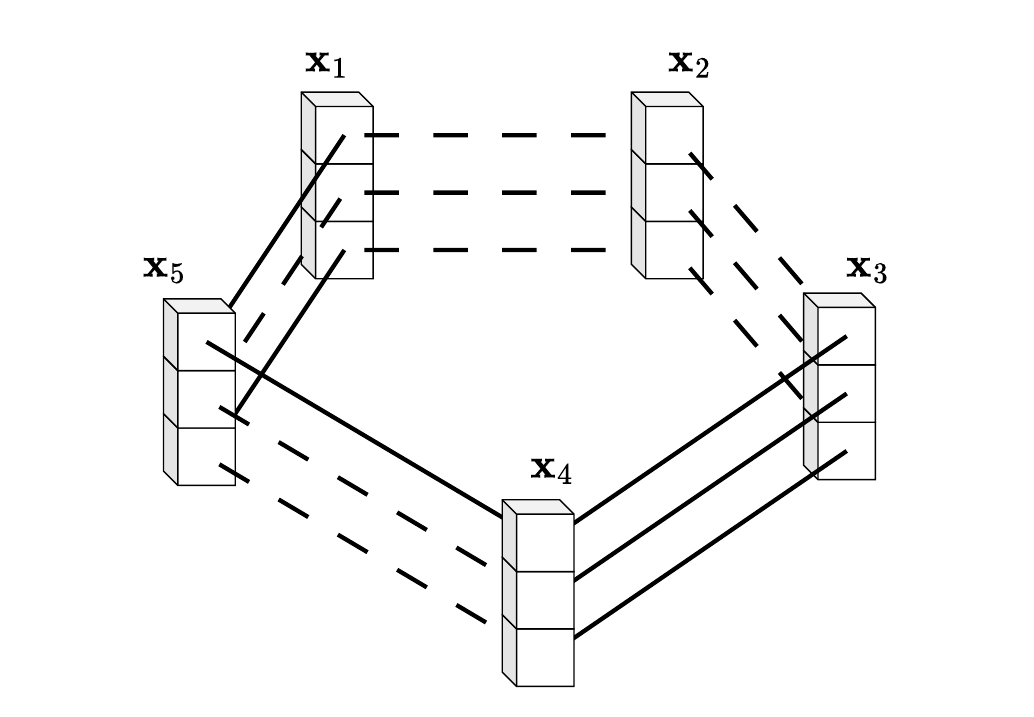}
    \caption{Illustration of a (time-varying) random graph ${\cal G}(\xi)$ for primal variable of dimension $d=3$ on a ring network of $n = 5$ nodes. Solid lines represent active edges while dashed lines represent disconnected edges. In this example, node 2 is considered as idle in an asynchronous environment. ${\bf C}_{15}(\xi)$ is a diagonal matrix such that ${\rm diag}({\bf C}_{15}(\xi)) = (1, 0, 1)$.}
    \label{fig:extended_graph}
\end{figure}

\newpage
\section{Detailed Convergence Rate Analysis for {\algnamesa}}
\label{app:detail_analysis}

Under the parameter choices of $\a = \frac{80\delta_1 L^2 n }{\sqrt{T \sigavg} \gamma_{\infty} \rho_{\min}}$, $\alpha = \sqrt{n / (T \sigavg)}$, we list several remarks on the convergence rate analysis of {\algnamesa} as follows.

\begin{itemize}[leftmargin=*]
\item In addition to the crude bound in \eqref{eq:crude_bound}, we can derive a fine grained characterization for the convergence rate of {\algnamesa}. Here, an interesting aspect is bounding the \emph{transient time} of the number of iterations for the decentralized algorithm to achieve the rate on par with CSGD in \eqref{eq:crude_bound}, independent of the network topology \citep{pu2021sharp}. Using the definition of $\bbC_{\sigma}$ in \eqref{eq:lyapunov_main}, it can be shown that {\algnamesa} has a transient time of
\begin{align}
T_{\rm trans} = \Omega \Big( \frac{\sigma_A^4}{\rho_{\min}^4} \cdot &\max\Big\{ n^6 \rho_{\max}^2, \min\{ \frac{\bar{\rho}_{\max}^4 \rho_{\max}^6 }{n \sigavg \rho_{\min}^3 } , \frac{n^{5/2} \bar{\rho}_{\max}^2 \rho_{\max}^4 }{\sigavg \rho_{\min}^2} \} \Big\} \Big)  
\end{align}
where we have hidden the dependence of $L$, $F(\avgprm^0) - f_\star$, $\| \prm^0 \|_{\wtk}^2$  in the $\Omega(\cdot)$ notation.
\item Theorem \ref{thm:main} also guarantees that  $\mathbb{E}[\| \prm^{\sf T}\|^2_{\wtk}] \to 0$ and $\mathbb{E}[\| {\bf v}^{\sf T}\|_{\wtk}^2] \to 0$ as $T \to \infty$. The latter ensures that the gradient tracking error \citep{lu2019gnsd} variable converges to zero and thus {\algnamesa} is stable at a global stationary solution (see more discussion around the definition of ${\bf v}$ in \eqref{eq:v_def_main}).
Notice that at each iteration, {\algnamesa} only communicates the local parameters $\prm^t$ once, while the gradient tracking algorithm \citep{lu2019gnsd} communicates twice for the gradient tracking vectors and the local parameters.
\item For sufficiently large $T$, the effect of noisy network only remains dominant in $\expec{\| \prm^{\sf T} \|_{\wtk}^2 } = \mathcal{O}(n^2 \sigma_A^2 \rho_{\max} / (T \rho_{\min}^2 ))$, keeping the rate in \eqref{eq:crude_bound} unaffected. When comparing against synchronous {\tt DSGD} \citep{lian2017can}, with $b_i(\xi) = 1$ for all $i\in [n]$, {\algnamesa} will converge faster in terms of stationarity during the post-transient stage due to the shorter iteration time under coordinate sparsification and random graph.
\item When $\bar{b}_i \neq \bar{b}_j$, i.e., when the local update rates are non-uniform, the asynchronous gradient induces an error that will be reflected on the convergence error $\mathcal{O}(\sigma / \sqrt{nT} )$ through $\sigma$. We remark that under a fully controlled environment, {\algnamesa} can accelerate through asynchronous gradient only when the iteration time speedup out-weights the induced variance error.


\end{itemize}
\vspace{.2cm}

\section{Proof of Theorem~\ref{thm:main}} \label{sec:proof}
As our goal is to develop bounds for $\| \nabla F( \avgprm^t ) \|^2$, a natural idea is to consider the descent of the primal objective value $F( \avgprm^t )$. This yields:
\begin{lemma} \label{lemma:descent}
    Under Assumption \ref{assm:lip} and \ref{assm:f_var}, and the step size condition $\alpha \leq \frac{1}{4L}$,
    \begin{align} 
        &\mathbb{E}_t\left[  F(\avgprm^{t+1}) \right] \leq F (\avgprm^{t}) - \frac{\alpha}{4} \left\| \nabla F(\avgprm^t) \right\|^2  + \frac{3\alpha L^2}{4n}  \| \prm^t \|^2_{\wtk} + \frac{\alpha^2 L}{2n} \sigavg. \label{eq:descent_lemma}
    \end{align}
\end{lemma}
\noindent See Appendix \ref{app:lemma_descent_proof} for the proof. 
It can be seen from the above lemma that if the consensus error satisfies $\| \prm^t \|_{\wtk}^2 = o(1)$, then setting $\alpha = 1 / \sqrt{T}$ suffices to yield 
\begin{equation}
\textstyle
(1/T) \sum_{t=1}^T \mathbb{E} [ \left\| \nabla F(\avgprm^t) \right\|^2 ] = {\cal O}( 1 / \sqrt{T} ).
\end{equation}
That said, the evolution of $\| \prm^t \|_{\wtk}^2$ tends to be complicated as the latter co-evolves with the dual variable $\dprm^t$ and the gradient. 
To simplify, we impose the following preliminary condition on the step sizes
\begin{equation} \label{eq:stepsizecond-consensus}
    \gamma \leq \min \left\{ \frac{\rho_{\min}}{\rho_{\max}^2}, \frac{ \rho_{\min} }{ 2 \sigma_A^2 \rho_{\max} } \right\},~ \alpha \le 1,~\eta \leq 1.
    \end{equation}
\begin{lemma} \label{lemma:consensus} Under Assumptions \ref{assm:lip}, \ref{assm:rand-graph}, \ref{assm:f_var}, \ref{assm:graph_var} and the step size condition \eqref{eq:stepsizecond-consensus}. The consensus error follows the recursive inequality
    \begin{align}
    \expec{ \| \prm^{t+1} \|_{\wtk}^2} &\leq \left[ 1 - \frac{\gamma}{2} \rho_{\min} + \alpha (1 + 3 L^2)  \right] \expec{\left\| \prm^t \right\|_{ \wtk }^2 } + 2\eta^2 \expec{\| {\bf v}^t \|_{  \wtk}^2} \\
    &\quad -2\eta \expec{ \dotp{\prm^t}{ {\bf v}^t }_{\wtk - \gamma \wta^\top \rlamb \rxii \wta } } + \alpha^2 n \sigavg .
\end{align}
\end{lemma}
\noindent See Appendix \ref{app:lemma_consensus_proof} for the proof.

The above lemma shows an intricate structure for the consensus error as the latter also depends on the violation of \eqref{eq:dual_meaning}, i.e., $\| {\bf v}^t \|_{ \wtk} ^2$, and the \emph{weighted} inner product between $\prm^t, {\bf v}^t$.
Naturally, we can further control the above terms:
\begin{lemma} \label{lemma:dual_err}
    Under Assumption \ref{assm:lip}, \ref{assm:f_var}. Let $\alpha \leq 1$, then for any constant $\c > 0$, the dual error satisfies
    \begin{align}
        & \expec{ \left\| {\bf v}^{t+1} \right\|_{\wtq + \c\wtk}^2} \leq \expec{\left\| {\bf v}^t \right\|_{\wtq + \c \wtk}^2} + 2\alpha (\rho_{\min}^{-1} + \c) \expec{ \left\| {\bf v}^t \right\|_{\wtk}^2} + 2\beta \expec{\dotp{{\bf v}^t}{ \prm^t}_{\wtk + \c \wta^\top  \rlamb \rxi \wta}} \notag \\
        & + \left( 2\beta^2 \bar{\rho}_{\max}^2 + \frac{10 \alpha^3 L^4}{\eta^2} \right) \big( \rho_{\min}^{-1} + \c \big) \expec{ \| \prm^t \|_{\wtk}^2 } + \frac{5 \alpha^3}{\eta^2} L^2 \big( \rho_{\min}^{-1} + \c \big) \left\{ \sigavg  + 2 n \expec{ \| \nabla F( \avgprm^t ) \|^2 }  \right\} 
    \end{align}
\end{lemma}
\noindent See Appendix~\ref{app:lemma_dual_err_proof} for the proof. Notice that we have considered the weighted norm $\left\| {\bf v}^t \right\|_{\wtq + \c \wtk}^2$ to induce a favorable inner product term for $\prm^t, {\bf v}^t$ below:

\begin{lemma} \label{lemma:xv_inner} Under Assumption \ref{assm:lip}, \ref{assm:f_var}, and the step size condition \eqref{eq:stepsizecond-consensus}, then 
    \begin{align}
     \expec{\dotp{\prm^{t+1}}{\mathbf{v}^{t+1}}_{\wtk}} &\leq \frac{\alpha - \eta }{2} \expec{ \| {\bf v}^t \|^2_{\wtk}} + \expec{\dotp{\prm^t}{\mathbf{v}^t}_{ \wtk - \gamma \wta^\top \rlamb \rxii \wta - \eta\beta  \wta^\top \rlamb  \rxi \wta }} \notag \\
    & + \left\{ \beta\rho_{\max} + \frac{\alpha + 3 \alpha^2 L^2 + \alpha \gamma^2 \sigma_A^2 \rho_{\max} + 20 \alpha^3 L^4 }{2\eta} + \frac{ \eta \beta^2 \rho_{\max}^2 }{2}   \right\} \expec{\|\prm^{t}\|^2_{ \wtk }} \notag \\
    & + \frac{10 \alpha^3 L^2 + \alpha^3 n}{2\eta} \sigavg  + \frac{10  \alpha^3 L^2 n}{\eta} \expec{ \| \nabla F( \avgprm^t ) \|^2 }. \notag
\end{align} 
\end{lemma}
\noindent See Appendix~\ref{app:lemma_xv_inner_proof} for the proof.

The above lemma shows that the recursion for the inner product $\dotp{\prm^{t}}{\mathbf{v}^{t}}_{\wtk}$ also depends on $\| \prm^t \|_{\wtk}^2$, $\| {\bf v}^t \|_{\wtk}^2$, etc., and can be similarly controlled.
Our final step is to construct a potential function $F_t$ whose recursive relation can guide us towards the convergence of {\algname}:
\begin{tcolorbox}[boxsep=2pt,left=4pt,right=4pt,top=3pt,bottom=3pt]
\begin{theorem} \label{lem:abcd_cond}
    For some constants $\a,\b,\c,\d>0$, we define the potential function
\begin{equation} \label{eq:ft_def_restated}
    F_t = \expec{ F(\avgprm^{t})  +  \a\| \prm^t \|_{\wtk}^2 + \b \| {\bf v}^t \|_{\wtq + \c \wtk}^2 + \d\dotp{\prm^t}{ {\bf v}^t }_{\wtk} }.
\end{equation}


Then, by the following choice of hyperparameters
\begin{align}
\b = \a \cdot \frac{\eta}{\beta}, \quad \c = \frac{(\eta \beta + \gamma )\d - 2 \eta \gamma \a}{2 \beta \b}, \quad \d = \delta_1 \eta \a, \label{eq:abcd_choice}
\end{align}
\begin{equation}
    \alpha \leq \alpha_\infty, \quad \eta \leq \eta_\infty, \quad \gamma \leq \gamma_\infty, \quad \beta = 1, \quad \delta_1 \ge 8, \label{eq:ss_cond}
\end{equation}
where
\begin{align}
    & \mkern-12mu \gamma_\infty := \frac{\rho_{\min}}{\rho_{\max}^2} \min \left\{ 1, \frac{\rho_{\max}}{ 2\sigma_A^{2}} \right\},~\eta_\infty := \frac{\rho_{\min}^2 }{64 \delta_1^2 \bar{\rho}_{\max}^2 \rho_{\max}^2 } \, \gamma_\infty ,
    \\
    & \alpha_\infty := \frac{\gamma_\infty  \rho_{\min}}{ 80 \delta_1  \sqrt{n}} \min\left\{ \frac{\a}{L^2}, \quad \eta_{\infty} \rho_{\min}, \sqrt{\frac{\eta_\infty \rho_{\min}}{L^2 \a}} \right\},
\end{align}
it holds that $F_t \geq F (\avgprm^t) \geq f_\star > -\infty$ for any $t \geq 0$, and the potential function follows the inequality
\begin{equation}
\begin{aligned}
F_{t+1} & \leq F_t - \frac{\alpha}{8} \expec{ \| \nabla F( \avgprm^t ) \|^2 } - \frac{ \a \gamma \rho_{\min}}{8} \expec{ \| \prm^t \|_{\wtk}^2 } - \a \eta^2 \expec{ \| {\bf v}^t \|_{\wtk}^2 } + \bbC_{\sigma} \alpha^2 \sigavg / n \label{eq:potential_end}
\end{aligned}
\end{equation}
such that
\begin{align}
\bbC_{\sigma} \leq \frac{ L }{2} + \a n^2 + \a n \alpha L^2 \left(  \frac{5}{\eta \beta \rho_{\min}} + \frac{5 \delta_1 \gamma }{ \eta \beta} + 5 \delta_1 \right) + \frac{\a \delta_1 n^2 \alpha }{2} \label{eq:bbC_bound}
\end{align}
\noindent See Appendix~\ref{app:abcd_cond} for the proof.
\end{theorem}
\end{tcolorbox}

By summing up \eqref{eq:potential_end} from $t=0$ to $t=T-1$ and rearranging terms, we obtain the bounds that lead to the theorem.
\hfill $\square$

\subsection{Proof of Lemma \ref{lemma:descent}}  \label{app:lemma_descent_proof}
By Assumption \ref{assm:lip}, it is straightforward to derive that
\begin{align}
    \mathbb{E}_t \left[ F (\avgprm^{t+1}) \right] 
    & \leq F(\avgprm^{t}) - \frac{\alpha}{n}\dotp{ \nabla F(\avgprm^{t})}{ \oneotimes \nabla {\bf f}(\prm^t)} + \frac{\alpha^2 L}{2n^2} \mathbb{E}_t \left[ \left\| \oneotimes  \nabla {\bf f}(\prm^t; \xi^{t})\right\|^2 \right] \label{proof:lemma_desc_step1}
\end{align}
where we have used the shorthand notation $\oneotimes := {\bf 1}^\top \otimes {\bf I}$.
The second term of \eqref{proof:lemma_desc_step1} can be bounded as
\begin{equation} 
\begin{aligned}
    & \frac{\alpha}{n} \dotp{ \nabla F(\avgprm^{t})}{ \oneotimes \nabla {\bf f} (\prm^t)} 
    \geq \frac{\alpha}{2}\left\| \nabla F(\avgprm^{t}) \right\|^2 - \frac{\alpha L^2}{2n} \left\| \prm^t \right\|_{\wtk}^2 
\end{aligned}
\end{equation}
The third term of \eqref{proof:lemma_desc_step1} can be bounded as
\begin{align}
    & \mathbb{E}_t \left[ \left\| \oneotimes \nabla {\bf f}(\prm^t; \xi^{t})\right\|^2 \right] 
    = \mathbb{E}_t\left[\left\| \oneotimes \left( \nabla {\bf f}(\prm^t; \xi^{t}) - \nabla {\bf f}(\prm^t) \right)\right\|^2\right] + \left\| \oneotimes \nabla {\bf f} (\prm^t) \right\|^2 \notag
\end{align}
Notice that due to the independence of gradient noise, we have $ \mathbb{E}_t\left[\left\|  \oneotimes \left( \nabla {\bf f}(\prm^t; \xi^{t}) - \nabla {\bf f}(\prm^t) \right)\right\|^2\right] \leq n \sigavg$. Furthermore, 
\begin{equation}
    \begin{aligned}
    \left\| \oneotimes \nabla {\bf f} (\prm^t) \right\|^2 & 
    \leq 2 n L^2 \| \prm^t \|_{\wtk}^2 + 2 n^2 \| \nabla F(\avgprm^t ) \|^2.
    \end{aligned}
\end{equation}
Substituting the above into \eqref{proof:lemma_desc_step1} and setting the step size $\alpha \leq 1 / (4L)$ concludes the proof of the lemma. \hfill $\square$

\subsection{Proof of Lemma \ref{lemma:consensus}} \label{app:lemma_consensus_proof}
Let $\oneotimesT := {\bf 1} \otimes {\bf I}$ and $\ARA := \wta^\top \rlamb \rxii \wta$, we introduce the following quantities to facilitate the analysis:
\begin{align}
    \serr^t &:= \alpha \left(\nabla {\bf f}(\prm^t) - \nabla {\bf f}(\prm^t; \xi^{t}) \right) + \gamma \left( \ARA - \wta(\xi^{t})^\top \wta \right)\prm^t \notag \\
    \gerr^t &:= \alpha \left( \nabla {\bf f}( \oneotimesT \avgprm^t) - \nabla {\bf f}(\prm^t) \right) \notag 
\end{align}
We can simplify the primal update as
\begin{align}
    \prm^{t+1} &= \left({\bf I} - \gamma \ARA \right) \prm^t - \eta \wta^\top \rlamb \dprm^t - \alpha \nabla {\bf f}( \oneotimesT \avgprm^t) + \serr^t + \gerr^t \notag \\
    &= \left({\bf I} - \gamma \ARA \right) \prm^t - \eta {\bf v}^t + \serr^t + \gerr^t
    \label{proof:x_update}
\end{align}
By \eqref{proof:x_update}, the consensus error can be measured by
\begin{equation}
\begin{aligned}
    & \expec{ \| \prm^{t+1} \|_{\wtk}^2} = \expec{\left\| \left({\bf I} - \gamma \ARA \right) \prm^t - \eta {\bf v}^t+ \serr^t + \gerr^t \right\|_{\wtk}^2} \\
    &\stackrel{(i)}{=} \expec{ \left\| \left({\bf I} - \gamma \ARA \right) \prm^t - \eta {\bf v}^t + \gerr^t  \right\|_{\wtk}^2} + \expec{\| \serr^t \|_{\wtk}^2} \\
    &= \expec{\left\| \left({\bf I} - \gamma \ARA \right) \prm^t \right\|_{\wtk}^2 } + \expec{\| \eta {\bf v}^t - \gerr^t \|_{\wtk}^2 }  - 2\expec{ \dotp{\left({\bf I} - \gamma \ARA \right) \prm^t}{\eta {\bf v}^t - \gerr^t }_{ \wtk } } + \expec{\| \serr^t \|_{\wtk}^2} \label{proof:lem1_step1}
\end{aligned}
\end{equation}
where $(i)$ uses the independence of random variables $\mathbb{E}[\serr^t | \prm^t, \dprm^t] = {\bf 0}$.
Let $\gamma \leq \rho_{\min} / \rho_{\max}^2$, by Assumption \ref{assm:rand-graph},
the first term of \eqref{proof:lem1_step1} can be bounded by 
\begin{equation} \label{eq:lem1_step2}
    \begin{aligned}
    \expec{\left\| \left({\bf I} - \gamma \ARA \right) \prm^t \right\|_{\wtk}^2 } \leq \big( 1 - \gamma \rho_{\min} \big) \expec{ \| \prm^t \|_{\wtk}^2 }.
    \end{aligned}
\end{equation}
The second term of \eqref{proof:lem1_step1} can be bounded by
\begin{align}
    \expec{\| \eta {\bf v}^t - \gerr^t \|_{\wtk}^2 } \leq 2 \eta^2 \expec{\| {\bf v}^t \|_{\wtk}^2} + 2 \expec{\| \gerr^t \|_{\wtk}^2 },
\end{align}
and the third term of \eqref{proof:lem1_step1} can be bounded by
\begin{equation}
\begin{aligned}
    &- 2\expec{ \dotp{\left({\bf I} - \gamma \ARA  \right) \prm^t}{\eta {\bf v}^t - \gerr^t }_{\wtk} }\\
    &\leq -2\eta \expec{\dotp{\prm^t}{ {\bf v}^t }_{\left({\bf I} - \gamma \ARA \right){\wtk}} } + \alpha \expec{\left\|\left({\bf I} - \gamma \ARA \right) \prm^t \right\|_{\wtk}^2 } + \frac{1}{\alpha} \expec{ \|\gerr^t \|_{\wtk}^2}
\end{aligned}
\end{equation}
Gathering the above inequalities, we obtain the following recursion on the consensus error:
\begin{align}
    &\expec{ \| \prm^{t+1} \|_{\wtk}^2} \leq \expec{ (1+\alpha) (1 - \gamma \rho_{\min}) \left\| \prm^t \right\|_{ \wtk }^2  +  2\eta^2 \| {\bf v}^t \|_{\wtk}^2 } \notag \\
    & + \expec{\| \serr^t \|_{\wtk}^2 + (2 + \frac{1}{\alpha}) \| \gerr^t \|_{\wtk}^2  -2\eta \dotp{\prm^t}{ {\bf v}^t }_{\wtk - \gamma \ARA } }. \notag
\end{align} 
Now we tackle the error term $\| \serr^t \|^2$,
by the independence between $\xi_a$ and $\xi_1,...,\xi_n$, it yields
\begin{equation}
\begin{aligned}
    \expec{ \| \serr^t\|^2_{\wtk} } 
    & \textstyle \leq \alpha^2 \sum_{i=1}^n \sigma_i^2 + \gamma^2 \sigma_A^2 \rho_{\max} \expec{ \| \prm^t \|_{\wtk}^2}, \label{eq:e_s_bound}
\end{aligned}
\end{equation}
where we have applied Assumptions~\ref{assm:rand-graph}, \ref{assm:f_var}, \ref{assm:graph_var} to obtain the above property.
Moreover, 
\begin{equation} 
\begin{aligned}
    \expec{ \| \gerr^t \|^2_{\wtk} } 
    & \leq 
    \alpha^2 L^2 \expec{\|\prm^t \|^2_{ \wtk }}.
\end{aligned} \label{eq:e_g_bound}
\end{equation}
To simplify expression, we assume $\alpha \leq 1$. Combining the upper bounds of $\expec{\| \serr^t \|^2_{\wtk}}$ and $\expec{\| \gerr^t \|^2_{\wtk}}$ gives us
\begin{equation} \notag
\begin{aligned}
    \expec{ \| \prm^{t+1} \|_{\wtk}^2} & \leq \left[ (1+\alpha) ( 1 - \gamma \rho_{\min}) + \gamma^2 \sigma_A^2 \rho_{\max} + 3 \alpha L^2  \right] \expec{\left\| \prm^t \right\|_{ \wtk }^2 } \\
    & \quad +  2\eta^2 \expec{\| {\bf v}^t \|_{\wtk}^2} -2\eta \expec{ \dotp{\prm^t}{ {\bf v}^t }_{\wtk - \gamma \ARA } } + \alpha^2 \sum_{i=1}^n \sigma_i^2 
\end{aligned}
\end{equation}
Using the step size condition \eqref{eq:stepsizecond-consensus} to simplify the first term completes the proof.
\hfill $\square$ 

\subsection{Proof of Lemma \ref{lemma:dual_err}} \label{app:lemma_dual_err_proof} 
Let $\oneotimesT := {\bf 1} \otimes {\bf I}$ and $\ARA := \wta^\top \rlamb \rxii \wta$, we observe that ${\bf v}^{t+1}$ is updated through the recursion:
\begin{align}
    & {\bf v}^{t+1} = \wta^\top \rlamb \dprm^{t+1} + \frac{\alpha}{\eta} \nabla {\bf f}( \oneotimesT \avgprm^{t+1}) \label{proof:v_update} \\
    &= {\bf v}^t + \underbrace{ \beta \wta^\top \rlamb  \wta(\xi^{t})\prm^t + \frac{\alpha}{\eta} \nabla {\bf f}( \oneotimesT \avgprm^{t+1}) - \frac{\alpha}{\eta} \nabla {\bf f}( \oneotimesT \avgprm^{t}) }_{ =: \Delta {\bf v}^t }
    \notag
\end{align}
Therefore, 
\begin{align}
    \expec{ \left\| {\bf v}^{t+1} \right\|_{\wtq + \c {\wtk}}^2} &\leq \expec{\left\| {\bf v}^t \right\|_{\wtq + \c \wtk}^2} + 2\expec{\dotp{{\bf v}^t }{ \Delta {\bf v}^t }_{\wtq + \c\wtk}} + 2 \beta^2 \expec{ \| \wta^\top  \rlamb \wta(\xi^{t})\prm^t \|_{\wtq + \c\wtk}^2 } \notag \\
    &\quad + \frac{2\alpha^2}{\eta^2} \expec{ \left\| \nabla {\bf f}( \oneotimesT \avgprm^{t+1}) -  \nabla {\bf f}( \oneotimesT \avgprm^{t}) \right\|_{\wtq + \c\wtk}^2}  \label{proof:lem2_step1}
\end{align}
The second term of \eqref{proof:lem2_step1} can be simplified as
\begin{align}
    \expec{\dotp{ {\bf v}^t }{ \Delta {\bf v}^t }_{\wtq + \c\wtk}}
    &\leq \beta \expec{\dotp{ {\bf v}^t }{ \prm^t}_{ \wtk + \c \ARA } } + \alpha \expec{ \left\| {\bf v}^t \right\|_{\wtq + \c\wtk}^2 } \notag \\
    & \quad +  \frac{\alpha}{ 4 \eta^2} \expec{\left\| \nabla {\bf f}( \oneotimesT \avgprm^{t+1}) - \nabla {\bf f}( \oneotimesT \avgprm^{t}) \right\|_{\wtq + \c \wtk}^2} \label{proof:lem2_step3}
\end{align}
where we used $\mathbb{E} [ (\wtq + \c\wtk) \wta^\top \wta(\xi^t) ] = \wtk + \c \ARA$
and note that $\expec{ \left\| {\bf v}^t \right\|_{\wtq + \c\wtk}^2 } \leq (\rho_{\min}^{-1} + \c) \expec{ \left\| {\bf v}^t \right\|_{\wtk}^2 }$ due to Assumption~\ref{assm:rand-graph}. The third term of \eqref{proof:lem2_step1} can be simplified as
\begin{align}
\expec{ \| \wta^\top  \rlamb \wta(\xi^{t})\prm^t \|_{\wtq + \c\wtk}^2 } &\leq \big( \rho_{\min}^{-1} + \c \big)  \expec{ \| \wta^\top  \rlamb \wta(\xi^{t}) \prm^t \|_{\wtk}^2 } \notag \\
&  = \big( \rho_{\min}^{-1} + \c \big)  \expec{ \| \prm^t \|_{\wta(\xi^t)^\top  \wta \wta^\top  \wta(\xi^{t})}^2 } \notag \\
&  \leq \rho(\wta \wta^\top )\bar{\rho}_{\max} \big( \rho_{\min}^{-1} + \c \big) \expec{ \| \prm^t \|_{\wtk}^2 } \notag \\
& = \bar{\rho}_{\max}^2 \big( \rho_{\min}^{-1} + \c \big) \expec{ \| \prm^t \|_{\wtk}^2 }. \notag
\end{align}
The fourth term of \eqref{proof:lem2_step1} can be simplified using Lemma~\ref{lemma:update_err} as
\begin{align}
    &\expec{ \left\| \nabla {\bf f}( \oneotimesT \avgprm^{t+1}) -  \nabla {\bf f}( \oneotimesT  \avgprm^{t}) \right\|_{\wtq + \c\wtk}^2} \notag \\
    & \leq \frac{ 4\alpha^2 n L^2 } { ( \rho_{\min}^{-1} + \c)^{-1} } \Big\{  \frac{1}{2 n^2}  \sum_{i=1}^n  \sigma_i^2  + \expec{ \| \nabla F( \avgprm^t ) \|^2 } + \frac{L^2}{n} \expec{ \| \prm^t \|_{\wtk}^2 }  \Big\}  .
\end{align}
Combining the above inequalities and using the condition $\alpha \leq 1$ to simplify constants yields the lemma.
\hfill $\square$

\subsection{Proof of Lemma \ref{lemma:xv_inner}}\label{app:lemma_xv_inner_proof}
Let $\oneotimesT := {\bf 1} \otimes {\bf I}$, $\ARA := \wta^\top \rlamb \rxii \wta$, and denote:
\[ 
\Delta {\bf v}^t := \beta \ARA \prm^t + \frac{\alpha}{\eta} (\nabla {\bf f}( \oneotimesT \avgprm^{t+1}) - \nabla {\bf f}( \oneotimesT \avgprm^{t})),
\]
then by the recursions in \eqref{proof:x_update} and \eqref{proof:v_update},
\begin{equation} 
\begin{aligned}
    &\expec{\dotp{\prm^{t+1}}{\mathbf{v}^{t+1}}_{\wtk}} \\
    & = \expec{\dotp{\prm^t}{\mathbf{v}^t}_{ \wtk - ( \gamma + \eta \beta ) \ARA }} + \expec{\|\prm^{t}\|^2_{\beta \left({\bf I} - \gamma \ARA \right) \ARA }} - \eta \expec{\| \mathbf{v}^t \|^2_{\wtk} } \\
    & \quad +\frac{\alpha}{\eta} \expec{\dotp{\prm^t}{\nabla {\bf f}( \oneotimesT \avgprm^{t+1}) - \nabla {\bf f}( \oneotimesT  \avgprm^{t})}_{\left({\bf I} - \gamma \ARA \right)\wtk}} -\alpha \expec{\dotp{\mathbf{v}^{t}}{\nabla {\bf f}( \oneotimesT \avgprm^{t+1}) - \nabla {\bf f}( \oneotimesT \avgprm^{t}) }_{\wtk}} \\
    & \quad + \expec{\dotp{ \gerr^t}{\mathbf{v}^{t} + \Delta {\bf v}^t }
    _{\wtk}} + \frac{\alpha}{\eta} \expec{ \dotp{\serr^t}{ \nabla {\bf f}( \oneotimesT \avgprm^{t+1}) - \nabla {\bf f}( \oneotimesT \avgprm^{t}) }_{\wtk} },
\end{aligned}
\end{equation}
where we have used the fact $\mathbb{E}[\serr | \prm^t] = {\bf 0}$. We first note that $\expec{\|\prm^{t}\|^2_{\beta \left({\bf I} - \gamma \ARA \right) \ARA }} \leq \beta (\rho_{\max} - \gamma \rho_{\min}^2 ) \expec{ \| \prm^t \|_{\wtk}^2 }$. By the Young's inequality, we get
\begin{equation}
\begin{aligned}
    & \expec{\dotp{\prm^t}{\nabla {\bf f}( \oneotimesT \avgprm^{t+1}) - \nabla {\bf f}( \oneotimesT \avgprm^{t})}_{\left({\bf I} - \gamma \ARA \right)\wtk}} \\
    &\leq \frac{1}{2} (1 - \gamma \rho_{\min} ) \expec{ \| \prm^t \|_{ \wtk }^2 } + \frac{1}{2} \expec{\|\nabla {\bf f}( \oneotimesT \avgprm^{t+1}) - \nabla {\bf f}( \oneotimesT \avgprm^{t}) \|^2 }.
\end{aligned}
\end{equation}
Moreover,
\begin{equation} \notag
\begin{aligned}
    & \expec{\dotp{\mathbf{v}^{t}}{\nabla {\bf f}( \oneotimesT \avgprm^{t+1}) - \nabla {\bf f}( \oneotimesT \avgprm^{t}) }_{\wtk}} \\
    & \leq \frac{1}{2} \expec{ \| {\bf v}^t \|^2_{\wtk}} + \frac{1}{2} \expec{\|\nabla {\bf f}( \oneotimesT \avgprm^{t+1}) - \nabla {\bf f}( \oneotimesT \avgprm^{t}) \|^2 }.
\end{aligned}
\end{equation}
Next, it is derived in \eqref{eq:e_g_bound} that $\expec{ \| \gerr^t \|^2_{\wtk} } \leq \alpha^2 L^2 \expec{ \| \prm^t \|^2_{\wtk} }$. Thus,
\begin{equation}
\begin{aligned}
    &\expec{\dotp{ \gerr^t}{\mathbf{v}^{t} + \Delta {\bf v}^t }_{\wtk} } \\
    &\leq \frac{1}{2}(\frac{3}{\eta}) \expec{\| \gerr^t \|^2_{\wtk}} + \frac{1}{2}(\frac{\eta}{3}) \expec{\| \mathbf{v}^{t} + \Delta {\bf v}^t \|^2_{\wtk}} \\
    &\leq \frac{3}{2\eta} \expec{\| \gerr^t \|^2_{\wtk}} + \frac{\eta}{2} \expec{\| \mathbf{v}^{t} \|^2_{\wtk}} + \frac{\eta\beta^2}{2} \expec{\|\ARA \prm^t \|^2_{\wtk}}  + \frac{\eta}{2} \cdot \frac{\alpha^2}{\eta^2} \expec{\| \nabla {\bf f}( \oneotimesT \avgprm^{t+1}) - \nabla {\bf f}( \oneotimesT \avgprm^{t}) \|^2_{\wtk}} \\
    & \leq \left( \frac{3 \alpha^2 L^2 }{2\eta} + \frac{ \eta \beta^2 \rho_{\max}^2 }{2} \right) \expec{ \| \prm^t \|_{\wtk}^2 } + \frac{\eta}{2} \expec{\| \mathbf{v}^{t} \|^2_{\wtk}} + \frac{\eta}{2} \cdot \frac{\alpha^2}{\eta^2} \expec{\| \nabla {\bf f}( \oneotimesT \avgprm^{t+1}) - \nabla {\bf f}( \oneotimesT \avgprm^{t}) \|^2_{\wtk}}.
\end{aligned}
\end{equation}
Similarly, it is derived in \eqref{eq:e_s_bound} that $ \textstyle \expec{ \| \serr^t \|^2_{\wtk} } \leq \alpha^2 \sum_{i=1}^n \sigma_i^2 + \gamma^2 \sigma_A^2 \rho_{\max} \expec{ \| \prm^t \|_{\wtk}^2}$. Thus,
\begin{equation}
\begin{aligned}
    &\expec{ \dotp{\serr^t}{ \nabla {\bf f}( \oneotimesT \avgprm^{t+1}) - \nabla {\bf f}( \oneotimesT \avgprm^{t}) }_{\wtk} } \\
    &\leq \frac{1}{2} \expec{\| \serr^t \|^2_{\wtk}} + \frac{1}{2} \expec{\|  \nabla {\bf f}( \oneotimesT \avgprm^{t+1}) - \nabla {\bf f}( \oneotimesT \avgprm^{t}) \|^2_{\wtk}} \\
    &\leq \frac{\alpha^2}{2}\sum_{i=1}^n \sigma_i^2 + \frac{\gamma^2 \sigma_A^2 \rho_{\max}}{2} \expec{ \| \prm^t \|_{\wtk}^2}  + \frac{1}{2} \expec{\|  \nabla {\bf f}( \oneotimesT \avgprm^{t+1}) - \nabla {\bf f}( \oneotimesT \avgprm^{t}) \|^2_{\wtk}} .
\end{aligned}
\end{equation}
Therefore, using the condition $\eta \leq 1$ and simplifying terms yield
\begin{equation}
\begin{aligned}
    \expec{\dotp{\prm^{t+1}}{\mathbf{v}^{t+1}}_{\wtk}} &\leq 
    \expec{\dotp{\prm^t}{\mathbf{v}^t}_{ \wtk - ( \gamma + \eta \beta ) \ARA }}  + \frac{\alpha - \eta }{2} \expec{ \| {\bf v}^t \|^2_{\wtk}} + \frac{\alpha^3}{2\eta} \sum_{i=1}^n \sigma_i^2  \\
    & \quad+ \left\{ \beta \rho_{\max} + \frac{\alpha}{2\eta} + \frac{3 \alpha^2 L^2 }{2\eta} + \frac{ \eta \beta^2 \rho_{\max}^2 }{2} + \frac{\alpha \gamma^2 \sigma_A^2 \rho_{\max}}{2 \eta} \right\} \expec{\|\prm^{t}\|^2_{ \wtk }} \\
    & \quad+ (\frac{2\alpha}{\eta} + \frac{\alpha^2}{2\eta}) \expec{\|\nabla {\bf f}( \oneotimesT \avgprm^{t+1}) - \nabla {\bf f}( \oneotimesT \avgprm^{t}) \|^2 }.
\end{aligned}
\end{equation}
The proof is concluded by applying Lemma~\ref{lemma:update_err} to bound the last term.
\hfill $\square$

\subsection{Proof of Lemma \ref{lem:abcd_cond}}\label{app:abcd_cond}
By combining the results of Lemma \ref{lemma:descent}, \ref{lemma:consensus}, \ref{lemma:dual_err}, \ref{lemma:xv_inner}, we obtain
\begin{equation} \label{eq:lyapunov_main}
    \begin{aligned}
    F_{t+1} & \leq F_t - \bbC_{\nabla F} \expec{ \| \nabla F( \avgprm^t ) \|^2 } + \bbC_{\sigma} \alpha^2 \sigavg / n + \expec{  \bbC_{ \prm } \| \prm^t \|_{\wtk}^2  + \bbC_{{\bf v}} \| {\bf v}^t \|_{\wtk}^2 + \dotp{ \prm^t }{ {\bf v}^t }_{ \bC_{xv} } }
    \end{aligned}
\end{equation}
where
    \begin{align}
    \bbC_{\nabla F}  & := \frac{\alpha}{4} - \b \frac{10 \alpha^3 L^2 n}{\eta^2} ( \rho_{\min}^{-1} + \c ) - \d \frac{10 \alpha^3 L^2 n}{\eta} \\
    \bbC_{\sigma} & := \frac{ L }{2} + \a n^2 + \b \frac{5 n \alpha L^2}{\eta^2} ( \rho_{\min}^{-1} + \c )  + \d \left( \frac{5 n \alpha L^2}{\eta} + \frac{\alpha n^2 }{2\eta} \right) \\
     \bbC_{ \prm } & := \frac{3\alpha L^2}{4n }+ \a \left( - \frac{\gamma}{2} \rho_{\min} + \alpha (1 + 3 L^2) \right) + \b \left( 2\beta^2 \bar{\rho}_{\max}^2 + \frac{10 \alpha^3 L^4}{\eta^2} \right) \big( \rho_{\min}^{-1} + \c \big) \notag \\
    & \quad + \d \left( \beta\rho_{\max} + \frac{\alpha}{2\eta} + \frac{3 \alpha^2 L^2 }{2\eta} + \frac{\alpha \gamma^2 \sigma_A^2 \rho_{\max}}{2\eta}\right) + \d \left( \frac{10 \alpha^3 L^4 }{\eta} + \frac{ \eta \beta^2 \rho_{\max}^2 }{2} \right) \\
    \bbC_{ {\bf v} } & := \a (2 \eta^2) + \b (2 \alpha (\rho_{\min}^{-1} +  \c )) + \d \frac{ \alpha - \eta }{2} \\
\bC_{xv} & := \left( 2 \gamma \eta \a + 2 \beta \b \c - \d ( \gamma + \eta \beta ) \right) \wta^\top \rxii \wta + \left( -2 \eta \a + 2 \beta \b \right) \wtk.
\end{align}

The equality $\bC_{xv} = {\bf 0}$ is ensured by the parameter choices in \eqref{eq:abcd_choice}.
We then observe that for any $\delta_0 > 0$, it holds
\begin{align*}
    F_t &\ge F(\avgprm^{t})  + \left( \a - \frac{ \d }{2 \delta_0} \right) \| \prm^t \|_{\wtk}^2 +  \| {\bf v}^t \|_{\b\wtq + ( \b\c - \frac{\d \delta_0}{2} ) \wtk}^2 \\
    & \geq \frac{1}{n} \sum_{i=1}^n f_i(\avgprm^{t}) + \| {\bf v}^t \|_{(\b \rho_{\max}^{-1} + \b\c - \frac{\d^2}{4\a}) \wtk}^2,  \label{eq:potential_lb}
\end{align*}
where the second inequality is achieved by setting $\delta = \frac{ \d }{ 2 \a }$ and using Assumption~\ref{assm:graph_var}.
Together with 
\[
\b \rho_{\max}^{-1} + \b\c - \frac{\d^2}{4\a} =
(\frac{\eta \a }{\beta \rho_{\max}} - \frac{\d^2}{4\a}) + (\frac{\gamma \d}{2 \beta} - \frac{\eta \gamma \a}{\beta}) + \frac{\eta \d}{2} \ge 0
\]
which is due to the step size choices
\begin{equation} \label{eq:pf_alpha0}
\d = \delta_1 \eta \a, \quad \delta_1 \ge 8, \quad \eta \beta \leq \frac{4}{\delta_1^2 \rho_{\max}}.
\end{equation}
We obtain that $F_t \geq F( \avgprm^t ) > -\infty$.
Second, to upper bound $\bbC_{\bf v}$, under the additional condition $\eta \beta \leq \gamma$, and provided that 
\begin{align} \label{eq:pf_alpha1}
    \begin{cases}
        \frac{2\alpha \eta}{\beta \rho_{\min}} \leq \frac{\eta^2}{3} & \Leftrightarrow \alpha \leq \frac{\eta \beta \rho_{\min}}{6} \\
        \frac{2 \delta_1 \alpha \eta \gamma}{\beta} \leq \frac{\eta^2}{3} &\Leftrightarrow \alpha \leq \frac{\eta \beta}{6 \delta_1 \gamma} \\
        \frac{\delta_1 \alpha \eta}{2} \leq \frac{\eta^2}{3} &\Leftrightarrow \alpha \leq \frac{2\eta}{3\delta_1 }
    \end{cases}
\end{align}
We obtain that 
\begin{align}
    \bbC_{\bf v} 
    &\leq \a \eta \left( 2\eta +  \frac{2\alpha}{\beta \rho_{\min}}  + \frac{2\delta_1 \alpha \gamma}{\beta} + \frac{\delta_1 \alpha}{2} - \frac{\delta_1 \eta}{2} \right) \leq - \a \eta^2. \notag
\end{align}
Third, to upper bound $\bbC_{\bf x}$, under the conditions $\alpha \leq \frac{1}{2 L^2}$, $\eta \beta \leq \gamma$, $\eta \beta \leq \rho_{\max}^{-1}$, $\gamma \leq \sigma_A^{-1} \rho_{\max}^{-1/2}$, we have 
\begin{align}
    & \bbC_{\bf x} \leq \frac{3 \alpha L^2}{4 n} + \Big(  - \frac{\gamma}{2} \rho_{\min} + \alpha (1 + 3L^2 + 4 \delta_1) \notag \\
    & \qquad \qquad \qquad \quad + (2 \eta\beta \bar{\rho}_{\max}^2 + \frac{5 \alpha}{2 \eta \beta}) (\rho_{\min}^{-1} + \delta_1 \gamma)  + \frac{3}{2} \delta_1 \eta \beta \rho_{\max}  \Big) \, \a \notag
\end{align}
With the step size conditions:
\begin{align} \label{eq:pf_alpha2}
    \begin{cases}
        \alpha (1 + 3L^2 + 4 \delta_1 ) \leq \frac{\gamma}{16} \rho_{\min} &\Leftrightarrow \alpha \leq \frac{\gamma}{16 (1 + 3L^2 + 4\delta_1 )} \rho_{\min} \\
        2\eta \beta \bar{\rho}_{\max}^2 \rho_{\min}^{-1} \leq \frac{\gamma}{32} \rho_{\min} &\Leftrightarrow \eta \beta \leq \frac{\gamma}{64} \cdot \frac{\rho_{\min}^2}{\bar{\rho}_{\max}^2}\\
        2 \delta_1 \gamma \eta \beta \bar{\rho}_{\max}^2 \rho_{\min}^{-1} \leq \frac{\gamma}{32}{\rho_{\min}} &\Leftrightarrow \eta \beta \leq \frac{1}{64 \delta_1 }\cdot \frac{\rho_{\min}^2}{\bar{\rho}_{\max}^2} \\
        \frac{5 \alpha}{2\eta \beta \rho_{\min}} \leq \frac{\gamma}{32} \rho_{\min} &\Leftrightarrow \alpha \leq \frac{\gamma \eta \beta}{80} \rho_{\min}^2\\
        \frac{5\delta_1 \gamma \alpha}{2\eta \beta} \leq \frac{\gamma}{32}\rho_{\min} &\Leftrightarrow \alpha \leq \frac{\eta \beta}{80 \delta_1 } \rho_{\min} \\
        \frac{3}{2}\delta_1 \eta \beta \rho_{\max} \leq \frac{\gamma}{16} \rho_{\min} &\Leftrightarrow \eta \beta \leq \frac{\gamma}{24 \delta_1} \cdot \frac{\rho_{\min}}{\rho_{\max}}
    \end{cases}
\end{align}
and $\alpha \leq \frac{ \a \gamma n \rho_{\min} }{ 6 L^2}$, it is guaranteed that 
\begin{equation}
    \bbC_{\bf x} \leq 3\alpha L^2 / (4n) - \a \gamma \rho_{\min} / 4 \leq - \a \gamma \rho_{\min} / 8.
\end{equation}
Fourth, to lower bound $\bbC_{ \nabla F }$, under the condition $\eta \beta \leq \gamma$ and the step size conditions 
\begin{align} \label{eq:pf_alpha3}
    \begin{cases}
        \a \frac{10 \alpha^3 L^2 n}{ \eta \beta \rho_{\min}} \leq \frac{\alpha}{24} & \Leftrightarrow \alpha \leq \sqrt{ \frac{\eta \beta \rho_{\min}}{240 L^2 n \a } }\\
        \a \frac{10 \delta_1 \alpha^3 \gamma L^2 n}{\eta \beta} \leq \frac{\alpha}{24} & \Leftrightarrow \alpha \leq \sqrt{\frac{\eta \beta}{240 \delta_1 \gamma L^2 n \a}} \\
        \a \cdot 10 \delta_1 \alpha^3 L^2 n \leq \frac{\alpha}{24} & \Leftrightarrow \alpha \leq \sqrt{\frac{1}{240 \delta_1 L^2 n \a }}
    \end{cases}
\end{align}
we have 
\begin{align}
    \bbC_{ \nabla F } & \geq \frac{\alpha}{4} - \a \left( \frac{10 \alpha^3 L^2 n}{ \eta \beta \rho_{\min}} + \frac{5 \alpha^3 L^2 n}{\eta^2 \beta} \left( 2\delta_1 \gamma \eta \right) + 10 \delta_1 \alpha^3 L^2 n \right) \geq \frac{\alpha}{8}. \notag
\end{align}
Finally, the condition $\eta \beta \leq \gamma$ guarantees:
\begin{equation}
\bbC_{\sigma} \leq \frac{ L }{2} + \a n^2 + \a n \alpha L^2 \left(  \frac{5}{\eta \beta \rho_{\min}} + \frac{5 \delta_1 \gamma }{ \eta \beta} + 5 \delta_1 \right) + \frac{\a \delta_1 n^2 \alpha }{2}.
\end{equation}
This shows the desired properties regarding $\bbC_{\prm}, \bbC_{\nabla F}, \bbC_{\sigma}, \bbC_{\bf v}$.
Gathering the step size conditions, i.e., \eqref{eq:pf_alpha0}, \eqref{eq:pf_alpha1}, \eqref{eq:pf_alpha2}, \eqref{eq:pf_alpha3} and simplifying yields \eqref{eq:ss_cond}.
\hfill $\square$

\subsection{Proof of Corollary \ref{cor:pl_main}}
\label{app:pl_proof}
The proof can be seen as a direct extension of \eqref{eq:potential_end} with Assumption \ref{assm:pl_main}. For instance, from \eqref{eq:potential_end} we get
\begin{equation}
\begin{aligned}
F_{t+1} - f_\star &\leq F_t - f_\star - \frac{\alpha}{8} \expec{ \| \nabla F( \avgprm^t ) \|^2 } + \bbC_{\sigma} \alpha^2 \sigavg  - \frac{ \a \gamma \rho_{\min}}{8} \expec{ \| \prm^t \|_{\wtk}^2 } - \a \eta^2 \expec{ \| {\bf v}^t \|_{\wtk}^2 } \\
&= (1-\delta)(F_t - f_\star) + \delta ( \expec{ F(\avgprm^t) } - f_\star) - \frac{\alpha}{8} \expec{ \| \nabla F( \avgprm^t ) \|^2 }  + \bbC_{\sigma} \alpha^2 \sigavg \\
&\quad + (\delta \a - \frac{\a \gamma \rho_{\min}}{8}) \expec{ \| \prm^t \|_{\wtk}^2 } + \delta \b \expec{ \| {\bf v}^t \|_{\wtq}^2 } + (\delta \b\c - \a \eta^2) \expec{ \| {\bf v}^t \|_{\wtk}^2 } \\
&\quad + \delta \d\expec{\dotp{\prm^t}{ {\bf v}^t }_{\wtk}},
\end{aligned}
\end{equation}
for any $\delta > 0$.
Applying Assumption \ref{assm:pl_main}, $\wtq \preceq \rho_{\min}^{-1} \wtk$ and $\dotp{\prm^t}{ {\bf v}^t }_{\wtk} \leq \frac{1}{2} \| \prm^t \|^2_\wtk + \frac{1}{2} \| {\bf v}^t \|^2_\wtk$ gives us
\begin{equation}
    \begin{aligned}
        F_{t+1} - f_\star &\leq (1-\delta)(F_t - f_\star) + ( \frac{\delta}{2\mu} - \frac{\alpha}{8} ) \expec{ \| \nabla F( \avgprm^t ) \|^2 } + \bbC_{\sigma} \alpha^2 \sigavg \\
        &\quad  + \left[\delta(\a + \frac{\d}{2}) - \frac{\a \gamma \rho_{\min}}{8}\right] \expec{ \| \prm^t \|_{\wtk}^2 } + \left[\delta(\frac{\b}{\rho_{\min}} + \b\c + \frac{\d}{2}) - \a\eta^2\right] \expec{\| {\bf v}^t \|^2_\wtk}.
    \end{aligned}
\end{equation}
Then choosing $\delta > 0$ such that $\delta \leq \min \{ {\alpha \mu} / {4}, ~{\gamma \rho_{\min}} / {16},~ {\eta \beta} / ({3 \rho_{\min}}), {\eta} / {12} \}$ enforces the coefficients of excessive terms to be non-positive, thus give rise to \eqref{eq:pl_ineq_main}.
\hfill $\square$ \\

\subsection{Auxiliary Lemma}
\begin{lemma} \label{lemma:update_err} 
Under Assumption \ref{assm:lip} and \ref{assm:f_var}, 
\begin{equation} \notag
\begin{aligned} 
& \expec{ \left\| \nabla {\bf f}( \oneotimesT \avgprm^{t+1}) - \nabla {\bf f}( \oneotimesT \avgprm^{t}) \right\|^2 }\leq 4 \alpha^2 n L^2 \left\{  \frac{1}{2n^2}  \sum_{i=1}^n  \sigma_i^2  + \expec{ \| \nabla F( \avgprm^t ) \|^2 } + \frac{L^2}{n} \| \prm^t \|_{\wtk}^2  \right\},
\end{aligned}
\end{equation}
where we have denoted $\oneotimesT := {\bf 1} \otimes {\bf I}$.
\end{lemma}

{\it Proof of Lemma \ref{lemma:update_err}}. By the Lipschitz gradient assumption on each local objective function $f_i$,
\begin{equation} \notag
\begin{aligned}
    &\expec{ \left\| \nabla {\bf f}( \oneotimesT \avgprm^{t+1}) - \nabla {\bf f}( \oneotimesT \avgprm^{t}) \right\|^2 } \leq n L^2 \expec{\|\avgprm^{t+1} - \avgprm^t\|^2} = \frac{ \alpha^2 L^2}{n}  \expec{ \left\| \oneotimes \nabla {\bf f}(\prm^t; \xi^{t})  \right\|^2} .
\end{aligned} 
\end{equation}
The latter can be further expanded as 
\begin{equation} \notag
\begin{aligned}
    & \frac{\alpha^2 L^2}{n} \expec{ \left\| \oneotimes \nabla {\bf f}(\prm^t; \xi^{t}) - \oneotimes \nabla {\bf f}(\prm^t) +  \oneotimes \nabla {\bf f}(\prm^t) \right\|^2} \\
    & \leq \frac{2 \alpha^2 L^2}{n} \left\{ \sum_{i=1}^n  \sigma_i^2  + \expec{ \left\|  \oneotimes \nabla {\bf f}(\prm^t) \right\|^2}  \right\} .
\end{aligned} 
\end{equation}
Lastly, we note that 
\begin{equation} \notag
    \begin{aligned}
    & \expec{ \left\|  \oneotimes \nabla {\bf f}(\prm^t) \right\|^2} \\
    & \leq 2 n^2 \expec{ \| \nabla F( \avgprm^t ) \|^2 } + 2 n \sum_{i=1}^n \expec{ \left\| \nabla f_i( \prm_i^t ) - \nabla f_i( \avgprm^t ) \right\|^2 } \\
    & \leq 2 n^2 \expec{ \| \nabla F( \avgprm^t ) \|^2 } + 2 n L^2 \expec{ \| \prm^t \|_{\wtk}^2}.
    \end{aligned} 
\end{equation}
This completes the proof.
\hfill $\square$ \\

\section{Proof of Theorem \ref{thm:main_vr}} \label{app:thm_main_vr_proof}
The proof structure of Theorem \ref{thm:main_vr} is similar to that of Theorem \ref{thm:main}. Below we summarize the Lemmas that control the error quantities governing the convergence of {\algnamevr}.

{\bf Notations.} We introduce following shorthand notations to indicate the additional error factors.
\begin{align}
    \nabla_{\prm} {\cal L}(\prm, \dprm; \xi) &:= \nabla {\bf f}(\prm; \xi) + \frac{\eta}{\alpha} \wta^\top \dprm + \frac{\gamma }{\alpha} \wta^\top \wta(\xi) \prm \\
    \nabla_{\prm} {\cal L}(\prm, \dprm) &:= \nabla {\bf f}(\prm) + \frac{\eta}{\alpha} \wta^\top \dprm + \frac{\gamma }{\alpha} \wta^\top {\bf R} \wta \prm \\
    \avgmom^t &:= \frac{1}{n} \oneotimes \mom^t
\end{align}

We recall that {\algnamevr} can be described by the following system:
\begin{align}
\prm^{t+1} &= \prm^t - \alpha \mom^t \\
\dprm^{t+1} & = \dprm^t + \beta \momdual^t \\
\mom^{t+1} &= \nabla_{\prm} {\cal L}(\prm^{t+1}, \dprm^{t+1}; \xi^{t+1}) + (1 - a_x) (\mom^t - \nabla_{\prm} {\cal L}(\prm^t, \dprm^{t}; \xi^{t+1}) ) \\
\momdual^{t+1} &= \wta(\xi^{t+1})\prm^{t+1} + (1 - a_\lambda) (\momdual^t - \wta(\xi^{t+1}) \prm^t)
\end{align}
\begin{lemma} \label{lemma:descent_storm}
    Under Assumption \ref{assm:lip} and the step size condition $\alpha \leq 1/(2L)$,
    \begin{align} 
        \mathbb{E}_t\left[  F(\avgprm^{t+1}) \right] &\leq F (\avgprm^{t}) - \frac{\alpha}{4} \left\| \avgmom^t \right\|^2 - \frac{\alpha}{2} \left\| \nabla F(\avgprm^t) \right\|^2 + \frac{\alpha L^2}{n}  \| \prm^t \|^2_{\wtk} + \alpha \left\| \avgmom^t - \frac{1}{n} \oneotimes \nabla {\bf f}(\prm^t) \right\|^2  \label{eq:descent_lemma_storm}
    \end{align}
\end{lemma}
See Appendix \ref{app:lemma_descent_storm_proof} for the proof.

Notice that as $\oneotimes \wta = {\bf 0}$, the network average primal STORM estimator $\avgmom^t$ evolve as
\begin{equation}
    \avgmom^{t+1} = \frac{1}{n} \oneotimes \nabla {\bf f}(\prm^{t+1}; \xi^{t+1}) + ( 1- a_x ) (\avgmom^t - \frac{1}{n} \oneotimes \nabla {\bf f}(\prm^t; \xi^{t+1}) ), \label{eq:net_avgmom}
\end{equation}
i.e., $\avgmom^t$ tracks the global objective function gradient $(1/n) \oneotimes \nabla {\bf f}(\prm^t)$. Meanwhile, the local estimator $\mom^t$ tracks the Lagrangian gradient $\nabla_{\prm} {\cal L}(\prm, \dprm) = \expec{\nabla_{\prm} {\cal L}(\prm, \dprm; \xi) }$.

With Assumption \ref{assm:sample_lip}, we are able to construct momentum estimation error bound where $a_x$ controls the variance noise as in Lemma \ref{lem:net_mom_err} and \ref{lem:mom_err}.

\begin{lemma} \label{lem:net_mom_err}
    Under Assumption \ref{assm:lip}, \ref{assm:rand-graph}, \ref{assm:f_var}, \ref{assm:sample_lip} and the step size condition $0 < a_x \leq 1$,
    \begin{align}
        &\expec{ \left\| \avgmom^{t+1} - \frac{1}{n} \oneotimes \nabla {\bf f}(\prm^{t+1}) \right\|^2} \leq (1-a_x )^2 \expec{ \left\| \avgmom^{t} - \frac{1}{n} \oneotimes \nabla {\bf f}(\prm^{t}) \right\|^2} \\ 
        &\quad + \frac{8 (\sampL^2 + L^2) \alpha^2 }{n} \expec{\| \mom^t - \nabla_{\prm} {\cal L}(\prm^t, \dprm^t) \|^2 } + 32(\sampL^2 + L^2) \alpha^2 \expec{\| \nabla F(\avgprm^t) \|^2} \\
        &\quad + \frac{32 (\sampL^2 + L^2) \eta^2}{n} \expec{\| {\bf v}^t \|^2_{\wtk}}  + \frac{32 (\sampL^2 + L^2) ( \alpha^2 L^2 + \gamma^2 \rho_{\max}^2)}{n} \expec{\| \prm^t \|_{\wtk}^2} + 2a_x^2 \sigavg \notag
    \end{align}
\end{lemma}

See Appendix \ref{app:lemma_net_mom_err_proof} for the proof.

\begin{lemma} \label{lem:mom_err}
Under Assumption \ref{assm:lip}, \ref{assm:rand-graph}, \ref{assm:f_var}, \ref{assm:graph_var}, \ref{assm:sample_lip} and the step size conditions $\alpha \leq \sqrt{a_x / (64 (L^2 + \sampL^2))}, \gamma \leq \sqrt{a_x / (64 (\rho_{\max}^2 + \bar{\rho}_{\max}^2))}$, $0 < a_x \leq 1$, $\gamma \leq 1/(8 \sqrt{a_x} \sigma_A^2)$,
\begin{align}
    &\expec{\left\| \mom^{t+1} - \nabla_{\prm} {\cal L}(\prm^{t+1}, \dprm^{t+1}) \right\|^2} \leq (1 - \frac{a_x}{4}) \expec{\left\| \mom^{t} - \nabla_{\prm} {\cal L}(\prm^{t}, \dprm^{t}) \right\|^2}    \\
    &\quad + 4a_x^2 n \sigavg + 64 \alpha^2 n G \cdot \expec{\| \nabla F(\avgprm^t) \|^2} + \left( \frac{16 a_x^2 \eta^2 \gamma^2 \sigma_A^2}{\alpha^2} + 64 \eta^2 G \right) \expec{\| {\bf v}^t \|_{\wtk}^2}  \\
    &\quad  + \left( \frac{20 a_x^2 \gamma^2 \sigma_A^2}{\alpha^2} + 64 (\alpha^2 L^2 + \gamma^2 \rho_{\max}^2) G  \right) \expec{\| \prm^t \|_{\wtk}^2}
\end{align}
where $G := \sampL^2 + L^2 + (\gamma^2 / \alpha^2) \bar{\rho}_{\max}^2 + (\gamma^2 / \alpha^2) \rho_{\max}^2 $.
\end{lemma}
See Appendix \ref{app:lemma_mom_err_proof} for the proof.

\begin{lemma} \label{lem:consensus_storm}
Under Assumption \ref{assm:lip}, \ref{assm:graph_var} and the step size condition $\alpha \leq \gamma \rho_{\min} \sqrt{1 - \gamma \rho_{\min}} / (\sqrt{32}L)$,
    \begin{align}
        &\expec{\| \prm^{t+1} \|^2_{\wtk}} \leq (1 - \frac{\gamma \rho_{\min}}{8}) \expec{ \|\prm^t\|_{\wtk}^2 } + \frac{\eta^2}{(1-\gamma \rho_{\min})(1- \gamma \rho_{\min}/2)} \expec{\| {\bf v}^t \|^2_{\wtk}} \\
        &\quad - \frac{2\eta(1 - \gamma\rho_{\min}/4)}{(1-\gamma \rho_{\min})} \dotp{\prm^t}{{\bf v}^t}_{\wtk - \gamma \wta^\top {\bf R} \wta } + \frac{2\alpha^2}{\gamma \rho_{\min}} \expec{\| \mom^t - \nabla_\prm {\cal L}(\prm^t, \dprm^t) \|^2}
    \end{align}
\end{lemma}
See Appendix \ref{app:lemma_consensus_storm_proof} for the proof.

\begin{lemma} \label{lem:xv_inner_storm}
    Under Assumption \ref{assm:lip}, \ref{assm:rand-graph} and the step size conditions $\alpha \leq \min\{ 1/\sqrt{12}, 1/(2L) \}$, $\eta \beta \leq 1/12$, $\beta \leq 1$,
    \begin{align}
        &\expec{\dotp{{\bf v}^{t+1}}{\prm^{t+1}}_{\wtk} }\leq \expec{ \dotp{{\bf v}^{t}}{\prm^{t}}_{\wtk - (\gamma + \eta \beta) \wta^\top {\bf R}\wta } } - \frac{\eta}{2} \expec{ \| {\bf v}^t \|_{\wtk}^2} \\
        &\quad  + \left( \frac{3\alpha^2 L^2}{\eta} +  2\beta \rho_{\max}^2 - \gamma \beta \rho_{\min}^2 + \frac{5 \beta}{2} + \frac{5\alpha^2}{2\eta} \right) \expec{\| \prm^t \|_{\wtk}^2}  \\
        &\quad + \left( (\frac{3}{\eta} + \frac{1}{2})\alpha^2 + 2 \alpha^2 \beta + \frac{2 \alpha^4}{\eta} \right) \expec{\| \mom^t - \nabla_\prm {\cal L}(\prm^t, \dprm^t) \|^2} \\
        &\quad + \frac{\beta}{2} \expec{\| \wta^\top \momdual^t - \wta^\top {\bf R}\wta \prm^t \|^2} + \frac{\alpha^2 L^2 n}{2\eta} \expec{\| \avgmom^t \|^2}
    \end{align}
\end{lemma}
See Appendix \ref{app:lemma_xv_inner_storm_proof} for the proof.

\begin{lemma} \label{lem:v_err_storm}
    Under Assumption \ref{assm:lip}, \ref{assm:graph_var} and the step size conditions $\alpha \leq \eta, \beta \leq 1$, for any constant $\b, \c > 0$,
    \begin{align}
        &\expec{\| {\bf v}^{t+1} \|_{\b{\bf Q} + \c \wtk}^2 } \leq (1+ \beta + \frac{\alpha}{\eta}) \expec{\| {\bf v}^{t} \|_{\b{\bf Q} + \c \wtk}^2 } + 2\beta \expec{\dotp{{\bf v}^t}{\prm^t}_{\b \wtk + \c \wta^\top {\bf R} \wta }} \\
        &\quad +  3 \beta^2 \rho_{\max}^2 (\b \rho_{\min}^{-1} + \c) \expec{\| \prm^t \|_{\wtk}^2 }  + 4 \beta  (\b \rho_{\min}^{-1} + \c) \expec{\| \wta^\top \momdual^t - \wta^\top {\bf R} \wta \prm^t \|^2}  \\
        &\quad + \frac{4 n \alpha^3 L^2}{\eta} (\b \rho_{\min}^{-1} + \c) \expec{\| \avgmom^t \|^2}
    \end{align}
\end{lemma}
See Appendix \ref{app:lemma_v_err_storm_proof} for the proof.

\begin{lemma} \label{lem:mom_dual_err}
Under Assumption \ref{assm:lip}, \ref{assm:rand-graph}, \ref{assm:graph_var} and the step size condition $a_\lambda \leq \sigma_A^{-1}$,
    \begin{align}
        &\expec{\| \wta^\top \momdual^{t+1} - \wta^\top {\bf R} \wta\prm^{t+1} \|^2} \leq (1 - a_\lambda)^2 \expec{\| \wta^\top \momdual^{t} - \wta^\top {\bf R} \wta\prm^t \|^2} \\
        &\quad + \left(10 a_\lambda^2 \sigma_A^2 + 32 (\rho_{\max}^2 + \bar{\rho}_{\max}^2 )(\alpha^2 L^2 + \gamma^2 \rho_{\max}^2 ) \right) \expec{\| \prm^{t} \|_{\wtk}^2} + 16 \eta^2 \expec{\| {\bf v}^t \|_{\wtk}^2} \\
        &\quad + 8 \alpha^2 ( 1+ \rho_{\max}^2 + \bar{\rho}_{\max}^2 ) \expec{\| \mom^t - \nabla_x {\cal L}(\prm^t, \dprm^t) \|^2} + 32 \alpha^2 n (\rho_{\max}^2 + \bar{\rho}_{\max}^2) \expec{\| \nabla F(\avgprm^t) \|^2}  \notag 
    \end{align}
\end{lemma}
See Appendix \ref{app:lemma_mom_dual_err_proof} for the proof.\footnote{Note that the proof can be easily extended to the case when $a_{\lambda} = 1$, i.e., {\algnamevr} without dual momentum, by ignoring the result of Lemma \ref{lem:mom_dual_err} and applying $\momdual^t = \wta^\top \wta(\xi^t) \prm^t$ together with Assumption \ref{assm:graph_var}.}

With the above Lemmas, we are ready to construct a potential function $F_t$ that balances the interaction between the error quantities:
\begin{tcolorbox}[boxsep=2pt,left=4pt,right=4pt,top=3pt,bottom=3pt]
\begin{theorem} \label{lem:potential_storm}
    For some constants $\a,\b,\c,\d,\e,\f,\g > 0$, we define the potential function
    \begin{align}
        F_t &= \mathbb{E} \Big[ F(\avgprm^t) + \a \| \prm^t \|_{\wtk}^2 + \| {\bf v}^t \|_{\b \wtq + \c \wtk}^2 + \d \dotp{\prm^t}{{\bf v}^t}_{\wtk} + \e \| \avgmom^t - \frac{1}{n} \oneotimes \nabla {\bf f}(\prm^t) \|^2 \\
        &\quad + \f \| \mom^t - \nabla_\prm {\cal L}(\prm^t, \dprm^t) \|^2 + \g \| \wta^\top \momdual^t - \wta^\top {\bf R} \wta \prm^t \|^2 \Big] \label{eq:potential_storm_def}
    \end{align}
    Then, by the following choice of hyperparameters
    \begin{align}
        &\a = n^{-1}, \quad \b = \a  \cdot \frac{\eta ( 1- \gamma \rho_{\min} / 4)}{\beta(1-\gamma \rho_{\min})} = \mathcal{O}(n^{-1} T^{2/3}), \label{eq:ss_choice_storm_start}\\
        &\c = \frac{\d}{2}\cdot (\frac{\gamma}{\beta} + \eta) - \a \cdot \frac{\eta \gamma}{\beta} (1- \gamma \rho_{\min} / 4) (1-\gamma \rho_{\min})^{-1} = \mathcal{O}(T^{2/3}), \\
        &\d = \mathcal{O}(T^{1/3}), \quad \e = \mathcal{O}(a_x^{-1/2}) = \mathcal{O}(\bar{\sigma}^{2/3} T^{1/3}), \\ 
        &\f = \frac{\a}{\gamma} = \mathcal{O}(n^{-1} T^{1/3}), \quad \g = \mathcal{O}(n^{-1} T^{1/3}), \\
        &\alpha = \mathcal{O}(\bar{\sigma}^{-2/3} T^{-1/3}), \quad \eta = \mathcal{O}(n), \\ 
        &\gamma = \mathcal{O}(T^{-1/3}), \quad \beta = \mathcal{O}( \frac{\a}{\d}\cdot \gamma ) = \mathcal{O}(n^{-1} T^{-2/3}), \\
        &a_x = \mathcal{O}(\bar{\sigma}^{-4/3} T^{-2/3}), \quad a_\lambda = \mathcal{O}( T^{-1/3} ), \label{eq:ss_choice_storm_end}
    \end{align}
    the potential function follows the inequality
    \begin{align}
        F_{t+1} &\leq F_t - \frac{\alpha}{4} \expec{\| \nabla F(\avgprm^t) \|^2} - \frac{\alpha}{8} \expec{\| \avgmom^t \|^2} + ( \e \cdot 2a_x^2 + \f \cdot 4 a_x^2 n) \sigavg \\
        &\quad - \a \cdot \frac{ \gamma \rho_{\min}}{8} \expec{\| \prm^t \|_{\wtk}^2} -\d \cdot \frac{ \eta}{4} \expec{\| {\bf v}^t \|^2_{\wtk}} - \e \cdot \frac{a_x}{2} \expec{\| \avgmom^t - \frac{1}{n} \oneotimes \nabla {\bf f}(\prm^t) \|^2} \\
        &\quad - \f \cdot \frac{a_x}{8} \expec{\| \mom^t - \nabla_\prm {\cal L}(\prm^t, \dprm^t) \|^2} - \g \cdot \frac{a_\lambda}{2} \expec{\| \wta^\top \momdual^t - \wta^\top {\bf R}\wta \prm^t \|^2} \label{eq:potential_conv}
    \end{align}
\end{theorem}
See Appendix \ref{app:lemma_potential_storm_proof} for the proof.
\end{tcolorbox}

Finally, summing up \eqref{eq:potential_conv} from $t=0$ to $t = T-1$ give us the convergence bound for the network average iterate and consensus error as
\begin{align}
    \frac{1}{T} \sum_{t=0}^{T-1} \expec{\| \nabla F(\avgprm^t) \|^2 } &\leq \frac{F_0 - F_T}{T \alpha / 4} + \frac{(\e \cdot 2a_x^2 + \f \cdot 4a_x^2 n) \sigavg}{\alpha / 4} = \mathcal{O}\left( \frac{F_0 - F_T + \bar{\sigma}^{2/3} }{T^{2/3}} \right) \\
    \frac{1}{T} \sum_{t=0}^{T-1} \expec{\| \prm^t \|_{\wtk}^2 } &\leq \frac{F_0 - F_T}{T \a \gamma \rho_{\min} / 8} + \frac{(\e \cdot 2a_x^2 + \f \cdot 4a_x^2 n) \sigavg}{\a \gamma \rho_{\min} / 8} = \mathcal{O}\left( \frac{F_0 - F_T + \bar{\sigma}^{2/3}  }{n^{-1} \rho_{\min} T^{2/3}} \right)
\end{align}
for large enough $T$.


\subsection{Proof of Lemma \ref{lemma:descent_storm}} \label{app:lemma_descent_storm_proof}
By Assumption \ref{assm:lip},
\begin{align}
    &\mathbb{E}_t[F(\avgprm^{t+1})] \\
    &= F(\avgprm^t) - \alpha \dotp{\nabla F(\avgprm^t) }{\avgmom^t} + \frac{\alpha^2 L}{2} \| \avgmom^t \|^2 \\
    &= F(\avgprm^t) - \left( \frac{\alpha}{2} + \frac{\alpha^2 L}{2} \right) \| \avgmom^t \|^2 - \frac{\alpha}{2} \| \nabla F(\avgprm^t) \|^2 + \frac{\alpha}{2} \| \avgmom^t - \nabla F(\avgprm^t) \|^2 \\
    &\stackrel{(i)}{\leq} F(\avgprm^t) -  \frac{\alpha}{4} \| \avgmom^t \|^2 - \frac{\alpha}{2} \| \nabla F(\avgprm^t) \|^2 + \alpha \left\| \avgmom^t - \frac{1}{n} \oneotimes \nabla {\bf f}(\prm^t) \right\|^2 + \alpha \left\| \frac{1}{n} \oneotimes \nabla {\bf f}(\prm^t) - \nabla F(\avgprm^t) \right\|^2 \notag
\end{align}
where $(i)$ uses the step size condition $\alpha \leq 1/(2L)$, and applying Assumption \ref{assm:lip} completes the proof. 

\hfill $\square$

\subsection{Proof of Lemma \ref{lem:net_mom_err}} \label{app:lemma_net_mom_err_proof}
By the network average momentum update in \eqref{eq:net_avgmom},
\begin{align}
    &\expec{\left\| \avgmom^{t+1} - \frac{1}{n} \oneotimes \nabla {\bf f}(\prm^{t+1}) \right\|^2} \\
    &=\mathbb{E}\big[ \big\| (1-a_x) (\avgmom^t - \frac{1}{n}\oneotimes \nabla {\bf f}(\prm^t)) + a_x \frac{1}{n}\oneotimes (\nabla {\bf f}(\prm^{t+1}; \xi^{t+1}) - \nabla {\bf f}(\prm^{t+1}) ) \\
    &\quad \quad + (1 - a_x) \frac{1}{n} \oneotimes (\nabla {\bf f}(\prm^t) - \nabla {\bf f}(\prm^t; \xi^{t+1}) - (\nabla {\bf f}(\prm^{t+1}) - \nabla {\bf f}(\prm^{t+1}; \xi^{t+1}))  \big\|^2 \big] \\
    &\leq (1 - a_x)^2 \expec{\| \avgmom^t - \frac{1}{n}\oneotimes \nabla {\bf f}(\prm^t) \|^2} + 2a_x^2 \expec{\left\| \frac{1}{n} \oneotimes (\nabla {\bf f}(\prm^{t+1}; \xi^{t+1}) - \nabla {\bf f}(\prm^{t+1}) )\right\|^2} \notag \\
    &\quad + 2(1-a_x)^2 \expec{\left\| \frac{1}{n} \oneotimes (\nabla {\bf f}(\prm^t) - \nabla {\bf f}(\prm^{t+1})) - (\nabla {\bf f}(\prm^t; \xi^{t+1})  - \nabla {\bf f}(\prm^{t+1}; \xi^{t+1}) ) \right\|^2}
\end{align}
Now observe that
\begin{align}
   2a_x^2 \expec{\left\| \frac{1}{n} \oneotimes (\nabla {\bf f}(\prm^{t+1}; \xi^{t+1}) - \nabla {\bf f}(\prm^{t+1}) )\right\|^2} \leq 2a_x^2 \sigavg ,
\end{align}
and 
\begin{align}
    &2(1-a_x)^2 \expec{\left\| \frac{1}{n} \oneotimes (\nabla {\bf f}(\prm^t) - \nabla {\bf f}(\prm^{t+1})) - (\nabla {\bf f}(\prm^t; \xi^{t+1})  - \nabla {\bf f}(\prm^{t+1}; \xi^{t+1}) ) \right\|^2} \\
    &\leq  4(1-a_x)^2 \expec{\left\| \frac{1}{n} \oneotimes (\nabla {\bf f}(\prm^t) - \nabla {\bf f}(\prm^{t+1})) \right\|^2} \\
    &\quad + 4(1-a_x)^2 \expec{\left\| \frac{1}{n} \oneotimes (\nabla {\bf f}(\prm^t; \xi^{t+1}) - \nabla {\bf f}(\prm^{t+1}; \xi^{t+1})) \right\|^2} \\
    &\leq \frac{4 (1-a_x)^2}{n} (L^2 + \sampL^2) \expec{\| \prm^{t+1} - \prm^t \|^2}
\end{align}
Finally, applying Lemma \ref{aux_lem:x_diff} and utilizing the step size condition $0 < a_x \leq 1$ to simplify the coefficients will complete the proof.
\hfill $\square$

\subsection{Proof of Lemma \ref{lem:mom_err}} \label{app:lemma_mom_err_proof}
\begin{align}
   &\expec{\left\| \mom^{t+1} - \nabla_{\prm} {\cal L}(\prm^{t+1},  \dprm^{t+1}) \right\|^2} \\
   &= \mathbb{E}\big[ \|(1-a_x) \left(\mom^{t} - \nabla_{\prm} {\cal L}(\prm^{t},  \dprm^{t}) \right) + a_x \left(\nabla_\prm {\cal L}(\prm^{t+1}, \dprm^{t+1}; \xi^{t+1}) - \nabla_\prm {\cal L}(\prm^{t+1}, \dprm^{t+1}) \right) \\
   &\qquad + (1-a_x) \left(\nabla_\prm {\cal L}(\prm^t, \dprm^t) - \nabla_\prm {\cal L}(\prm^t, \dprm^t, \xi^{t+1}) \right) \\
   &\qquad + (1-a_x) \left(\nabla_\prm {\cal L}(\prm^{t+1}, \dprm^{t+1}; \xi^{t+1}) - \nabla_\prm {\cal L}(\prm^{t+1}, \dprm^{t+1}) \right) \|^2 \big] \\
   &\leq (1 - a_x)^2 \expec{\| \mom^t - \nabla_\prm {\cal L}(\prm^t, \dprm^t) \|^2} \\
   &\quad + 2a_x^2 \expec{\|\nabla_\prm {\cal L}(\prm^{t+1}, \dprm^{t+1}; \xi^{t+1}) - \nabla_\prm {\cal L}(\prm^{t+1}, \dprm^{t+1}) \|^2} \\
   &\quad + 2(1-a_x)^2 \mathbb{E}\big[ \|\nabla_\prm {\cal L}(\prm^{t+1}, \dprm^{t+1}; \xi^{t+1}) - \nabla_\prm {\cal L}(\prm^t, \dprm^t, \xi^{t+1}) \\
   &\qquad \qquad \qquad \qquad  - ( \nabla_\prm {\cal L}(\prm^{t+1}, \dprm^{t+1}) - \nabla_\prm {\cal L}(\prm^t, \dprm^t)) \|^2 \big]
\end{align}
Now observe that
\begin{align}
    &2a_x^2 \expec{\|\nabla_\prm {\cal L}(\prm^{t+1}, \dprm^{t+1}; \xi^{t+1}) - \nabla_\prm {\cal L}(\prm^{t+1}, \dprm^{t+1}) \|^2} \\
    &= 2a_x^2 \expec{\left\| \nabla {\bf f}(\prm^{t+1}; \xi^{t+1}) - \nabla {\bf f}(\prm^{t+1}) + \frac{\gamma}{\alpha} \wta^\top \wta(\xi^{t+1})\prm^{t+1} - \frac{\gamma}{\alpha} \wta^\top {\bf R}\wta\prm^{t+1} \right\|^2 }\\
    &\leq 4 a_x^2 n \sigavg + \frac{4 a_x^2 \gamma^2 \sigma_A^2}{\alpha^2} \expec{\| \prm^{t+1} \|^2_{\wtk}} \\
    &\leq  4 a_x^2 n \sigavg + \frac{20 a_x^2 \gamma^2 \sigma_A^2}{\alpha^2}\expec{\| \prm^{t} \|^2_{\wtk}} + \frac{16 a_x^2 \eta^2 \gamma^2 \sigma_A^2}{\alpha^2} \expec{\|{\bf v}^t \|_{\wtk}^2} \\
    &\quad + 16 a_x^2 \gamma^2 \sigma_A^2 \expec{\| \mom^t - \nabla_\prm {\cal L}(\prm^t, \dprm^t) \|^2}
\end{align}
where the last inequality is expanded from the primal update rule $\prm^{t+1} = (1- \gamma \wta^\top {\bf R} \wta) \prm^t - \eta {\bf v}^t - \alpha(\nabla {\bf f}(\prm^t) - \nabla {\bf f}(\oneotimesT \avgprm^t)) - \alpha (\mom^t - \nabla_\prm {\cal L}(\prm^t, \dprm^t) )$ and the step size condition $\alpha \leq 1/(2L)$. On the other hand,
\begin{align}
    &\expec{ \|\nabla_\prm {\cal L}(\prm^{t+1}, \dprm^{t+1}; \xi^{t+1}) - \nabla_\prm {\cal L}(\prm^t, \dprm^t, \xi^{t+1}) - ( \nabla_\prm {\cal L}(\prm^{t+1}, \dprm^{t+1}) - \nabla_\prm {\cal L}(\prm^t, \dprm^t))\|^2} \notag \\
    &= \mathbb{E}\big[ \|\nabla {\bf f}(\prm^{t+1}; \xi^{t+1}) - \nabla {\bf f}(\prm^t; \xi^{t+1}) - (\nabla {\bf f}(\prm^{t+1}) - \nabla {\bf f}(\prm^t) )\\
    &\qquad \quad + \frac{\gamma}{\alpha} (\wta^\top \wta(\xi^{t+1})( \prm^{t+1} - \prm^t ) - \wta^\top {\bf R}\wta (\prm^{t+1} - \prm^t))\|^2 \big] \\
    &\leq 4 \left(\sampL^2 + L^2 + \frac{\gamma^2}{\alpha^2} \bar{\rho}_{\max}^2 + \frac{\gamma^2}{\alpha^2} \rho_{\max}^2 \right) \expec{\| \prm^{t+1} - \prm^t \|^2}
\end{align}
Finally, applying Lemma \ref{aux_lem:x_diff} and utilizing the step size conditions $\alpha \leq \sqrt{a_x / (64 (L^2 + \sampL^2))}$, $\gamma \leq \sqrt{a_x / (64 (\rho_{\max}^2 + \bar{\rho}_{\max}^2))}$, $0 < a_x \leq 1$, $\gamma \leq 1/(8 \sqrt{a_x} \sigma_A^2)$ to simplify the coefficients will complete the proof.
\hfill $\square$

\subsection{Proof of Lemma \ref{lem:consensus_storm}} \label{app:lemma_consensus_storm_proof}
\begin{align}
    &\expec{\| \prm^{t+1} \|_{\wtk}^2} = \expec{\| \prm^t - \alpha \mom^t \|_{\wtk}^2} \\
    &= \expec{\| \prm^t - \alpha \nabla_\prm {\cal L}(\prm^t, \dprm^t) - \alpha (\mom^t - \nabla_x {\cal L}(\prm^t, \dprm^t)) \|_{\wtk}^2 } \\
    &\leq (1 + z_1) \expec{\| \prm^t - \alpha \nabla {\bf f}(\prm^{t}) - \eta \wta^\top \dprm^t - \gamma \wta^\top {\bf R} \wta \prm^t \|_{\wtk}^2} \\
    &\quad + (1 + z_1^{-1}) \alpha^2 \expec{\|\mom^t - \nabla_x {\cal L}(\prm^t, \dprm^t)\|^2} \quad \text{for any } z_1 > 0,
\end{align}
Now observe that
\begin{align}
    &\expec{\| \prm^t - \alpha \nabla {\bf f}(\prm^{t}) - \eta \wta^\top \dprm^t - \gamma \wta^\top {\bf R} \wta \prm^t \|_{\wtk}^2} \\
    &= \expec{ \| ({\bf I} - \gamma \wta^\top {\bf R} \wta) \prm^t - \eta {\bf v}^t - \alpha (\nabla {\bf f}(\prm^t) - \nabla {\bf f}(\oneotimesT \avgprm^t)) \|_{\wtk}^2} \\
    &\leq (1+ z_2 ) \expec{\| ({\bf I} - \gamma \wta^\top {\bf R} \wta) \prm^t - \eta {\bf v}^t \|^2_{\wtk}} \\
    &\quad + ( 1 + z_2^{-1}) \alpha^2 \expec{\| \nabla {\bf f}(\prm^t) - \nabla {\bf f}(\oneotimesT \avgprm^t) \|^2 } \quad \text{for any } z_2 > 0, \\
    &\leq (1+ z_2)(1 - \gamma \rho_{\min}) \expec{\| \prm^t \|_{\wtk}^2} + (1+ z_2) \eta^2 \expec{\| {\bf v}^t \|_{\wtk}^2 } \\
    &\quad -2 (1 + z_2)\eta \dotp{\prm^t}{{\bf v}^t}_{\wtk - \gamma \wta^\top {\bf R}\wta} + (1 + z_2^{-1}) \alpha^2 L^2 \expec{\| \prm^{t} \|_{\wtk}^2}
\end{align}
Then, we choose $z_1 = \frac{\gamma \rho_{\min}/2}{1 - \gamma \rho_{\min}}$ so that $(1 + z_1)(1- \gamma \rho_{\min}) = 1 - \gamma \rho_{\min}/2$ and similarly $z_2 = \frac{\gamma \rho_{\min}/4}{1 - \gamma \rho_{\min}/2}$ so that $(1 + z_1)(1 + z_2)(1- \gamma \rho_{\min}) = 1 - \gamma \rho_{\min}/4$. Finally, observing that
\begin{align}
    (1+ z_1)(1 + z_2^{-1}) \alpha^2 L^2 \leq \frac{4 \alpha^2 L^2}{\gamma \rho_{\min} (1 - \gamma \rho_{\min})}  \leq \frac{\gamma \rho_{\min}}{8}
\end{align}
when imposing $\alpha \leq \gamma \rho_{\min} \sqrt{1 - \gamma \rho_{\min}} / (\sqrt{32} L)$ completes the proof.
\hfill $\square$

\subsection{Proof of Lemma \ref{lem:xv_inner_storm}} \label{app:lemma_xv_inner_storm_proof}
By the update rules, we have
\begin{align}
    {\bf v}^{t+1} &= {\bf v}^t + \beta \wta^\top {\bf R} \wta \prm^t + \beta \wta^\top(\momdual^t - {\bf R}\wta \prm^t) + \frac{\alpha}{\eta}(\nabla {\bf f}(\oneotimesT \avgprm^{t+1}) - \nabla {\bf f}(\oneotimesT \avgprm^{t})) \label{eq:v_storm_update_decomp} \\
    \prm^{t+1} &= (1- \gamma \wta^\top {\bf R} \wta) \prm^t - \eta {\bf v}^t - \alpha(\nabla {\bf f}(\prm^t) - \nabla {\bf f}(\oneotimesT \avgprm^t)) - \alpha (\mom^t - \nabla_\prm {\cal L}(\prm^t, \dprm^t) ) \label{eq:prm_storm_update_decomp}
\end{align}
Therefore, for any constants $z_1, z_2 > 0$,
\begin{align}
    &\expec{\dotp{{\bf v}^{t+1}}{\prm^{t+1}}_{\wtk}} \leq \expec{\dotp{{\bf v}^{t}}{\prm^{t}}_{\wtk - (\gamma + \eta \beta) \wta^\top {\bf R} \wta }} + (-1 + \frac{z_1 }{2} + \frac{z_2 }{2})  \eta \cdot   \expec{\| {\bf v}^t \|_{\wtk}^2 } \\
    &\quad + \left( \frac{\alpha^2 L^2}{2 z_1 \eta} + \beta \rho_{\max}^2 + \frac{\beta^2}{2} \rho_{\max}^2 + \frac{\alpha^2 L^2}{2} - \gamma \beta \rho_{\min}^2 + \frac{\beta^2}{2} \rho_{\max}^2 \right) \expec{\| \prm^t \|_{\wtk}^2} \\
    &\quad + \left(  \frac{\alpha^2 }{2 z_2 \eta} + \frac{\alpha^2}{2}\right) \expec{\| \mom^t - \nabla_\prm {\cal L}(\prm^t, \dprm^t) \|^2 }+ \frac{\beta}{2} \expec{\| \wta^\top \momdual^t - \wta^\top {\bf R}\wta \prm^t \|^2_{\wtk}} \\
    &\quad + \left( \frac{\beta }{2} + \frac{\alpha^2}{2\eta} \right) \expec{\| \prm^{t+1} \|_{\wtk}^2} + \frac{L^2 n }{2 \eta} \expec{ \| \avgprm^{t+1} - \avgprm^t \|^2 }
\end{align}
Note that by \eqref{eq:prm_storm_update_decomp}, we have
\begin{align}
    \|\prm^{t+1} \|_{\wtk}^2 \leq 4 (1-\gamma \rho_{\min})^2 \| \prm^t \|^2_{\wtk} + 4 \eta^2 \|{\bf v}^t \|_{\wtk}^2 + 4 \alpha^2 L^2 \|\prm^t \|_{\wtk}^2 + 4 \alpha^2 \| \mom^t - \nabla_\prm {\cal L}(\prm^t, \dprm^t) \|^2
\end{align}
Now choose $z_1 = z_2 = 1/6$, then applying $\| \avgprm^{t+1} - \avgprm^t \|^2 = \alpha^2 \| \avgmom^t \|^2$ 
and utilizing the step size conditions $\alpha \leq \min\{1/\sqrt{12}, 1/(2L) \}$, $\eta \beta \leq 1/12$, $\beta \leq 1$ to simplify the coefficients will complete the proof.
\hfill $\square$

\subsection{Proof of Lemma \ref{lem:v_err_storm}} \label{app:lemma_v_err_storm_proof}
By the dual update rule, we have
\begin{align}
    &\expec{\| {\bf v}^{t+1} \|_{\b {\bf Q} + \c \wtk}^2} \\
    &= \expec{\| {\bf v}^t + \beta \wta^\top {\bf R} \wta \prm^t + \beta \wta^\top (\momdual^t - {\bf R}\wta \prm^t) + \frac{\alpha}{\eta} (\nabla {\bf f}(\oneotimesT \avgprm^{t+1}) - \nabla {\bf f}(\oneotimesT \avgprm^t)) \|_{\b {\bf Q} + \c \wtk}^2 } \notag \\
    &\leq \expec{ \| {\bf v}^t \|_{\b {\bf Q} + \c \wtk}^2 } + \expec{\| {\bf v}^{t+1} - {\bf v}^t \|_{\b {\bf Q} + \c \wtk}^2} \\
    &\quad + 2 \expec{\dotp{{\bf v}^t}{\beta \wta^\top {\bf R} \wta \prm^t + \beta \wta^\top (\momdual^t - {\bf R}\wta \prm^t) + \frac{\alpha}{\eta} (\nabla {\bf f}(\oneotimesT \avgprm^{t+1}) - \nabla {\bf f}(\oneotimesT \avgprm^t))}_{\b {\bf Q} + \c \wtk}} \notag \\
    &\leq  (1+ \beta + \frac{\alpha}{\eta}) \expec{\| {\bf v}^{t} \|_{\b {\bf Q} + \c \wtk}^2 } + 2\beta \expec{\dotp{{\bf v}^t}{\prm^t}_{\b \wtk + \c \wta^\top {\bf R} \wta }} \\
    &\quad + (\frac{\alpha }{\eta} + \frac{3 \alpha^2}{\eta^2})  \expec{\| \nabla {\bf f}(\oneotimesT \avgprm^{t+1}) - \nabla {\bf f}(\oneotimesT \avgprm^t) \|^2_{\b {\bf Q} + \c \wtk}} + 3 \beta^2 \rho_{\max}^2 \expec{ \| \prm^t \|_{\b {\bf Q} + \c\wtk}^2 } \\
    &\quad + (\beta + 3\beta^2) \expec{\| \wta^\top \momdual^t - \wta^\top {\bf R}\wta \prm^t \|^2_{\b {\bf Q} + \c \wtk}} \\
    &\leq  (1+ \beta + \frac{\alpha}{\eta}) \expec{\| {\bf v}^{t} \|_{\b {\bf Q} + \c \wtk}^2 } + 2\beta \expec{\dotp{{\bf v}^t}{\prm^t}_{\b \wtk + \c \wta^\top {\bf R} \wta }} \\
    &\quad + \frac{4 n \alpha L^2}{\eta} (\b \rho_{\min}^{-1} + \c) \expec{ \| \avgprm^{t+1} - \avgprm^t \|^2} + 3 \beta^2 \rho_{\max}^2 (\b \rho_{\min}^{-1} + \c) \expec{ \| \prm^t \|_{\wtk}^2 } \\
    &\quad + 4\beta  (\b \rho_{\min}^{-1} + \c) \expec{\| \wta^\top \momdual^t - \wta^\top {\bf R}\wta \prm^t \|^2}
\end{align}
Then, applying $\| \avgprm^{t+1} - \avgprm^t \|^2 = \alpha^2 \| \avgmom^t \|^2$ 
and utilizing the step size conditions $\alpha \leq \eta$, $\beta \leq 1$ to simplify the coefficients will complete the proof.
\hfill $\square$

\subsection{Proof of Lemma \ref{lem:mom_dual_err}} \label{app:lemma_mom_dual_err_proof}
\begin{align}
    &\expec{\| \wta^\top \momdual^{t+1} - \wta^\top {\bf R} \wta \prm^{t+1} \|^2} \\
    & = \mathbb{E}\big[ \| 
    (1- a_\lambda) (\wta^\top \momdual^t - \wta^\top {\bf R} \wta \prm^t)  + a_\lambda (\wta^\top \wta(\xi^{t+1}) \prm^{t+1} - \wta^\top {\bf R}{\wta} \prm^{t+1}) \\
    &\qquad \quad + (1-a_\lambda) (\wta^\top {\bf R} \wta \prm^t - \wta^\top \wta(\xi^{t+1}) \prm^t - (\wta^\top {\bf R} \wta \prm^{t+1} - \wta^\top \wta(\xi^{t+1}) \prm^{t+1})) 
    \|^2\big] \notag \\
    & \leq (1- a_\lambda)^2 \expec{\| \wta^\top \momdual^t - \wta^\top {\bf R} \wta \prm^t \|^2} + 2 a_\lambda^2 \expec{ \| \wta^\top \wta(\xi^{t+1}) \prm^{t+1} - \wta^\top {\bf R} \wta \prm^{t+1} \|^2 } \\
    &\quad +2(1 - a_\lambda)^2 \expec{\| \wta^\top {\bf R} \wta(\prm^{t} - \prm^{t+1} ) - \wta^\top \wta(\xi^{t+1}) (\prm^t - \prm^{t+1}) \|^2} \\
    &\leq  (1- a_\lambda)^2 \expec{\| \wta^\top \momdual^t - \wta^\top {\bf R} \wta \prm^t \|^2} + 2 a_\lambda^2 \sigma_A^2 \expec{\| \prm^{t+1} \|_{\wtk}^2} \\
    &\quad + 4 (1-a_\lambda)^2 (\rho_{\max}^2 + \bar{\rho}_{\max}^2) \expec{\| \prm^{t+1} - \prm^t \|^2 }
\end{align}
Now note that by the primal update rule $\prm^{t+1} = (1- \gamma \wta^\top {\bf R} \wta) \prm^t - \eta {\bf v}^t - \alpha(\nabla {\bf f}(\prm^t) - \nabla {\bf f}(\oneotimesT \avgprm^t)) - \alpha (\mom^t - \nabla_\prm {\cal L}(\prm^t, \dprm^t) )$, we have
\begin{align}
    \|\prm^{t+1} \|_{\wtk}^2 \leq 4 (1-\gamma \rho_{\min})^2 \| \prm^t \|^2_{\wtk} + 4 \eta^2 \|{\bf v}^t \|_{\wtk}^2 + 4 \alpha^2 L^2 \|\prm^t \|_{\wtk}^2 + 4 \alpha^2 \| \mom^t - \nabla_\prm {\cal L}(\prm^t, \dprm^t) \|^2
\end{align}
Finally, applying Lemma \ref{aux_lem:x_diff} and the step size condition $a_\lambda \leq \sigma_A^{-1}$ to simplify the coefficients will complete the proof.
\hfill $\square$

\subsection{Proof of Theorem \ref{lem:potential_storm}} \label{app:lemma_potential_storm_proof}
Combining the results of Lemma \ref{lemma:descent_storm}, \ref{lem:consensus_storm}, \ref{lem:v_err_storm}, \ref{lem:xv_inner_storm}, \ref{lem:net_mom_err}, \ref{lem:mom_err}, \ref{lem:mom_dual_err}, when the step sizes satisfy
\begin{align}
&\alpha \leq \min\left\{ \frac{1}{2L}, \sqrt{\frac{a_x}{64 (L^2 + \sampL^2)}}, \frac{\gamma \rho_{\min} \sqrt{1 - \gamma \rho_{\min}}}{\sqrt{32} L }, \frac{1}{\sqrt{12}}, \frac{1}{2L}, \eta \right\}, \label{eq:storm_ss_cond_start} \\
&\gamma \leq \min\left\{ \sqrt{\frac{a_x}{64(\rho_{\max}^2 + \bar{\rho}_{\max}^2)}}, \frac{1}{8\sqrt{a_x} \sigma_A^2} \oscarhide{transient} \right\}, \\
&\eta \beta \leq \frac{1}{12}, \quad \beta \leq 1, \quad 0 < a_x \leq 1, \quad  a_\lambda \leq \frac{1}{\sigma_A},
\end{align}
we obtain
\begin{align}
    F_{t+1} &\leq F_t + \mathbb{C}_{\nabla F}  \expec{\| \nabla F(\avgprm^t) \|^2} + \mathbb{C}_{\sigma} \sigavg + \mathbb{C}_\prm \expec{\| \prm^t \|_{\wtk}^2} + \mathbb{C}_{\bf v} \expec{\| {\bf v}^t \|_{\wtk}^2} \\
    &\quad + \dotp{\prm^t}{{\bf v}^t}_{{\bf C}_{xv}} + \mathbb{C}_{\Delta \avgmom}  \expec{\| \avgmom^t - \frac{1}{n} \oneotimes \nabla {\bf f}(\prm^t) \|^2} \\
    &\quad   + \mathbb{C}_{\Delta\mom} \expec{\| \mom^t - \nabla_\prm {\cal L}(\prm^t, \dprm^t) \|^2} + \mathbb{C}_{\Delta\momdual} \expec{\| \wta^\top \momdual^t - \wta^\top {\bf R}\wta \prm^t \|^2} \\
    &\quad + \mathbb{C}_{\avgmom} \expec{\| \avgmom^t \|^2}
\end{align}
where
\allowdisplaybreaks
\begin{align}
\mathbb{C}_{\nabla F} &= - \frac{\alpha}{2} + \e \cdot 32 (\sampL^2 + L^2) \alpha^2 + \f \cdot 64 \alpha^2 n G + \g \cdot 32 \alpha^2 n (\rho_{\max}^2 + \bar{\rho}_{\max}^2 )\\
    &\leq - \frac{\alpha}{4} \quad \text{when } \begin{cases}
        \alpha \leq \min\left\{\e^{-1} \cdot \frac{1}{512 (\sampL^2 + L^2)}, \quad \g^{-1} \cdot \frac{1}{512 n (\rho_{\max}^2 + \bar{\rho}_{\max}^2)} \right\} \\[0.2cm]
        \f \cdot 64 \alpha^2 n G \leq \frac{\alpha}{24} \Leftarrow \begin{cases} \alpha \leq \f^{-1} \cdot \frac{1}{1024 n (\sampL^2 + L^2)} \\[0.2cm] \gamma^2 \leq \f^{-1} \cdot \frac{\alpha}{1024 n (\rho_{\max}^2 + \bar{\rho}_{\max}^2)} \end{cases}
    \end{cases} \\
\mathbb{C}_{\sigma} &= \e \cdot 2a_x^2 + \f \cdot 4 a_x^2 n \\
\mathbb{C}_{\prm} &= - \a \cdot \frac{\gamma \rho_{\min}}{8} + \frac{\alpha L^2}{n} + 3\beta^2 \rho_{\max}^2(\b \rho_{\min}^{-1} + \c) \\
    &\quad + \d \cdot \left( \frac{3 \alpha^2 L^2}{\eta} + 2 \beta \rho_{\max}^2 - \gamma \beta \rho_{\min}^2 + \frac{5\beta}{2} + \frac{5\alpha^2}{2\eta} \right) \\
    &\quad + \e \cdot 32(\sampL^2 + L^2)(\alpha^2 L^2 + \gamma^2 \rho_{\max}^2) / n \\
    &\quad + \f \cdot\left(\frac{20 a_x^2 \gamma^2 \sigma_A^2}{\alpha^2} + 64 (\alpha^2 L^2 + \gamma^2 \rho_{\max}^2) G \right)\\
    &\quad + \g \cdot \left( 10 a_\lambda^2 \sigma_A^2 + 32 (\rho_{\max}^2 + \bar{\rho}_{\max}^2) (\alpha^2 L^2 + \gamma^2 \rho_{\max}^2) \right) \\
    &\leq -\a \cdot \frac{\gamma\rho_{\min}}{16} \quad \text{when } \begin{cases}
        \alpha \leq  \a \cdot \frac{\gamma \rho_{\min} n }{256 L^2} \\[0.2cm]
        \alpha^2 \leq  \min\left\{ \frac{\a}{\d} \cdot \frac{\eta \gamma \rho_{\min}}{768 L^2}, \quad \frac{\a}{\d} \cdot \frac{\eta \gamma \rho_{\min}}{640} \right\} \\[0.2cm]
        \beta^2 \leq \min\left\{ \frac{\a}{\b} \cdot \frac{\gamma \rho_{\min}^2}{1536 \rho_{\max}^2}, \quad \frac{\a}{\c} \cdot \frac{\gamma \rho_{\min}}{1536 \rho_{\max}^2} \right\} \oscarhide{transient} \\[0.2cm]
        \alpha^3 \leq \frac{\a}{\d} \cdot \frac{\eta \gamma \rho_{\min}}{768 L^4} \oscarhide{transient} \\[0.2cm]
        \beta \leq \min\left\{ \frac{\a}{\d} \cdot \frac{\gamma \rho_{\min}}{512 \rho_{\max}^2},  \quad \frac{\a}{\d} \cdot \frac{\gamma \rho_{\min}}{640} \right\}  \oscarhide{$\beta$ cond.}\\[0.2cm]
        \alpha^2 \leq \min\left\{ \frac{\a}{\e} \cdot \frac{\gamma \rho_{\min} n}{8192 L^2 (\sampL^2 + L^2)}, \quad  \frac{\a}{\g} \cdot \frac{\gamma \rho_{\min}}{8192 L^2 (\bar{\rho}_{\max}^2 + \rho_{\max}^2)} \right\}\\[0.2cm]
        \gamma \leq \min\left\{ \frac{\a}{\e} \cdot \frac{\rho_{\min} n}{8192 \rho_{\max}^2 (\sampL^2 + L^2)} , \quad \frac{\a}{\g} \cdot \frac{\rho_{\min}}{8192 \rho_{\max}^2 (\bar{\rho}_{\max}^2 + \rho_{\max}^2) }  \right\}\\[0.2cm]
        a_x^2 \leq \frac{\a}{\f} \cdot \frac{\alpha^2 \rho_{\min}}{5120 \gamma \sigma_A^2} \oscarhide{transient} \\[0.2cm]
        \f \cdot 64 \alpha^2 L^2 G \leq \a \cdot \frac{\gamma \rho_{\min}}{256} \Leftarrow \begin{cases} \alpha^2 \leq \frac{\a}{\f} \cdot \frac{\gamma \rho_{\min}}{32768 L^2(\sampL^2 + L^2)} \\[0.2cm] 
                        \gamma \leq \frac{\a}{\f} \cdot \frac{\rho_{\min}}{32768 L^2 (\bar{\rho}_{\max}^2 + \rho_{\max}^2) }\end{cases} \\[0.2cm]
        \f \cdot 64 \gamma^2 \rho_{\max}^2 G \leq \a \cdot \frac{\gamma \rho_{\min}}{256} \Leftarrow \begin{cases} \gamma \leq \frac{\a}{\f} \cdot \frac{\rho_{\min}}{32768 \rho_{\max}^2 (\sampL^2 + L^2)} \\[0.2cm] \gamma^3 \leq \frac{\a}{\f} \cdot \frac{\alpha^2 \rho_{\min}}{32768 \rho_{\max}^2 (\bar{\rho}_{\max}^2 + \rho_{\max}^2) } \end{cases} \\[0.2cm]
        a_\lambda^2 \leq \frac{\a}{\g} \cdot \frac{\gamma \rho_{\min}}{2560 \sigma_A^2} \oscarhide{$a_\lambda$ cond.} \\[0.2cm]
    \end{cases} \\
\mathbb{C}_{\bf v}&= -\d \cdot \frac{\eta}{2} + \a \cdot \frac{\eta^2}{(1-\gamma \rho_{\min})(1- \gamma \rho_{\min}/2)} + (\beta + \frac{\alpha}{\eta}) (\b \rho_{\min}^{-1} + \c) \\
    &\quad + \e \cdot \frac{32(\sampL^2 + L^2) \eta^2}{n} + \f \cdot \left(\frac{16 a_x^2 \eta^2 \gamma^2 \sigma_A^2}{\alpha^2} + 64 \eta^2 G \right) + \g \cdot 16\eta^2  \\
    &\leq -\d \cdot \frac{\eta}{4} \quad \text{when }
    \begin{cases}
        \eta \leq \frac{\d}{\a} \cdot (28 (1-\gamma\rho_{\min})(1-\gamma\rho_{\min} / 2))^{-1} \oscarhide{transient} \\[0.2cm]
        \beta + \frac{\alpha}{\eta} \leq \min \left\{ \frac{\d}{\b} \cdot \frac{\eta \rho_{\min}}{28}, \quad \frac{\d}{\c} \cdot \frac{\eta}{28} \right\} \\[0.2cm]
        \eta \leq \min \left\{ \frac{\d}{\e}\cdot \frac{n}{896 (\sampL^2 + L^2)}, \quad \frac{\d}{\g} \cdot (1/448) \right\} \\[0.2cm]
        a_x^2 \eta \gamma^2 \leq \frac{\d}{\f} \cdot \frac{\alpha^2}{448 \sigma_A^2} \oscarhide{transient} \\[0.2cm]
        \f \cdot 64 \eta^2 G \leq \d \cdot \frac{\eta}{28} \Leftarrow \begin{cases} \eta \leq \frac{\d}{\f} \cdot (3584 (\sampL^2 + L^2))^{-1} \\[0.2cm] \eta \gamma^2 \leq \frac{\d}{\f} \cdot \frac{\alpha^2}{3584 (\bar{\rho}_{\max}^2 + \rho_{\max}^2) } \oscarhide{$\gamma$'s $n$ dependence}\end{cases}
    \end{cases} \\
\mathbf{C}_{xv}&= \left(-\a \cdot \frac{2\eta (1 - \gamma \rho_{\min}/4)}{1 - \gamma \rho_{\min}} + \b \cdot 2\beta \right) \wtk \\
    &\quad + \left( \a \cdot \frac{2\eta \gamma (1- \gamma \rho_{\min}/4)}{1- \gamma \rho_{\min}} + \c \cdot 2\beta - \d \cdot (\gamma + \eta \beta) \right) \wta^\top {\bf R} \wta \\
    &= {\bf 0} \quad \text{when } \begin{cases}
         \b = \frac{\eta ( 1- \gamma \rho_{\min} / 4)}{\beta(1-\gamma \rho_{\min})} \a, \\[0.2cm]
         \c = \frac{\d(\gamma + \eta\beta) - \a \cdot 2 \eta \gamma (1- \gamma \rho_{\min} / 4) (1-\gamma \rho_{\min})^{-1}}{2 \beta}
    \end{cases} \\
\mathbb{C}_{\Delta\avgmom}&= -\e \cdot a_x + \alpha  \\
    &\leq - \e \cdot \frac{a_x}{2} \quad \text{when } \alpha \leq \e \cdot \frac{a_x}{2} \oscarhide{$\e$ cond.} \label{eq:e_necess_cond} \\
\mathbb{C}_{\Delta\mom}&= -\f \cdot \frac{a_x}{4} + \a \cdot \frac{2 \alpha^2}{\gamma \rho_{\min}} + \d \cdot ((\frac{3}{\eta} + \frac{1}{2})\alpha^2 + 2\alpha^2\beta + \frac{2 \alpha^4}{ \eta}) \\
    &\quad + \e \cdot \frac{8 (\sampL^2 + L^2) \alpha^2}{n} + \g \cdot 8\alpha^2 ( 1 + \rho_{\max}^2 + \bar{\rho}_{\max}^2) \\
    &\leq -\f \cdot \frac{a_x}{8} \quad \text{when } \begin{cases}
        \alpha^2 \leq \min\left\{\frac{\f}{\a} \cdot \frac{a_x \gamma \rho_{\min}}{64}, \quad \frac{\f}{\d} \cdot \frac{a_x \eta}{384}, \quad \frac{\f}{\d} \cdot \frac{a_x}{64}\right\} \oscarhide{$\d, \f$ cond.} \\[0.2cm] 
        \alpha^2 \beta \leq \frac{\f}{\d} \cdot \frac{a_x}{256} \oscarhide{transient} \\[0.2cm]
        \alpha^4 \leq \frac{\f}{\d} \cdot \frac{a_x \eta}{256} \oscarhide{transient} \\[0.2cm]
        \alpha^2 \leq \min \left\{ \frac{\f}{\e} \cdot \frac{a_x n}{256 (\sampL^2 + L^2)}, \quad \frac{\f}{\g} \cdot \frac{a_x}{256(1+ \rho_{\max}^2 + \bar{\rho}_{\max}^2 )} \right\}
    \end{cases}\\
\mathbb{C}_{\Delta\momdual}&= -\g \cdot a_\lambda + 4\beta(\b \rho_{\min}^{-1} + \c) + \d \cdot \frac{\beta}{2} \\
    &\leq -\g \cdot \frac{a_\lambda}{2} \quad \text{when }
        \beta \leq \min\left\{  \frac{\g}{\b} \cdot \frac{a_\lambda \rho_{\min}}{24} , \quad \frac{\g}{\c}\cdot \frac{a_\lambda}{24}, \quad \frac{\g}{\d} \cdot \frac{a_\lambda}{3} \oscarhide{transient }\right\} \oscarhide{$\g$ cond.}\\
\mathbb{C}_{\avgmom} &= -\frac{\alpha}{4} + \frac{4 n \alpha^3 L^2}{\eta} (\b \rho_{\min}^{-1} + \c) + \d \cdot \frac{\alpha^2 L^2 n}{2\eta} \\
&\leq -\frac{\alpha}{8} \quad\text{when } \begin{cases}
    \alpha^2 \leq \min \left\{  \b^{-1} \cdot \frac{\eta \rho_{\min}}{96 nL^2}, \quad \c^{-1} \cdot \frac{\eta}{96 nL^2} \right\}\\[0.2cm]
    \alpha \leq \d^{-1} \frac{\eta}{12 nL^2} \oscarhide{$\eta$ kills $n$ for $\alpha$.}
\end{cases} \label{eq:storm_ss_cond_end}
\end{align}
Notice that the hyperparameter choices in \eqref{eq:ss_choice_storm_start} - \eqref{eq:ss_choice_storm_end} will satisfy all of the above conditions, therefore completes the proof to \eqref{eq:potential_conv}.
\hfill $\square$

\subsection{Auxiliary Lemma}
\begin{lemma} \label{aux_lem:x_diff} Under Assumption \ref{assm:lip}, \ref{assm:rand-graph},
    \begin{align}
        \| \prm^{t+1} - \prm^t \|^2 &\leq 2\alpha^2 \| \mom^t - \nabla_{\prm} {\cal L}(\prm^t, \dprm^t) \|^2 + 8 \alpha^2 n \| \nabla F(\avgprm^t) \|^2 \\
        &\quad + 8 \eta^2 \| {\bf v}^t \|_{\wtk}^2  + (8 \alpha^2 L^2 + 8 \gamma^2 \rho_{\max}^2)  \| \prm^t \|_{\wtk}^2
    \end{align}
\end{lemma}

{\it Proof of Lemma \ref{aux_lem:x_diff}.}
\begin{align}
    &\|\prm^{t+1} - \prm^t \|^2 = \alpha^2 \| \mom^t \|^2 \leq 2\alpha^2 \| \mom^t - \nabla_\prm {\cal L}(\prm^t, \dprm^t) \|^2 + 2 \alpha^2 \| \nabla_\prm {\cal L}(\prm^t, \dprm^t) \|^2
\end{align}
Upon noting that $\| {\bf v}^t \|^2_{\wtk} = \left\| \wta^\top \dprm^t - (\alpha/\eta)(n^{-1} {\bf 1}{\bf 1}^{\top} - {\bf I}_n) \otimes {\bf I}_d \nabla {\bf f}(\oneotimesT \avgprm^t) \right\|^2$, we expand
\begin{align}
     &\| \nabla_\prm {\cal L}(\prm^t, \dprm^t) \|^2 = \left\| \nabla {\bf f}(\prm^t) + \frac{\eta}{\alpha} \wta^\top \dprm^t + \frac{\gamma}{\alpha} \wta^\top {\bf R}\wta \prm^t \right\|^2 \\
     & = \left\| \oneotimesT \nabla F(\avgprm^t) + \frac{\eta}{\alpha} \wta^\top \dprm^t - ( \oneotimesT \nabla F(\avgprm^t) - \nabla {\bf f}(\oneotimesT \avgprm^t)) - (\nabla {\bf f}(\oneotimesT \avgprm^t) - \nabla {\bf f}(\prm^t) ) + \frac{\gamma}{\alpha} \wta^\top {\bf R}\wta \prm^t \right\|^2 \notag \\
     &\leq 4 \| \oneotimesT \nabla F(\avgprm^t) \|^2 + \frac{4 \eta^2}{\alpha^2} \left\| \wta^\top \dprm^t - \frac{\alpha}{\eta} (\oneotimesT \nabla F(\avgprm^t) - \nabla {\bf f}(\oneotimesT \avgprm^t)) \right\|^2 + 4 \| \nabla {\bf f}(\oneotimesT \avgprm^t) - \nabla {\bf f}(\prm^t)  \|^2 \notag \\
     &\quad + \frac{4 \gamma^2 \rho_{\max}^2}{\alpha^2} \| \prm^t \|_{\wtk}^2 \\
     &\leq 4n \| \nabla F(\avgprm^t) \|^2 + \frac{4 \eta^2}{\alpha^2} \| {\bf v}^t \|_{\wtk}^2 + (4 L^2 +  \frac{4 \gamma^2 \rho_{\max}^2}{\alpha^2} ) \| \prm^t \|_{\wtk}^2
\end{align}
Combining the above inequalities completes the proof.
\hfill $\square$


\section{Ablation Study} \label{app:ablation}
In this section, we investigate how {\algname} behaves under different problem configurations, for instance, the different levels of data heterogeneity, random graph sparsity, graph topology, gradient noise and dual momentum. Unless specified explicitly, we assume the experiment adopts ${\cal G}$ as the fully connected (complete) graph topology.

\subsection{Data Heterogeneity}
\begin{figure}[h]
    \centering
    \includegraphics[width=0.85\textwidth]{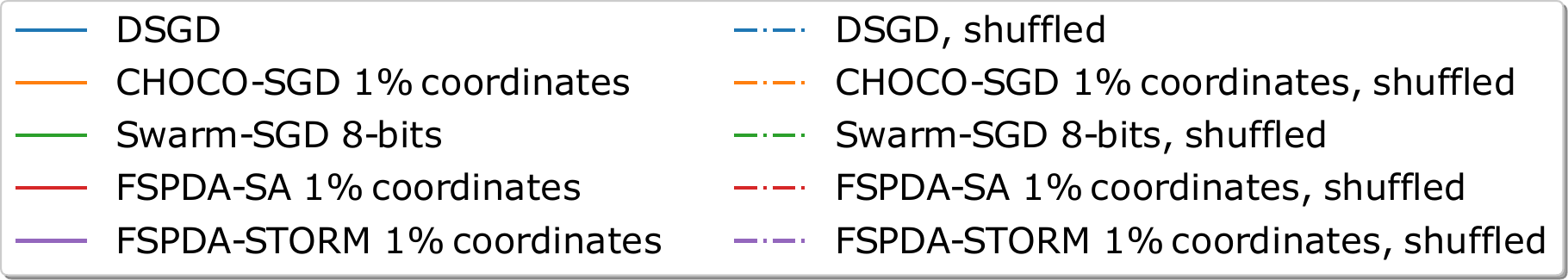} \\
    \vspace{0.25cm}
    \includegraphics[width=0.32\textwidth]{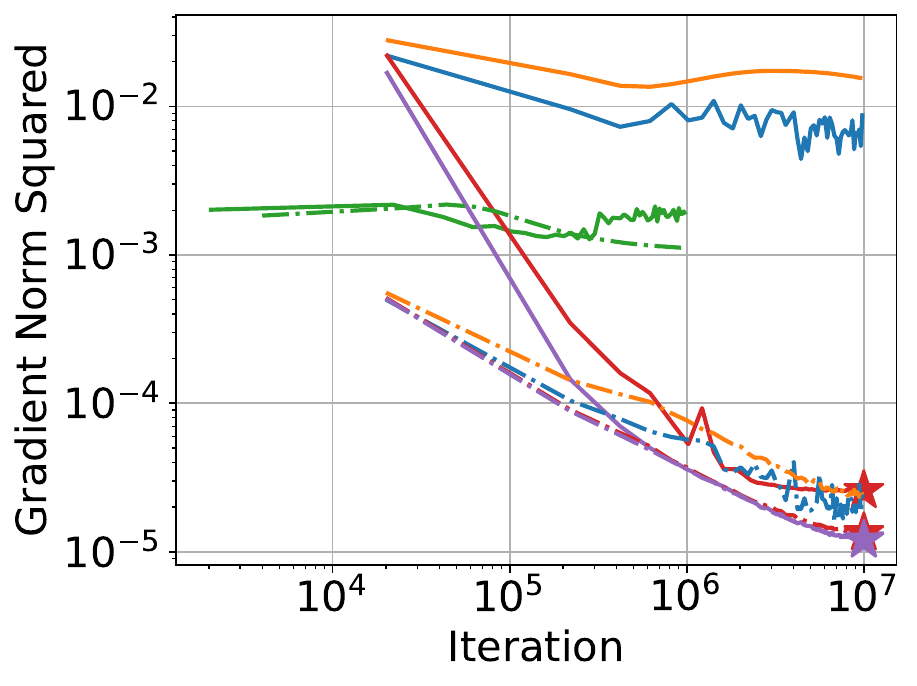} 
    \includegraphics[width=0.32\textwidth]{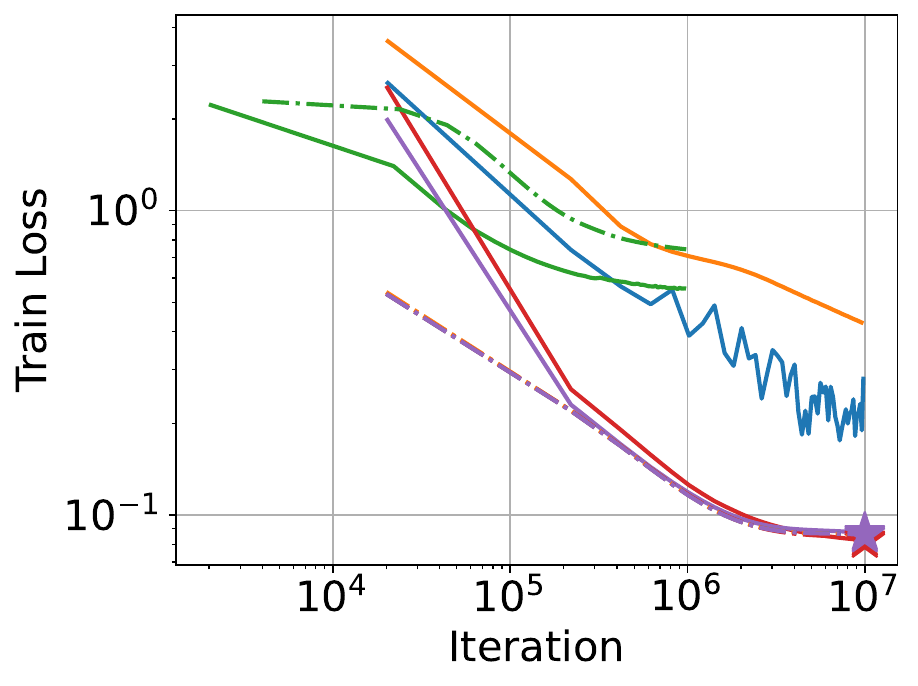}
    \includegraphics[width=0.32\textwidth]{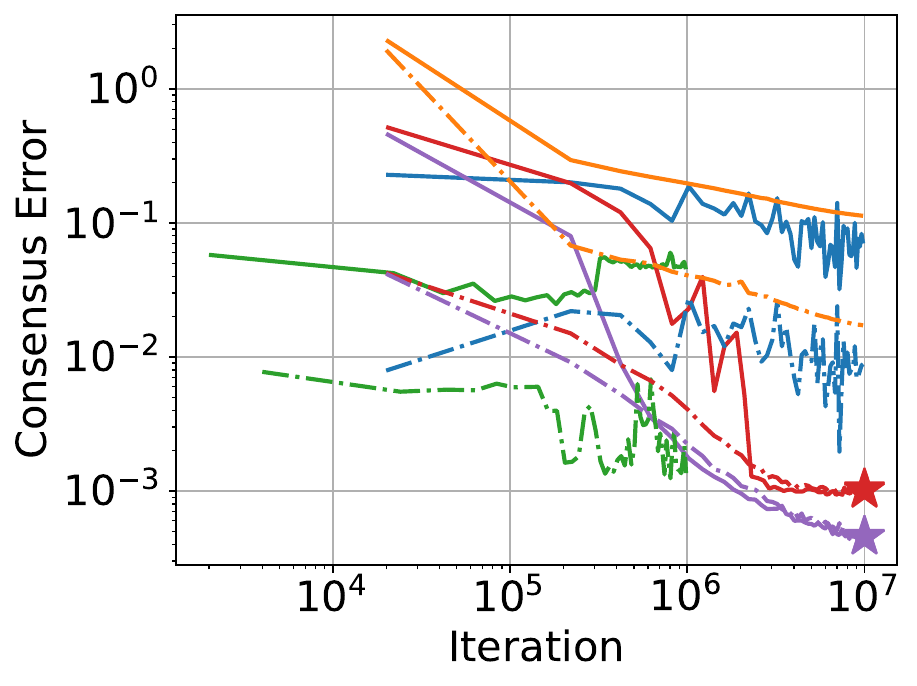}
    \\
    \includegraphics[width=0.32\textwidth]{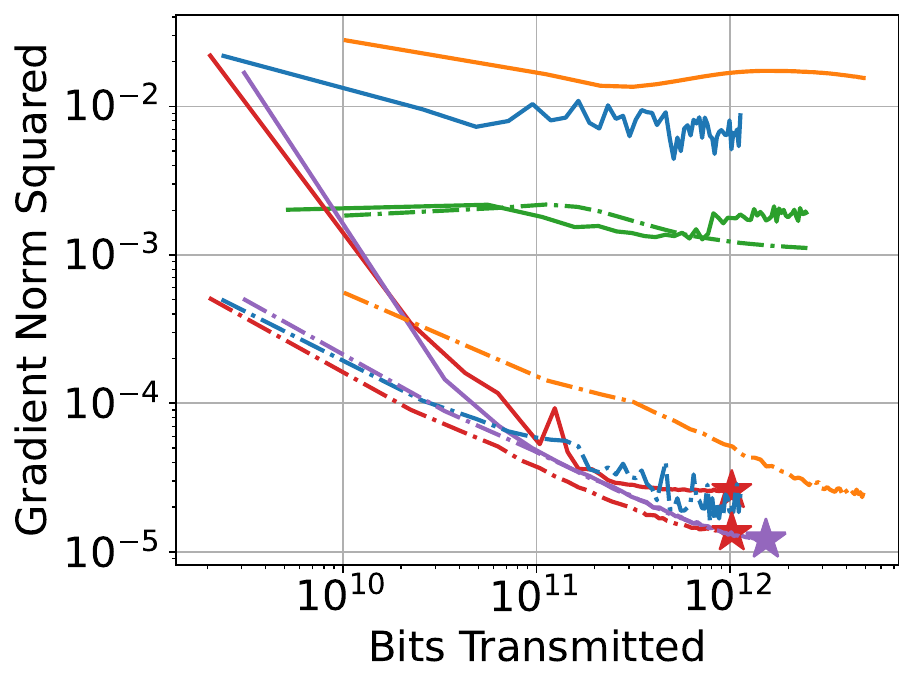}
    \includegraphics[width=0.32\textwidth]{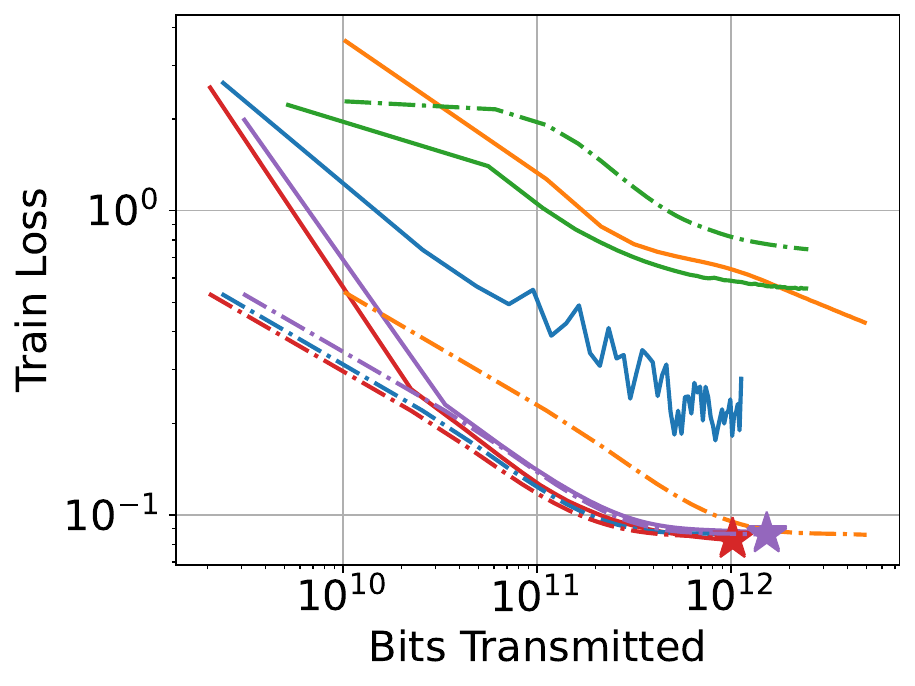}
    \includegraphics[width=0.32\textwidth]{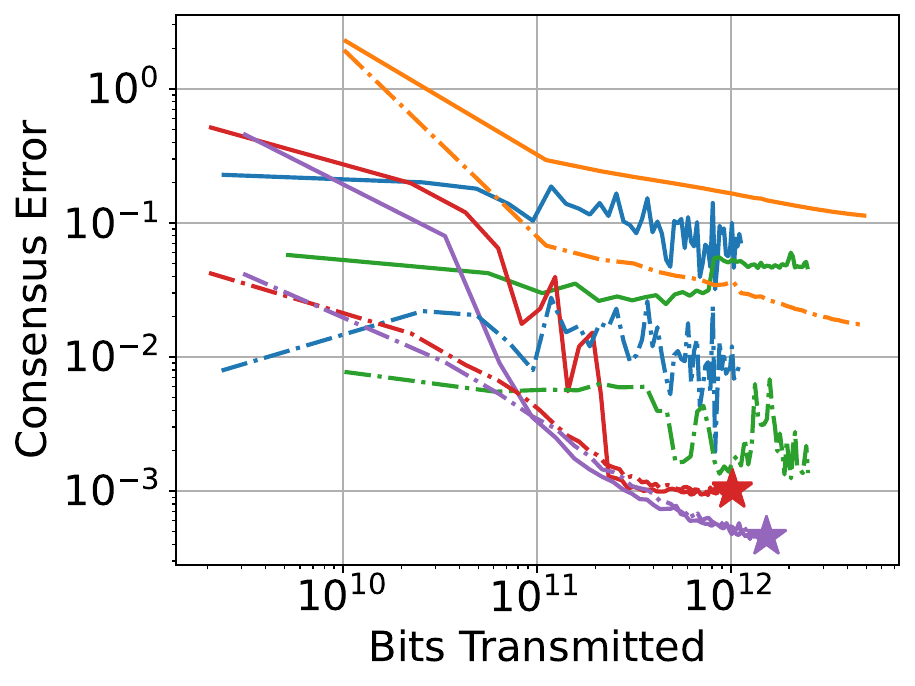}
    \\
    \includegraphics[width=0.32\textwidth]{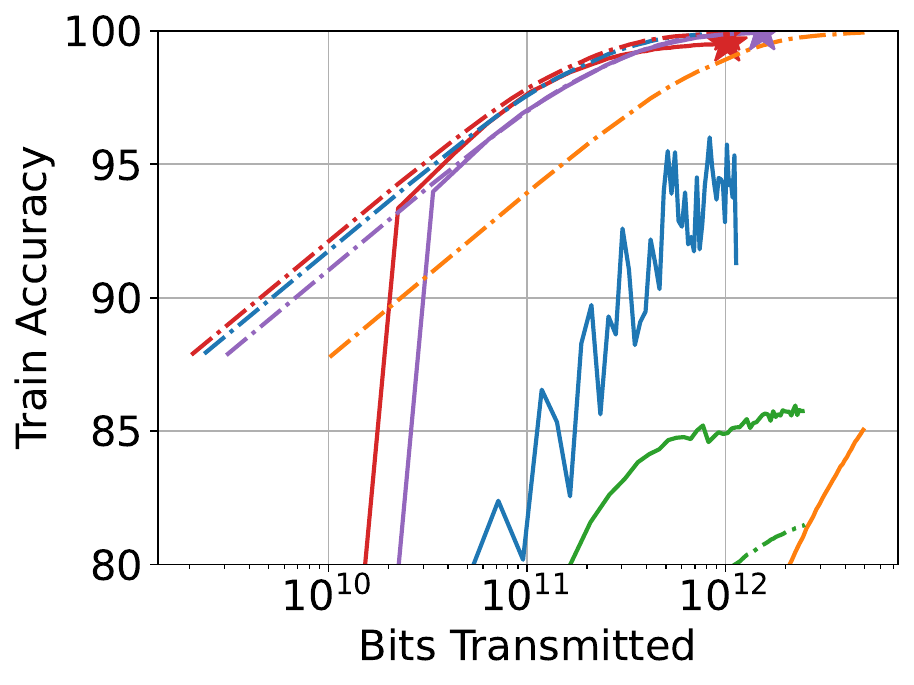}
    \includegraphics[width=0.32\textwidth]{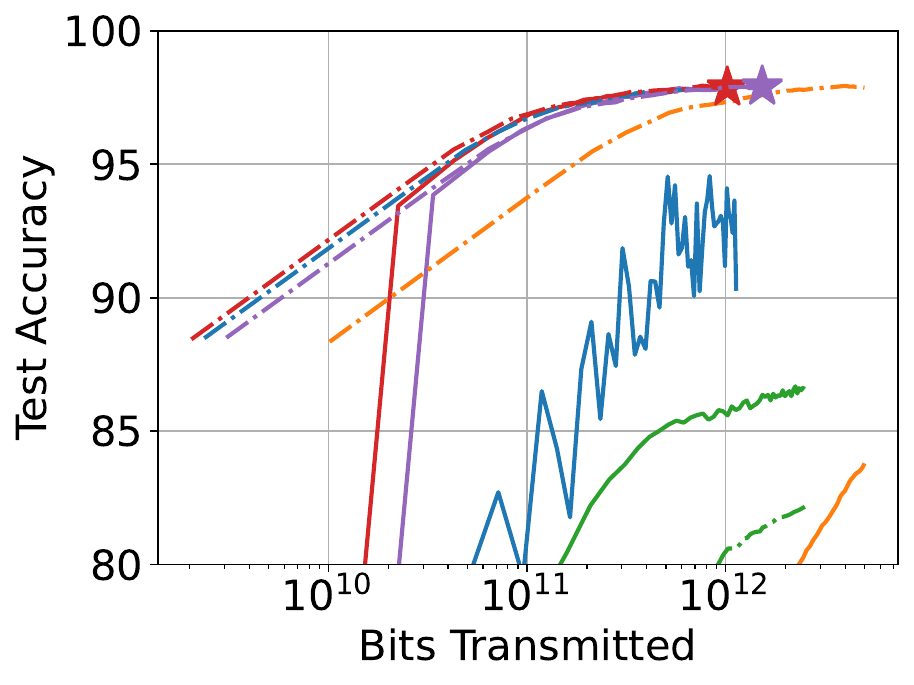}
    \caption{Feed-forward neural network classification training on MNIST with two levels of data heterogeneity.
    } 
    \label{fig:mnist}
\end{figure}

To study the effect of heterogeneity of data distribution across agents, we experiment with two types of data splitting for MNIST: ($i$) the dataset is split into $n=10$ disjoint sets according to the class labels, or ($ii$) the dataset is split into $n$ evenly distributed disjoint sets by shuffling. Setup ($i$) creates a large discrepancy across local objective functions, i.e., a larger data heterogeneity. Figure \ref{fig:mnist} compares the performance of {\algname} and the benchmark algorithms under the above setup which demonstrates the robustness of {\algname} under heterogeneous data distribution. For instance, the performance of {\tt CHOCO-SGD} and {\tt DSGD} hugely degrades in the heterogeneous setup ($i$), while that of {\algnamesa} is only affected by a small margin and {\algnamevr} is able to converge to the same low error in both setup ($i$) and ($ii$). We list the hyperparameters used in Figure \ref{fig:mnist} by Table \ref{tab:hypprm_mnist} in Appendix \ref{app:exp}.

\subsection{Random Graph Sparsity}
\begin{figure}[h]
    \centering
    \includegraphics[width=0.85\textwidth]{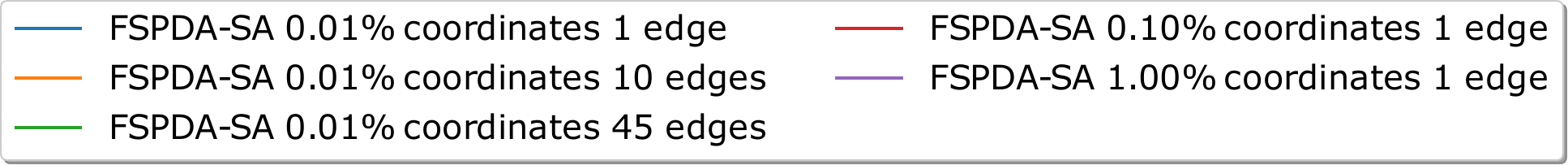} \\
    \vspace{0.25cm}
    \includegraphics[width=0.35\textwidth]{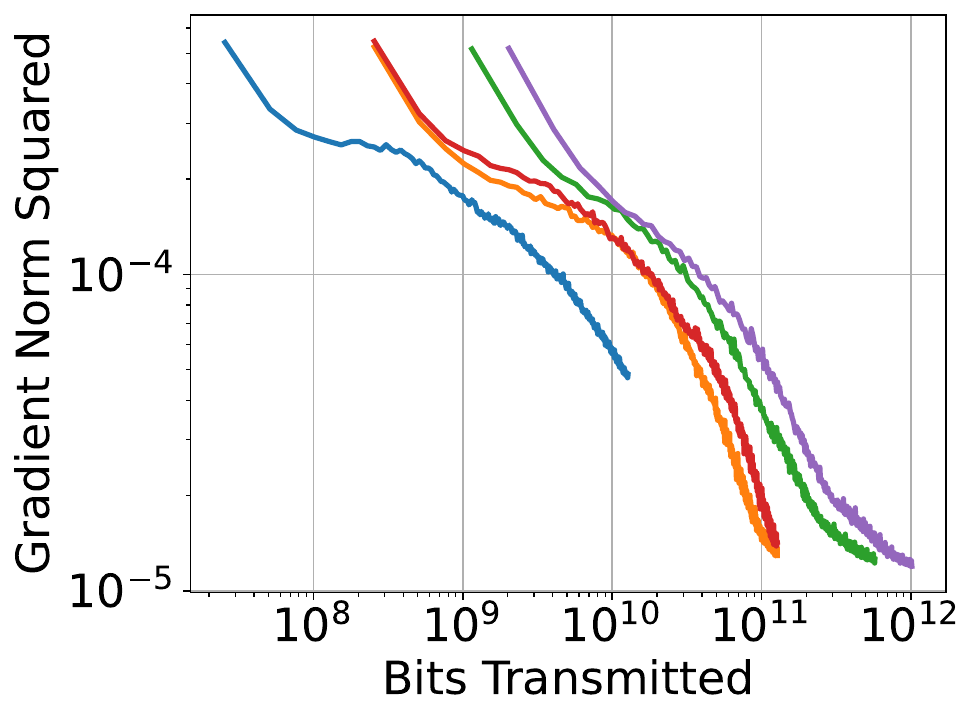}
    \includegraphics[width=0.35\textwidth]{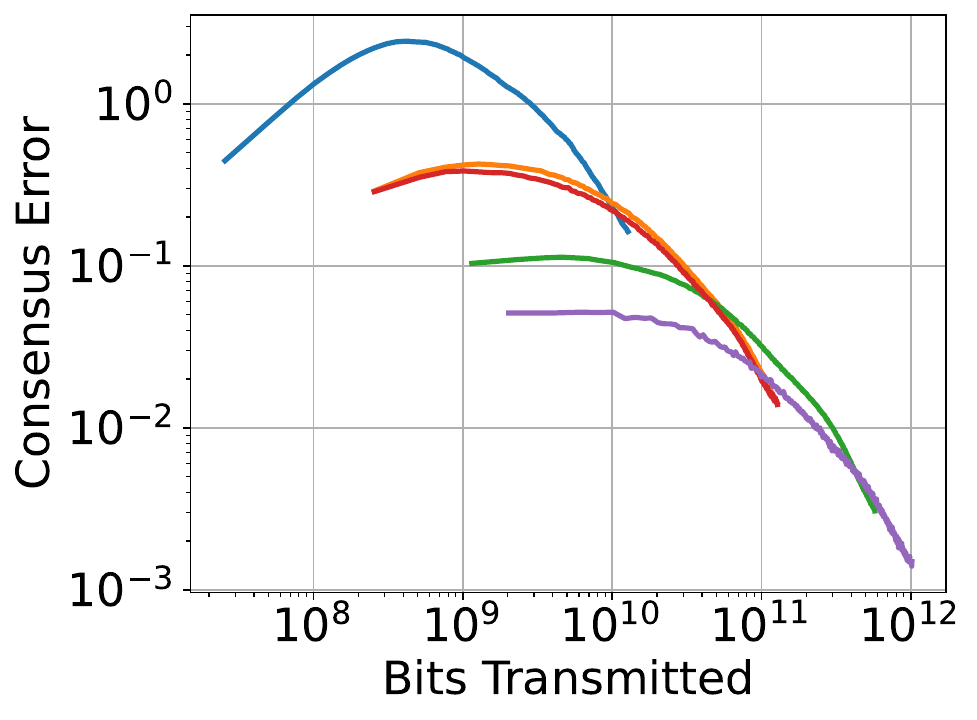} \\
    \includegraphics[width=0.35\textwidth]{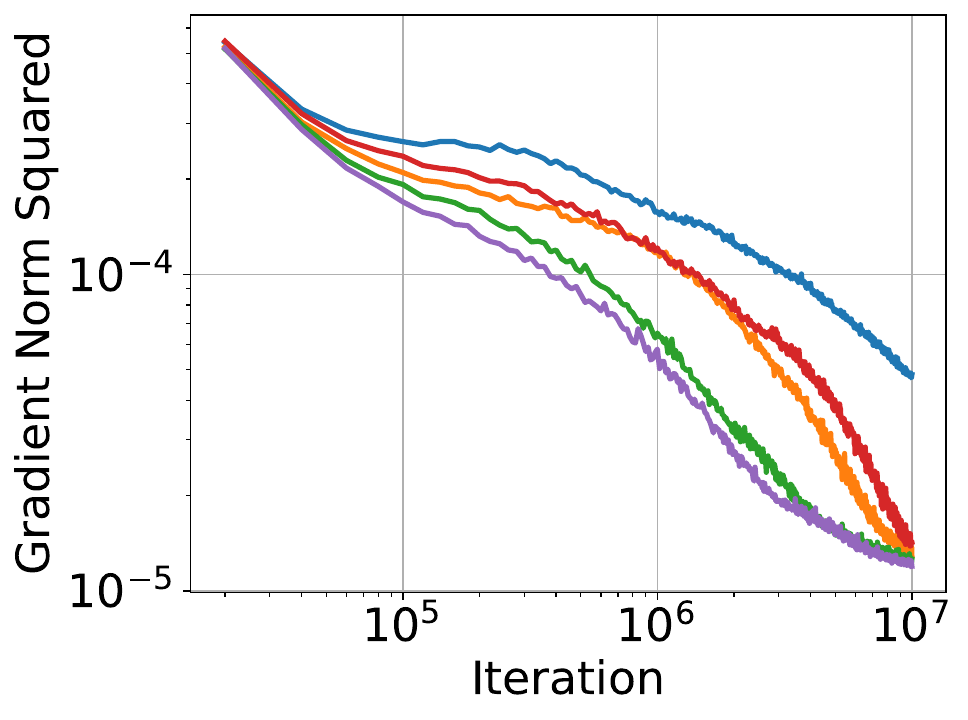}
    \includegraphics[width=0.35\textwidth]{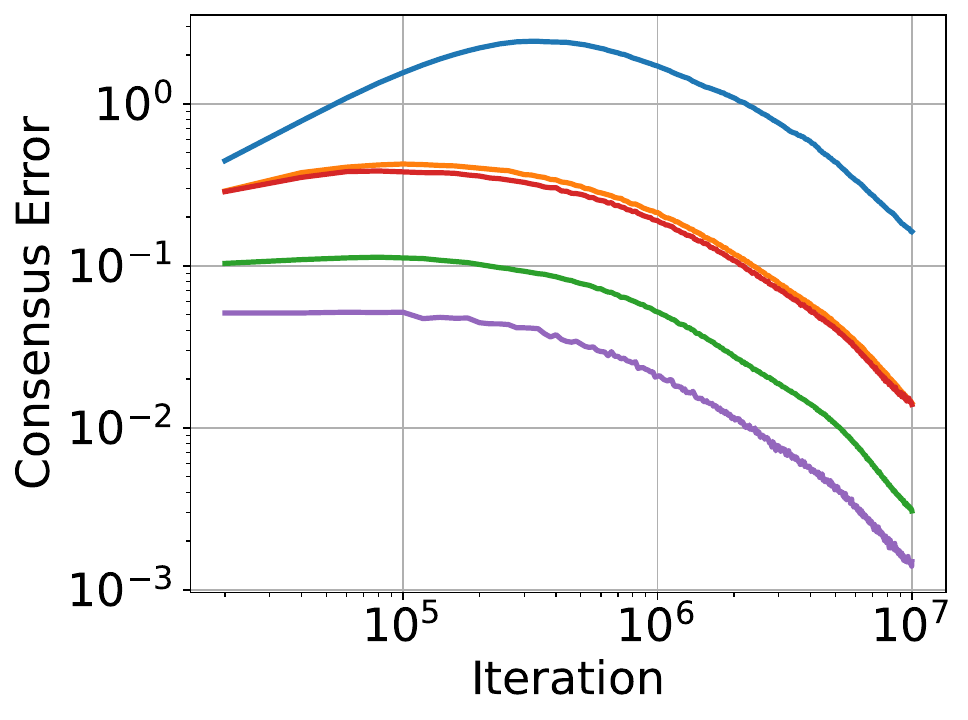}
    \caption{Feed-forward neural network classification training on shuffled MNIST. Random graph of $k$ edges in expectation ($k \in \{1, 10, 45\}$) is drawn from a complete topology per iteration. 
    }
    \label{fig:mnist_sparsity}
\end{figure}

Next, we study the effects of random graph sparsity by the experiments shown in Figure \ref{fig:mnist_sparsity}. As the random graph sparsity decreases, the random graph variance $\sigma_A^2$ in Assumption~\ref{assm:graph_var} decreases which will significantly improve the consensus error. We observe from the figure that despite the different levels of consensus error, 4 out of 5 configurations eventually converge to the same stationarity. In the extreme case of 0.01\% sparsity with one-edge random graph, we see the dominance of sparsity error which leads to a longer transient time as expected from Theorem \ref{thm:main}. We conclude that for a fixed number of iteration, a certain amount of communication sparsity can be tolerated in {\algname} without degradation in optimization error. For Figure \ref{fig:mnist_sparsity}, we tuned {\algnamesa} to use the step sizes $\alpha = 10^{-4}, \eta = 10^{-6}, \gamma = 0.5, \beta = 1$ for all configurations.


\subsection{Graph Topology}
\begin{figure}[h]
    \centering
    \includegraphics[width=0.75\textwidth]{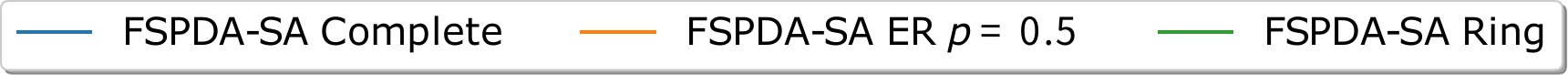} \\
    \vspace{0.25cm}
    \includegraphics[width=0.35\textwidth]{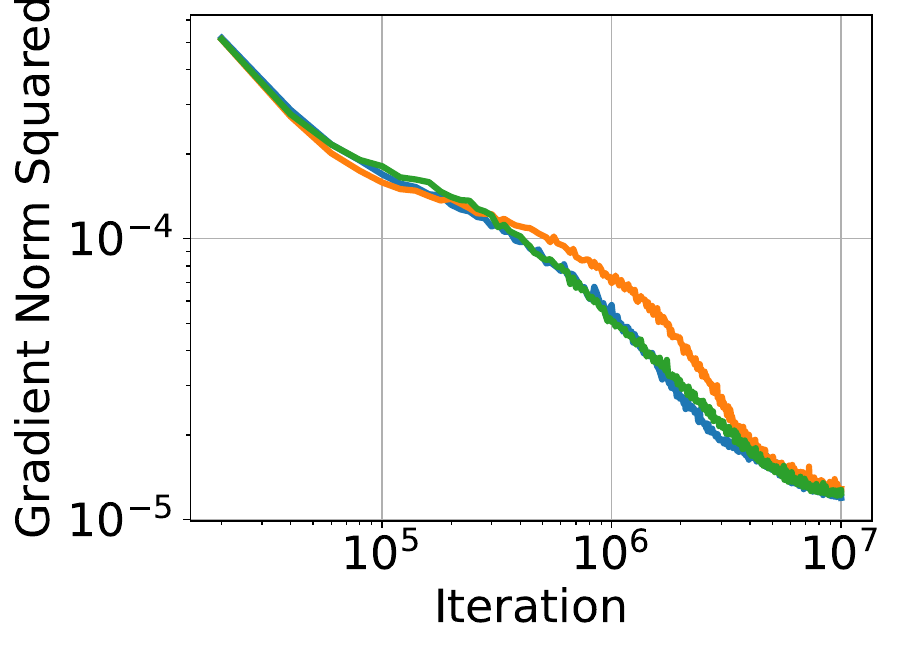}
    \includegraphics[width=0.35\textwidth]{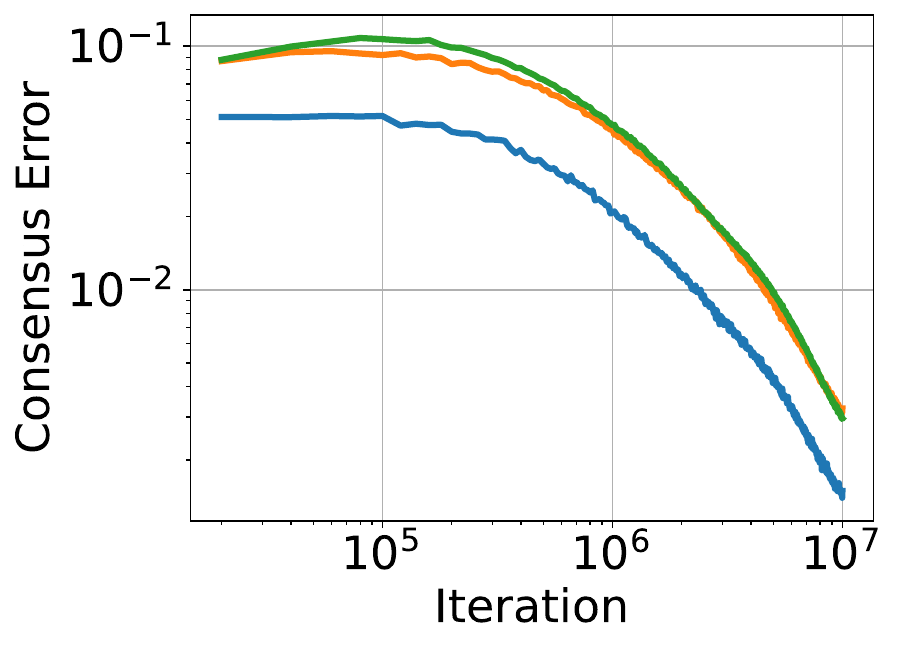}
    \caption{Feed-forward neural network classification training on shuffled MNIST with different graph topology ${\cal G}$. Only one edge is activated per iteration, exchanging 1\% coordinates of model parameters. 
    } 
    \label{fig:mnist_topology}
\end{figure}

We then study the effects of network topology $\mathcal{G}$ in {\algnamesa} by drawing one-edge random graphs from a complete graph, an ER graph with probability $p = 0.5$ and a ring graph in Figure \ref{fig:mnist_topology}. Note that the communication cost per iteration is the same across different topologies due to the use of one-edge random graph. The result indicates the transient effect of topology where it only slow down the convergence of consensus error while converging to the same level of stationarity. For Figure \ref{fig:mnist_topology}, we tuned {\algnamesa} to use the step sizes $\alpha = 10^{-4}, \eta = 10^{-6}, \gamma = 0.5, \beta = 1$ for all configurations.

\subsection{Deterministic Gradient}
\begin{figure}[h]
    \centering
    \includegraphics[width=0.7\textwidth]{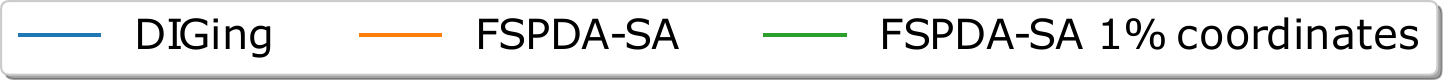} \\
    \vspace{0.25cm}
    \includegraphics[width=0.35\textwidth]{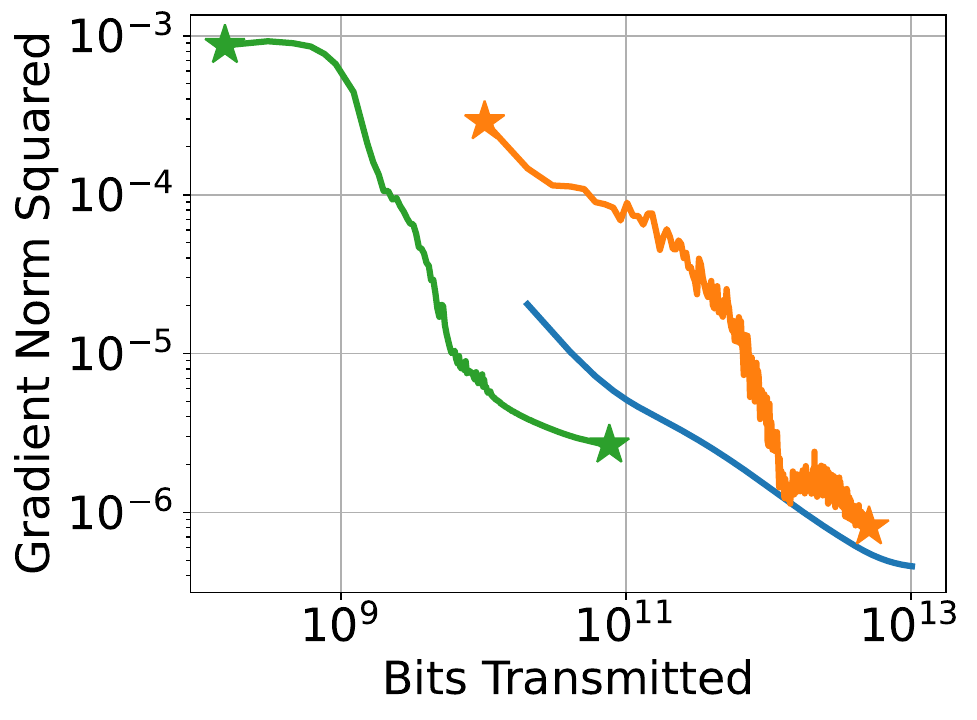}
    \includegraphics[width=0.35\textwidth]{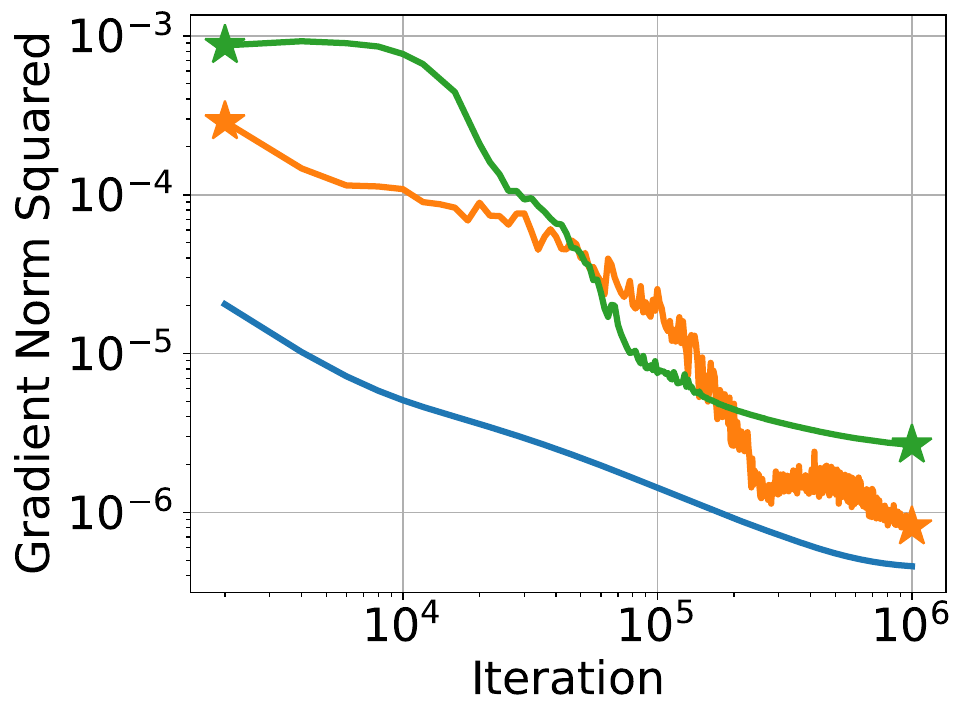} \\
    \includegraphics[width=0.35\textwidth]{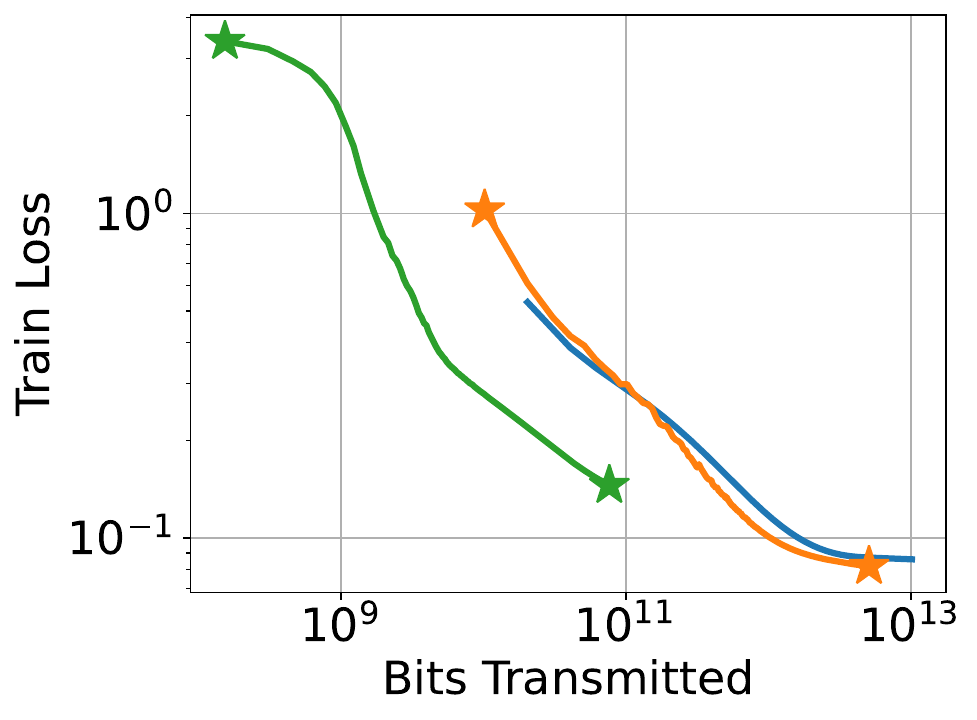}
    \includegraphics[width=0.35\textwidth]{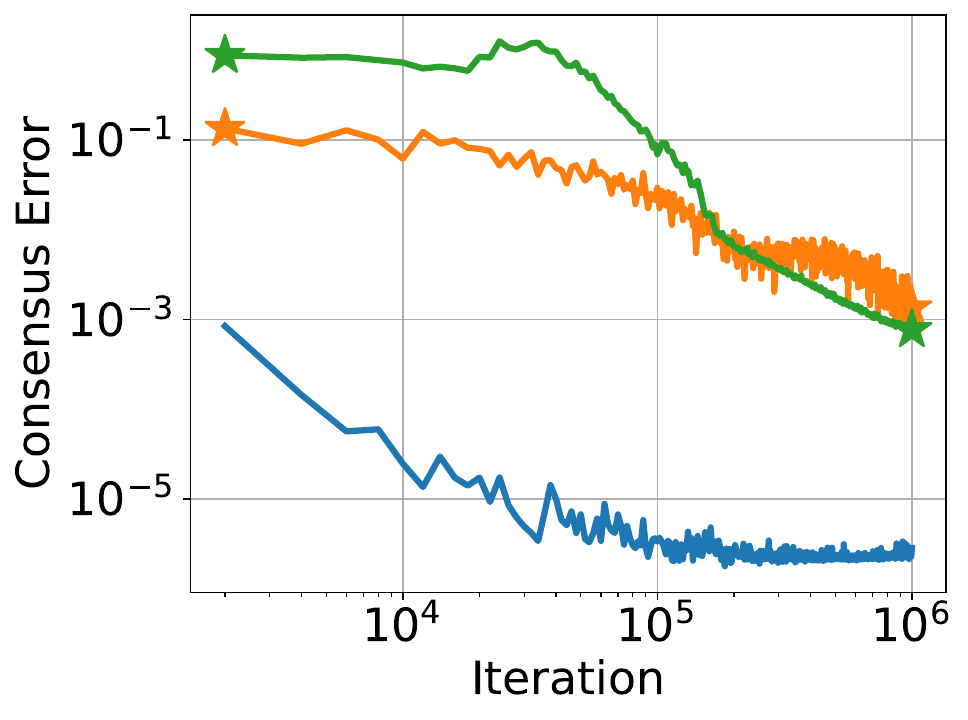}
    \caption{Feed-forward neural network classification training on class separated MNIST using exact local gradient, i.e., $\sigma_i = 0~\forall i$. One-edge random graph is drawn from a complete topology per iteration.
     } 
    \label{fig:mnist_exact_grad}
\end{figure}

For the case when the gradient estimate is exact, i.e., $\sigavg = 0$, we compare the performance of {\algnamesa} against deterministic gradient algorithm {\tt DIGing} \citep{nedic2017achieving} in Figure \ref{fig:mnist_exact_grad}. Despite {\algnamesa} only performs model parameter gossip while {\tt DIGing} performs an extra step of gradient tracker gossip, {\algnamesa} shows comparative performance and reduced the communication cost by half.

Also, in the absence of stochastic gradient noise, notice that the effect of parameter sparsification immediately transfer to a slower optimization convergence. This is in line with our theorem by observing the convergence rate in 
Theorem \ref{thm:main}, where $\sigma_A^4$ remains dominant in the convergence bound. We list the hyperparameters used in Figure \ref{fig:mnist_exact_grad} by Table \ref{tab:hypprm_mnist_exact_grad} in Appendix \ref{app:exp}.

\subsection{Dual Momentum}
\begin{figure}[h]
    \centering
    \includegraphics[width=0.97\textwidth]{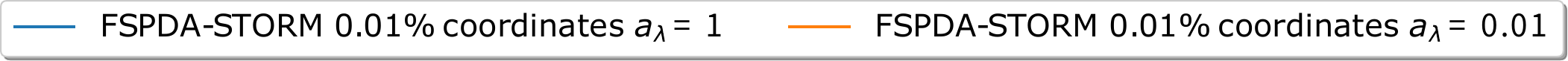} \\
    \vspace{0.25cm}
    \includegraphics[width=0.32\textwidth]{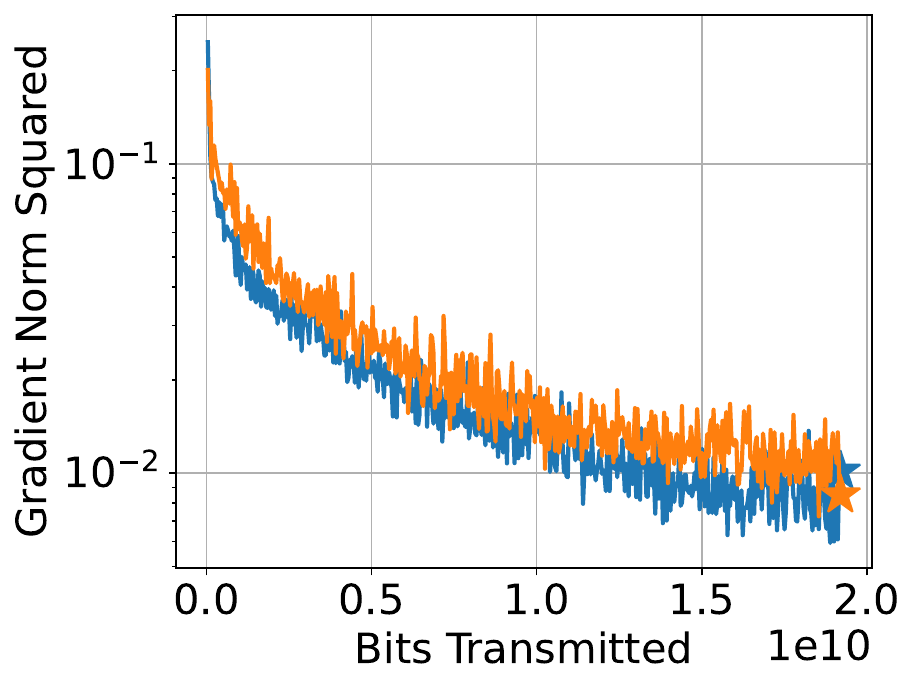}
    \includegraphics[width=0.32\textwidth]{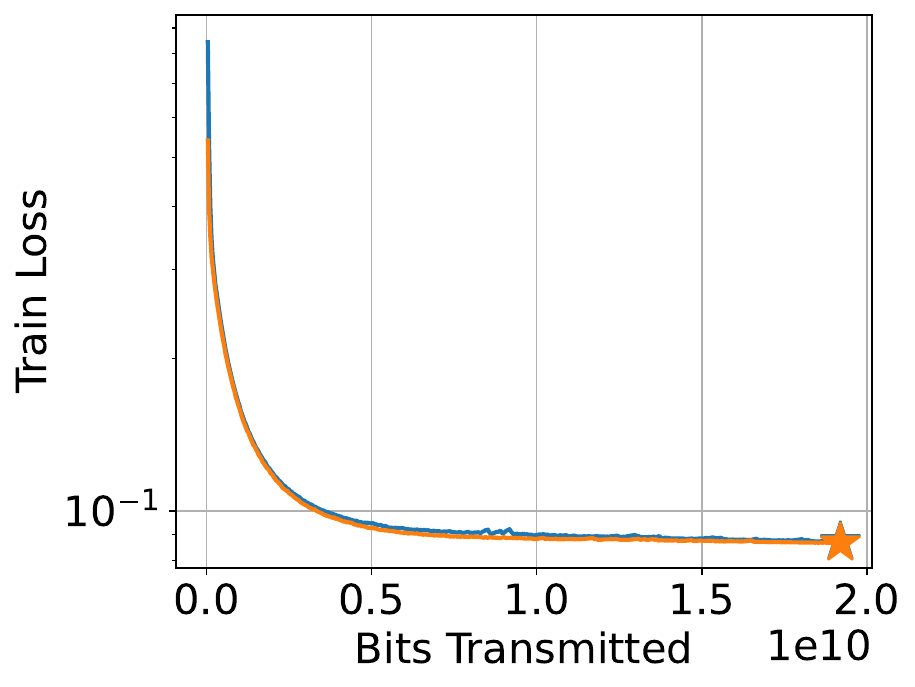}
    \includegraphics[width=0.32\textwidth]{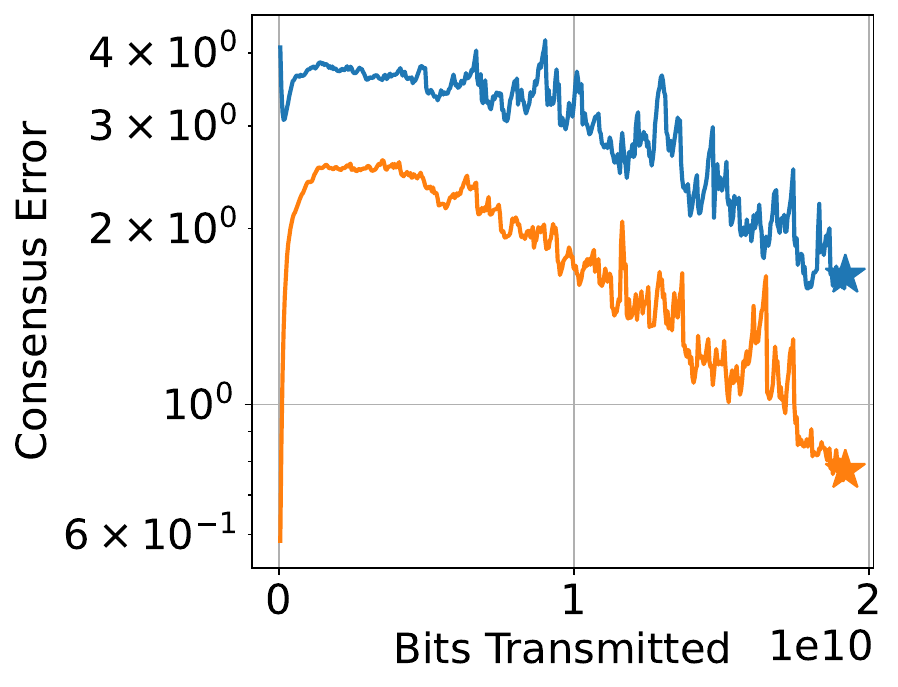}
    \caption{Feed-forward neural network classification training on shuffled MNIST using $10^7$ iterations. One-edge random graph is drawn from a complete topology per iteration.}
    \label{fig:mnist_dual_mom}
\end{figure}

To investigate the benefits of dual momentum in {\algnamevr}, we construct a case where both the primal and dual stochastic gradients carry large variance error. In Figure \ref{fig:mnist_dual_mom}, the local objective function gradient is estimated by batch size 1, and one-edge random graphs with 0.01\% coordinate sparsification is adopted for communication. We can observe that applying dual momentum ($a_\lambda = 0.01$) outperforms not applying dual momentum ($a_\lambda = 1$) in terms of consensus error convergence. We list the hyperparameters used in Figure \ref{fig:mnist_dual_mom} by Table \ref{tab:hypprm_mnist_dual_mom} in Appendix \ref{app:exp}.

\section{Experiment Hyperparameters} \label{app:exp}
\begin{table}[h]
    \centering
    \begin{tabular}{l|c c c c c c}
    \toprule 
     {\tt FSPDA} & $\alpha$  & $\eta$ & $\gamma$ & $\beta$ & $a_x$ & $a_\lambda$ \\[0.05cm]
    {\tt -SA} (10\% sparse coor.) & $10^{-4}$ & $10^{-5}$ & $0.5$ & $1$ & - & - \\
    {\tt -STORM} (6.7\% sparse coor.) & $10^{-3}$ & $10^{-2}$ & $0.5$ & $0.1$ & $10^{-2}$ & $10^{-2}$ \\
    \midrule
     & {\bfseries Opt. S.S.} & {\bfseries Local Steps}  & $\eta_s$ & - & - & - \\[0.05cm]
    {\tt K-GT} & $10^{-4}$ & $150$ & $1$ & - & - & -\\
    \midrule
     & {\bfseries Opt. S.S.} & {\bfseries Local Steps} & - & - & - & - \\[0.05cm]
    {\tt LED} & $10^{-4}$ & $75$ & - & - & - & -\\
    \midrule
     & {\bfseries Opt. S.S.} & $\tau$ & $p$ & - & - & - \\[0.05cm]
    {\tt Decen-Scaffnew} & $10^{-4}$ & $130$ & $0.013$ & - & - & -\\
      \midrule
     & {\bfseries Opt. S.S.} & {\bfseries Edge Prob.} & - & - & - & - \\[0.05cm]
    {\tt DSGD} & $10^{-4}$ & $0.013$ & - & - & - & -\\
    \midrule
     & {\bfseries Opt. S.S.} & {\bfseries Quant. Side Length} & - & - & - & - \\[0.05cm]
    {\tt Swarm-SGD} & $5 \times 10^{-5}$ & $10^{-4}$ & - & - & - & -\\
    \midrule
     & {\bfseries Opt. S.S.} & {\bfseries Consensus S.S.} & {\bfseries Active Prob.} & - & - & - \\[0.05cm]
    {\tt CHOCO-SGD} & $10^{-4}$ & $10^{-3}$ & $0.03$ & - & - & -\\
    \bottomrule
    \end{tabular}
    \caption{Hyperparameter values used in Figure \ref{fig:mnist_new}.}
    \label{tab:hypprm_mnist_new}
\end{table}

\begin{table}[h]
    \centering
    \begin{tabular}{l|c c c c}
    \toprule 
    {\algnamesa} & $\max_t \alpha_t$  & $\max_t \eta_t$ & $\gamma$ & $\beta$ \\[0.05cm]
    10\% coordinates & $0.1$ & $5 \times 10^{-9}$ & $0.5$ & $1$ \\
    1\% coordinates & $0.1$ & $10^{-9}$ & $0.5$ & $1$ \\
    0.1\% coordinates & $0.05$ & $5 \times 10^{-10}$ & $0.5$ & $1$ \\
    \midrule
    {\tt CHOCO-SGD} & {\bfseries Max. Opt. S.S.} & {\bfseries Consensus S.S.} & {\bfseries Active Prob.} & - \\[0.05cm]
    10\% coordinates & $0.1$ & $0.05$ & $0.1$ & - \\
    1\% coordinates & $0.1$ & $0.005$ & $0.1$ & - \\
    \midrule
    {\tt Swarm-SGD} & {\bfseries Max. Opt. S.S.} & {\bfseries Quant. Side Length} & - & - \\[0.05cm]
    8-bits quantization & $10^{-3}$ & $3 \times 10^{-5}$ & - & - \\
    \bottomrule
    \end{tabular}
    \caption{Hyperparameter values used in Figure \ref{fig:imagenet}.}
    \label{tab:hypprm_imagenet}
\end{table}

\begin{table}[h]
    \centering
    \begin{tabular}{l|c c c c c c}
    \toprule 
    {\tt DSGD} & {\bfseries Opt. S.S.} & {\bfseries Edge Prob.} & - & -  & - & -\\[0.05cm]
    ($i$) hetero. & $10^{-4}$ & $5 \times 10^{-4}$ & - & -  & - & -\\
    ($ii$) homo. & $10^{-4} $ & $5 \times 10^{-4}$ & - & -  & - & -\\
    \midrule
    {\tt CHOCO-SGD} & {\bfseries Opt. S.S.} & {\bfseries Consensus S.S.} & {\bfseries Active Prob.} & -  & - & -\\[0.05cm]
    ($i$) hetero. & $10^{-4}$ & $10^{-3}$ & $0.1$ & -  & - & -\\
    ($ii$) homo. & $10^{-4}$ & $10^{-3}$ & $0.1$ & -  & - & -\\
    \midrule
    {\tt Swarm-SGD} & {\bfseries Opt. S.S.} & {\bfseries Quant. Side Length} & - & -  & - & -\\[0.05cm]
    ($i$) hetero. & $5 \times 10^{-5}$ & $10^{-4}$ & - & -  & - & -\\
    ($ii$) homo. & $5 \times 10^{-5}$ & $10^{-4}$ & - & -  & - & -\\
    \midrule
    {\algnamesa} & $\alpha$  & $\eta$ & $\gamma$ & $\beta$ & - & - \\[0.05cm]
    ($i$) hetero. & $10^{-4}$ & $10^{-4}$ & $0.5$ & $1$ & - & -\\
    ($ii$) homo. & $10^{-4}$ & $10^{-5}$ & $0.5$ & $1$ & - & - \\
    \midrule
     {\algnamevr} & $\alpha$  & $\eta$ & $\gamma$ & $\beta$  & $a_x$ & $a_\lambda$ \\[0.05cm]
    ($i$) hetero. & $10^{-3}$ & $10^{-3}$ & $0.5$ & $0.1$ & $0.1$ & $0.1$ \\
    ($ii$) homo. & $10^{-3}$ & $10^{-4}$ & $0.5$ & $0.1$  & $0.1$ & $0.1$ \\
    \bottomrule
    \end{tabular}
    \caption{Hyperparameter values used in Figure \ref{fig:mnist}.}
    \label{tab:hypprm_mnist}
\end{table}



\begin{table}[h]
    \centering
    \begin{tabular}{l|c c c c}
    \toprule 
    {\algnamesa} & $\alpha$  & $\eta$ & $\gamma$ & $\beta$ \\[0.05cm]
    no sparse & $10^{-3}$ & $5 \times 10^{-6}$ & 0.5 & 1 \\
    1\% coordinates & $10^{-4}$ & $5 \times 10^{-4}$ & 0.5 & 1 \\
    \midrule
    {\tt DIGing} & {\bfseries Opt. S.S.} & - & - & - \\[0.05cm]
    no sparse & $10^{-3}$ & - & - & - \\
    \bottomrule
    \end{tabular}
    \caption{Hyperparameter values used in Figure \ref{fig:mnist_exact_grad}.}
    \label{tab:hypprm_mnist_exact_grad}
\end{table}

\begin{table}[h]
    \centering
    \begin{tabular}{l|c c c c c c}
    \toprule 
      & $\alpha$  & $\eta$ & $\gamma$ & $\beta$ & $a_x$ & $a_\lambda$ \\[0.05cm]
    {\algnamevr} & $10^{-3}$ & $5 \times 10^{-6}$ & $0.5$ & $1$ & $10^{-3}$ & $1$ \\
    {\algnamevr} & $10^{-3}$ & $5 \times 10^{-6}$ & $0.5$ & $1$ & $10^{-3}$ & $10^{-2}$ \\
    \bottomrule
    \end{tabular}
    \caption{Hyperparameter values used in Figure \ref{fig:mnist_dual_mom}.}
    \label{tab:hypprm_mnist_dual_mom}
\end{table}

\begin{table}[h]
    \centering
    \begin{tabular}{l|c c c c c c}
    \toprule 
      & Time & Peak Memory & Machine  \\[0.05cm]
    \midrule
    Figure \ref{fig:mnist_new} & 27 hours & 339 MB & Intel Xeon Gold 6148
CPU \\
    Figure \ref{fig:imagenet} & 38 hours & \begin{tabular}{@{}c@{}} {(GPU) 112 GB} \\ { (CPU) 6650 MB}\end{tabular} & 8$\times$ NVIDIA V100 GPU \\
    Figure \ref{fig:mnist} & 197 hours & 372 MB & Intel Xeon Gold 6148 CPU\\
    Figure \ref{fig:mnist_sparsity} & 201 hours & 360 MB & Intel Xeon Gold 6148 CPU\\
    Figure \ref{fig:mnist_topology} & 72 hours & 1146 MB & Intel Xeon Gold 6148 CPU\\
    Figure \ref{fig:mnist_exact_grad} & 354 hours & 1066 MB & Intel Xeon Gold 6148 CPU \\
    Figure \ref{fig:mnist_dual_mom} & 46 hours & 350 MB & Intel Xeon Gold 6148 CPU \\
    \bottomrule
    \end{tabular}
    \caption{Statistics of experiment compute time (per algorithm run) and compute instance.}
    \label{tab:compute_stats}
\end{table}

\begin{figure}[h]
    \centering
    \includegraphics[width=0.38\linewidth]{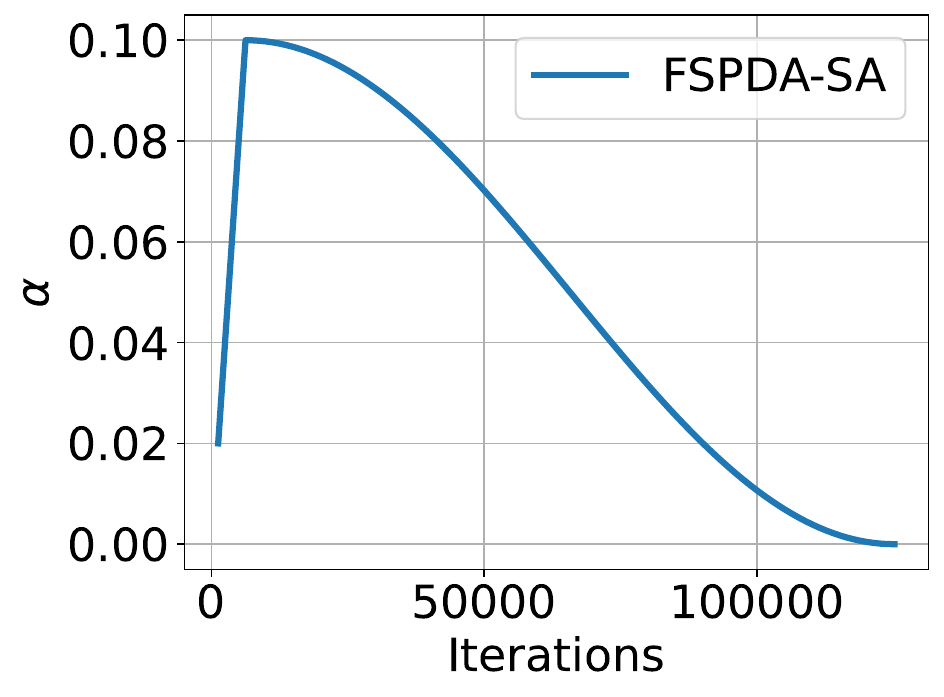}
    \includegraphics[width=0.36\linewidth]{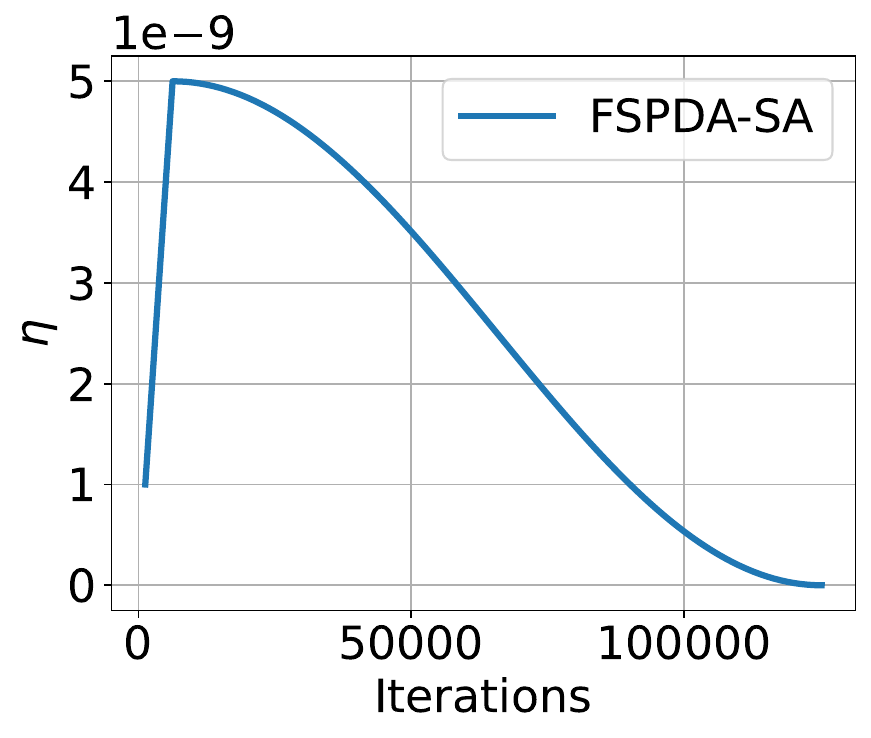}
    \caption{Illustration of step size cosine scheduling used in Figure \ref{fig:imagenet} for {\algnamesa} with 10\% sparse coordinates.}
    \label{fig:ss_scheduling}
\end{figure}

\end{document}